\documentclass[11pt,leqno]{article}
\usepackage{amsmath,amsfonts,amssymb,amsthm}
\newtheorem{definition}{Definition}[section]
\newtheorem{example}[definition]{Example}
\newtheorem{lemma}[definition]{Lemma}
\newtheorem{proposition}[definition]{Proposition}
\newtheorem{theorem}[definition]{Theorem}
\newtheorem{corollary}[definition]{Corollary}
\newtheorem{remark}[definition]{Remark}
\def\proof{\medskip\noindent {\it Proof \hspace{2mm}}}
\def\finpf{\hfill $\square$ \bigskip}
\def\finex{\hfill $\diamond$ \bigskip}

\def\ie{{\it i.e. }}
\def\eg{{\it e.g. }}

\def\N{{\bf N}}
\def\Z{{\bf Z}}
\def\C{{\bf C}}

\def\Q{{\bf Q}}
\def\K{{\bf K}}
\def\P{{\bf P}}
\def\PP{{\bf P}}
\def\F{{\bf F}}
\def\M{{\bf M}}

\def\bra{\langle}
\def\ket{\rangle}
\def\vr{\varnothing}
\def\w{\wedge}

\def\Ia{{\boldsymbol {1}}^+_a}
\def\Ib{{\boldsymbol {1}}^-_b}
\def\kb{{\boldsymbol {k}}}
\def\lb{{\boldsymbol {l}}}

\def\We{{\sf\Lambda}}
\def\aWeb{\We^r(a,b)}

\def\affW{\,\widehat{{\!W\!}}\,}
\def\coxW{\,\widetilde{{\!W\!}}\,}
\def\adW{W_a}
\def\auW{{}^a W}
\def\auaffW{{}^a \affW}
\def\Wbd{W_b}
\def\Wbu{W^b}
\def\affWbu{\affW^b}
\def\He{{\mathcal {H}}}
\def\Gc{{\mathcal {G}}}
\def\Lc{{\mathcal {L}}}
\def\Bc{{\mathcal {B}}}
\def\SIW{\Lambda^{s+\frac{\infty}{2}}}
\def\FN{\F_q[\sal]}
\def\FL{\F_p[\san]}
\def\MN{\M_q[\sal]}
\def\ML{\M_p[\san]}
\def\VN{V_q}

\renewcommand{\leq}{\leqslant}
\renewcommand{\geq}{\geqslant}
\def\lal{{\boldsymbol{\lambda}}_l}
\def\lan{{\boldsymbol{\lambda}}_n}
\def\mul{{\boldsymbol{\mu}}_l}
\def\mun{{\boldsymbol{\mu}}_n}
\def\nul{{\boldsymbol{\nu}}_l}

\def\ol{{\boldsymbol{\varnothing}}_l}
\def\on{{\boldsymbol{\varnothing}}_n}

\def\sal{{\boldsymbol{s}}_l}
\def\san{{\boldsymbol{s}}_n}
\def\tal{{\boldsymbol{t}}_l}

\def\affH{\widehat{H}}
\def\affR{\widehat{{R}}}
\def\gl{ {\mathfrak {gl}}}
\def\sll{ {\mathfrak {sl}}}
\def\asl{\widehat{\sll}}
\def\Ep{E}
\def\Em{\widetilde{E}}
\def\Hp{H}
\def\Hm{\widetilde{H}}

\def\UN{U_q(\asl_n)}
\def\UNp{U'_q(\asl_n)}
\def\UL{U_p(\asl_l)}
\def\ULp{U'_p(\asl_l)}

\def\wt{{\mathrm {wt}}}
\def\res{{\mathrm {res}}}
\def\wtd{\dot{\wt}}
\def\dLambda{\dot{\Lambda}}
\def\dalpha{\dot{\alpha}}
\def\ddelta{\dot{\delta}}
\def\Wt{{\mathrm {Wt}}}

\def\en{e}
\def\fn{f}
\def\tn{t}
\def\dn{\partial}
\def\el{\dot{e}}
\def\fl{\dot{f}}
\def\tl{\dot{t}}
\def\dl{\dot{\partial}}

\def\ep{\varepsilon}
\def\tf{ {\mathfrak {t}} }
\def\hf{ {\mathfrak {h}} }
\def\vd{\dot{v}}
\def\dplus{\Delta^{+}}
\def\dminus{\Delta^{-}}

\def\ov{\overline}
\def\un{\underline}

\begin{document}

\begin{titlepage}

\title{ {\bf Canonical bases of  higher-level $q$-deformed Fock spaces and Kazhdan-Lusztig polynomials}}

\author{Denis Uglov
\thanks{Research Institute for Mathematical Sciences, Kyoto University, 606 Kyoto, Japan.}}
\date{}
%\email{duglov@kurims.kyoto-u.ac.jp}

%\subjclass{Primary 17B67; Secondary 17B37, 17B81}
%
% Abstract
%

\maketitle 
\begin{abstract}
The aim of this paper is to  generalize several aspects of the recent work of Leclerc-Thibon and Varagnolo-Vasserot on the canonical bases of the level $1$ $q$-deformed Fock spaces due to  Hayashi. Namely, we define canonical bases for the higher-level $q$-deformed Fock spaces of Jimbo-Misra-Miwa-Okado and establish a relation between these bases and (parabolic) Kazhdan-Lusztig polynomials for the affine Weyl group of type $A^{(1)}_{r-1}.$ As an application we derive an inversion formula for a sub-family of these polynomials.

\end{abstract}

\thispagestyle{empty}
\end{titlepage}

\section{Introduction}
For any symmetrizable Kac-Moody Lie algebra ${\mathfrak{g}}$ Kashiwara introduces in \cite{Ka1} the notion of a lower global crystal basis of an integrable module $M$ of the universal quantum enveloping algebra $U_{q}({\mathfrak{g}}).$ Furthermore, he proves the existence and uniqueness of such a basis when $M$ is irreducible.    

The lower global crystal bases of irreducible modules were recently recognized to be closely related with the representation theory of Hecke algebras. It was conjectured in \cite{LLT} that the vacuum irreducible module of $\asl_n$ at level $1$ can be identified with the direct sum $\oplus_m K_m$  where $K_m$ is the  complexified Grothendieck group of finitely generated projective modules of the Hecke algebra $H_m(v)$ at $v$ a complex $n$th root of unity. 
Furthermore, it was conjectured that under this identification the specialization of the lower global crystal basis at $q=1$ corresponds to the basis of $K$ formed by the indecomposable summands of the Hecke algebra. 
Proofs of these conjectures, and of a more general result providing a similar interpretation of the lower global crystal basis for any irreducible integrable module of $\asl_n,$ were given in \cite{A}.

If an integrable module of a quantum enveloping algebra is not irreducible, it can have more than one lower global crystal basis. An example  is the $q$-deformed Fock space of Hayashi \cite{H,MM} which is a reducible integrable module of $\UN$ at level $1$ spanned by the set of all partitions. One may then ask whether among the many lower global crystal bases of the $q$-deformed Fock space one can single out a canonical one with particularly favorable properties.   

Leclerc and Thibon gave a definition of such a canonical basis in \cite{LT1}. The essential point of their definition was to introduce a bar-involution of the $q$-deformed Fock space by using the semi-infinite wedge construction of the latter given by Stern \cite{S,KMS}. With the involution in hand, the  canonical basis is defined as the unique bar-invariant basis with a certain congruence property with respect to the $\Q[q]$-lattice spanned by partitions.  

It was conjectured in \cite{LT1} that the specialization at $q=1$ of the transition matrix between the canonical basis of the $q$-deformed Fock space and the basis formed by partitions coincides with the decomposition matrix of the Weyl modules of the $v$-Schur algebra at $v$ a complex $n$th root of unity. This conjecture has been proved in \cite{VV} where it was shown that the entries of the transition matrix are given by certain parabolic Kazhdan-Lusztig polynomials for  affine Weyl groups of type $A^{(1)}_{r-1}.$ An excellent review of these developments can be found in \cite{LT2}.     

The higher-level $q$-deformed Fock spaces generalizing that of Hayashi have been introduced  in \cite{J} in order to compute the crystal graph of an arbitrary irreducible integrable module of $\UN.$ For each sequence of $l$ integers, called charge, there is a  $q$-deformed Fock space which is  an integrable module of $\UN$ at level $l$ having as a basis the set of all $l$-tuples of partitions.  

The main aim of this article is to give a construction of canonical bases for these $q$-deformed Fock spaces generalizing the construction of Leclerc and Thibon for level $1.$ 

As in that case, the first step of this construction is to give a semi-infinite wedge realization of each $q$-deformed Fock space. This already has been done in \cite{TU}, and we follow that paper except in some minor details. The wedge realization allows to define on a $q$-deformed Fock space a natural  bar-involution. 
Then the definition of the canonical basis proceeds exactly as in the level $1$ case. 
In particular, as in \cite{VV}, we find that the entries of the transition matrix between the canonical basis and the basis formed by the $l$-tuples of partitions are parabolic Kazhdan-Lusztig polynomials, the type of the parabolic subgroup being determined by the charge of the $q$-deformed Fock space.      

For a  $q$-deformed Fock space of charge $(s_1,\dots,s_l)$ the $\UN$-submodule $M$ generated by the $l$-tuple of empty partitions is isomorphic to the irreducible module with the highest weight $\sum_{b=1}^l \Lambda_{s_b\bmod n}.$ The canonical basis of the Fock space is, as in the case of level $1,$ a lower global crystal basis and contains the lower global crystal basis of $M$ as a subset.     

The definition of the canonical bases given in this article is constructive. It  provides one with a simple algorithm for computation of the transition matrices between these bases and the bases formed by $l$-tuples of partitions. Due to the main result of \cite{A}, the sub-matrix of this transition matrix that corresponds to the expansion of the lower global crystal basis of $M$ is known to be a $q$-analogue of the decomposition matrix of a certain  Ariki-Koike algebra  at the $n$th root of unity. We get, therefore, an algorithm for computing this decomposition matrix.  

Now let us give an outline of the present article. In Section 2 we give the definitions of the $q$-deformed Fock spaces, of their crystal bases and of the lower global crystal bases of their irreducible submodules. This part is entirely expository, and follows mostly the work \cite{Fo}.  
In Section 3, after giving the necessary background on affine Weyl groups and affine Hecke algebras, we introduce  wedge products and their canonical bases. Also we establish the connection between these canonical bases and parabolic Kazhdan-Lusztig polynomials. Most of the results of this part are straightforward generalizations of results of \cite{KMS} and of a small subset of results in \cite{VV}. 
In the first part of Section 4 we give the realizations of the $q$-deformed Fock spaces as subspaces of the semi-infinite wedge products. Here we follow, with some deviations,  the paper \cite{TU}. In the second part of this section we introduce the bar-involution and define the canonical bases. The content of this  part is very close to the content of \cite{LT1}.         
In Section 5 we describe a symmetry of the bar-involution, and derive, by using this symmetry an inversion formula for certain parabolic Kazhdan-Lusztig polynomials. When the level $l$ equals 1, this formula has already been established in \cite{LT2}.   
Finally, in Section 6 we give some tables of the canonical bases.

Since most of the results of the present article are to be found, in the special case of level $1,$ in  \cite{LT2}, we tried to organize the exposition so that it parallels the relevant parts of that work. 

The preliminary version of this article appeared as the preprint \cite{U}.

\section{The $q$-deformed Fock spaces}
\subsection{Definitions} \label{ss:DEF}
 Let $\hf = (\oplus_{i=0}^{n-1}\Q h_i)\oplus \Q\dn$ be the Cartan subalgebra of $\asl_n,$ and let $\hf^* = (\oplus_{i=0}^{n-1}\Q \Lambda_i)\oplus \Q\delta$ be its dual. Here $\Lambda_0,\dots,\Lambda_{n-1}$ and $\delta$ are, respectively, the fundamental weights and the null root defined in terms of the pairing between $\hf^*$ and $\hf$  by 
$$ \bra \Lambda_i, h_j \ket = \delta_{ij},\quad \bra \Lambda_i,\dn \ket = \bra \delta , h_i \ket = 0, \quad \bra \delta , \dn \ket = 1.$$
The space $\hf^*$ is equipped with a bilinear symmetric form defined by 
$$ \textstyle{(\Lambda_i | \Lambda_j) = \min(i,j) - \frac{ij}{n},\quad (\Lambda_i | \delta) =1,\quad (\delta | \delta) =0.} $$
For $\Lambda \in \hf^*$ we shall write $|\Lambda|^2$ to mean $(\Lambda | \Lambda).$ It will be convenient to extend the index set of the fundamental weights to all integers by setting $\Lambda_i = \Lambda_{i\bmod n}$ for $i \in \Z.$ Then the simple roots are defined for all integer $i$  as $\alpha_i = 2\Lambda_i - \Lambda_{i+1} - \Lambda_{i-1} + \delta_{i\equiv 0\bmod n} \delta,$ where for a statement $S$ we put $\delta_S =1$ if $S$ is true and $\delta_S =0$ if otherwise.   

\bigskip

\noindent Let $\UN$ be the $q$-deformed universal enveloping algebra of $\asl_n.$ This is an algebra over $\K = \Q(q)$ with generators $\en_i, \fn_i , \tn_i, (\tn_i)^{-1}$ $(i=0,\dots,n-1)$ and $\dn.$ The relations between $\en_i, \fn_i , \tn_i, (\tn_i)^{-1}$ are standard (see \eg \cite{Ka}). We define the relations between the degree generator $\dn$ and the rest of the generators by 
$$ [\dn,\en_i] = \delta_{i,0} \en_i,\qquad [\dn,\fn_i] = -\delta_{i,0} \fn_i, \qquad [\dn,\tn_i]=0.  $$
For $l\in \Z,$ a module $M$ of $\UN$ is said to have {\em level} $l$ if the canonical central element $\tn_0\tn_1\cdots \tn_{n-1}$ of $\UN$ acts on $M$ as the multiplication by $q^l.$ Let $\UNp$ be the subalgebra of $\UN$ generated by $\en_i, \fn_i , \tn_i, (\tn_i)^{-1}.$  

\bigskip 

\noindent For a non-negative integer $k,$ let $\Pi_k$ be the set of partitions of $k,$ \ie the set of all non-decreasing sequences of non-negative integers $\lambda = (\lambda_1,\lambda_2,\dots )$ summing to $k.$ Let $\Pi = \sqcup_{k \geq 0} \Pi_k$ be the set of all partitions. For $l \in \N,$ an element $\lal = (\lambda^{(1)},\dots , \lambda^{(l)})$ of $\Pi^l$ is called  {an \em $l$-multipartition}. It will be convenient to identify a multipartition $\lal$ with its diagram defined as the set 
$$ \{ (i,j,b) \in \N^3\:|\: 1 \leq b \leq l,\; 1 \leq j \leq \lambda^{(b)}_i\}.$$
An element of the diagram of a multipartition $\lal$ is called {\em a node} of $\lal,$ and the total number of nodes of $\lal$ is denoted by $|\lal|.$ 

\bigskip 

\noindent Let $\sal = (s_1,\dots,s_l)$ be a sequence of $l$ integers. With any such sequence one associates {\em the $q$-deformed Fock space} $\FN$ defined as 
$$ \FN = \bigoplus_{\lal \in \Pi^l} \K \, |\lal,\sal \ket.$$ 
In other words, $\FN$ is a $\K$-linear space with a distinguished basis $|\lal,\sal \ket$ labelled by the set of all $l$-multipartitions. The number $l$ is called {\em the level} of $\FN,$ and the sequence $\sal$ is called {\em the charge} of $\FN.$  

It was shown in \cite{J} that $\FN$ can be endowed with a structure of  an integrable $\UN$-module. We shall describe this structure following the exposition  given in \cite{Fo}. To do this we introduce some notations. For a node $\gamma = (i,j,b)$ of a multipartition $\lal$ one defines its {\em $n$-residue} as $\res_n(\gamma) = (s_b + j -i)\bmod n.$ For $i$ between $0$ and $n-1,$ we say that $\gamma \in \lal$ is {\em an $i$-node} of $\lal$ if $\res_n(\gamma) = i.$ Given two nodes $\gamma = (i,j,b)$ and $\gamma' = (i',j',b')$ of a multipartition $\lal,$ we write $\gamma < \gamma'$ if either $(s_b + j - i) <  (s_{b'} + j' - i')$ or $(s_b + j - i) =  (s_{b'} + j' - i')$ and $b < b'.$ If $\mul$ and $\lal$ are two multipartitions such that $\mul \supset \lal,$ and $\gamma = \mul\setminus \lal$ is an $i$-node of $\mul,$ we say that $\gamma$ is a removable $i$-node of $\mul,$ and is an addable $i$-node of $\lal.$ In this case we define 
\begin{align*}         
N_i^{>}(\lal,\,\mul |\sal,n) &= \sharp\{ \text{ addable $i$-nodes $\gamma'$ of $\lal$ such that $\gamma' > \gamma$}\}  - \\ 
& - \sharp\{ \text{ removable $i$-nodes $\gamma'$ of $\lal$ such that $\gamma' > \gamma$}\}, \\
N_i^{<}(\lal,\,\mul |\sal,n) &= \sharp\{ \text{ addable $i$-nodes $\gamma'$ of $\lal$ such that $\gamma' < \gamma$}\}  - \\ 
& - \sharp\{ \text{ removable $i$-nodes $\gamma'$ of $\lal$ such that $\gamma' < \gamma$}\}.
\end{align*}
Also, for a multipartition $\lal$ and $i$ between $0$ and $n-1$ we define 
\begin{align*}
&N_i(\lal |\sal,n) = \sharp\{ \text{ addable $i$-nodes of $\lal$}\}   
 - \sharp\{ \text{ removable $i$-nodes of $\lal$} \}, \\
&M_i(\lal |\sal,n) = \sharp\{ \text{  $i$-nodes of $\lal$}\},  
\end{align*}
and for $\sal =(s_1,\dots,s_l) \in \Z^l$ we set 
$$
\Delta(\sal|n) =  \frac{1}{2} \sum_{b=1}^l|\Lambda_{s_b}|^2 +  \frac{1}{2} \sum_{b=1}^l ( \frac{s_b^2}{n} - s_b).
$$
Now we can state 
\begin{theorem}[\cite{J,Fo}]  \label{t:ACT}
The following formulas define on $\FN$ a structure of an integrable $\UN$-module. 
\begin{align*} 
&\fn_i |\lal,\sal\rangle  = \sum_{{\mathrm {res}}_n(\mul/\lal) = i} q^{\; N_i^{>}(\lal,\,\mul |\sal,n)} \: \,\, |\mul,\sal\rangle, \\
&\en_i |\mul,\sal\rangle  = \sum_{{\mathrm {res}}_n(\mul/\lal) = i} q^{- N_i^{<}(\lal,\,\mul |\sal,n)} \:  |\lal,\sal\rangle,  \\
%\end{alignat*}
%\begin{align*} 
&\tn_i |\lal,\sal\rangle  = q^{N_i(\lal |\sal,n)} |\lal,\sal\rangle,  \\ 
& \dn |\lal,\sal\rangle  =  -(\Delta(\sal|n) + M_0(\lal |\sal,n)) |\lal,\sal\rangle .   
\end{align*}
\end{theorem}
\begin{remark} {\rm 
 Our  labeling of the basis vectors differs from that of \cite{Fo} by the transformation  reversing the order of components in $\lal = (\lambda^{(1)},\dots,\lambda^{(l)})$ and $\sal = (s_1,\dots,s_l).$ Also, Theorem \ref{t:ACT}, as well as  Theorem \ref{t:CRYSTAL} below, are stated in \cite{Fo} only for $\sal$ such that $n>s_1 \geq s_2 \geq \cdots \geq s_l \geq 0.$ Generalizations for all $\sal \in \Z^l$ are straightforward.\

}
\end{remark}

\noindent Note that the vector $|\ol,\sal \ket,$ where $\ol$ denotes the $l$-tuple of empty partitions, is a highest weight vector of $\FN.$ Since $\FN$ is an integrable module, it follows that $$\MN = \UN\,|\ol,\sal\ket $$ is an irreducible submodule of $\FN.$ Computing the weight of $|\ol,\sal \ket$ in accordance with Theorem \ref{t:ACT}, we see that $\MN$ is isomorphic to the irreducible $\UN$-module $\VN(\Lambda)$ with highest weight $\Lambda = - \Delta(\sal|n)\delta + \Lambda_{s_1} + \cdots + \Lambda_{s_l}.$     

\subsection{Crystal bases.} \label{ss:CRYSTAL}
The $q$-deformed Fock spaces were introduced in \cite{J} in order to compute the crystal graphs of irreducible integrable modules of $\UN.$ We have seen that any such module is embedded into a $q$-deformed Fock space as the component generated by the highest weight vector labelled by the empty multipartition. From the crystal base theory it follows that the crystal graph of an irreducible module is embedded into the crystal graph of the corresponding Fock space. The last  crystal graph was described in \cite{J}. To recall how the arrows of this graph are obtained, we introduce, following \cite{Fo}, the notion of {\em a good node} of a multipartition $\lal.$   

First, observe that for each $i$ between $0$ and $n-1$ the relation $\gamma < \gamma'$ defines a total order on the set of all $i$-addable and $i$-removable nodes. 
\begin{example} \label{ex:AR}{\rm 
Let $n=3,$ $l=4$ and $\sal = (5,0,2,1).$ Then, marking the  $0$-addable and the $0$-removable nodes of the multipartition 
$$\lal = ((5,3^2,1), (3,2), (4,3,1), (2^3,1))$$ 
on the diagram of $\lal$ by $\bullet,$ we get 
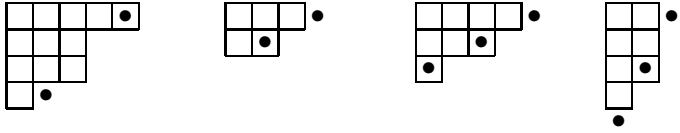
\[
\begin{picture}(50,50)(0,-50)
\setlength{\unitlength}{10pt}
\put(0,0){\line(1,0){5}}
\put(0,-1){\line(1,0){5}}
\put(0,-2){\line(1,0){3}}
\put(0,-3){\line(1,0){3}}
\put(0,-4){\line(1,0){1}}
%%%%%%%%%%%%%%%%%%%%%%%%
\put(0,0){\line(0,-1){4}}
\put(1,0){\line(0,-1){4}}
\put(2,0){\line(0,-1){3}}
\put(3,0){\line(0,-1){3}}
\put(4,0){\line(0,-1){1}}
\put(5,0){\line(0,-1){1}}
%%%%%%%%%%%%%%%%%%%%%%%%
\put(4,-1){\makebox(1,1){$\bullet$}}
\put(1,-4){\makebox(1,1){$\bullet$}}
\end{picture}
\qquad \quad 
%% 2 %%%%%%%%%%%%%%%%%%%%%%%%%%%%%%%%%%%%%%%%%%%%%%%%%%%%%
\begin{picture}(50,50)(0,-50)
\setlength{\unitlength}{10pt}
\put(0,0){\line(1,0){3}}
\put(0,-1){\line(1,0){3}}
\put(0,-2){\line(1,0){2}}
%%%%%%%%%%%%%%%%%%%%%%%%
\put(0,0){\line(0,-1){2}}
\put(1,0){\line(0,-1){2}}
\put(2,0){\line(0,-1){2}}
\put(3,0){\line(0,-1){1}}
%%%%%%%%%%%%%%%%%%%%%%%%
\put(3,-1){\makebox(1,1){$\bullet$}}
\put(1,-2){\makebox(1,1){$\bullet$}}
\end{picture}
\qquad
%% 3 %%%%%%%%%%%%%%%%%%%%%%%%%%%%%%%%%%%%%%%%%%%%%%%%%%%%%
\begin{picture}(50,50)(0,-50)
\setlength{\unitlength}{10pt}
\put(0,0){\line(1,0){4}}
\put(0,-1){\line(1,0){4}}
\put(0,-2){\line(1,0){3}}
\put(0,-3){\line(1,0){1}}
%%%%%%%%%%%%%%%%%%%%%%%%
\put(0,0){\line(0,-1){3}}
\put(1,0){\line(0,-1){3}}
\put(2,0){\line(0,-1){2}}
\put(3,0){\line(0,-1){2}}
\put(4,0){\line(0,-1){1}}
%%%%%%%%%%%%%%%%%%%%%%%%
\put(4,-1){\makebox(1,1){$\bullet$}}
\put(2,-2){\makebox(1,1){$\bullet$}}
\put(0,-3){\makebox(1,1){$\bullet$}}
\end{picture}
\qquad
%% 4 %%%%%%%%%%%%%%%%%%%%%%%%%%%%%%%%%%%%%%%%%%%%%%%%%%%%%
\begin{picture}(50,50)(0,-50)
\setlength{\unitlength}{10pt}
\put(0,0){\line(1,0){2}}
\put(0,-1){\line(1,0){2}}
\put(0,-2){\line(1,0){2}}
\put(0,-3){\line(1,0){2}}
\put(0,-4){\line(1,0){1}}
%%%%%%%%%%%%%%%%%%%%%%%%
\put(0,0){\line(0,-1){4}}
\put(1,0){\line(0,-1){4}}
\put(2,0){\line(0,-1){3}}
%%%%%%%%%%%%%%%%%%%%%%%%
\put(2,-1){\makebox(1,1){$\bullet$}}
\put(1,-3){\makebox(1,1){$\bullet$}}
\put(0,-5){\makebox(1,1){$\bullet$}}
\end{picture}
. \]
 
So, these nodes are ordered as
$$ A_{-3,4} < R_{0,2} < R_{0,3} < R_{0,4} < A_{3,1} < A_{3,2} < R_{3,3} < A_{3,4} < A_{6,3} < R_{9,1}, $$ 
where $A_{d,b}$ ($R_{d,b}$)  denotes an addable (removable) node $(i,j,b)$ with $s_b+j-i = d.$
}
\finex \end{example}

 Next, for a multipartition $\lal$ write the sequence of its addable and removable $i$-nodes ordered as explained above. Then, remove from this sequence recursively all  pairs $R A$ until no such pairs remain. The resulting sequence then has the form $ A \dots A R \dots R.$ The rightmost $R$-node in this sequence is called {\em the good removable $i$-node} of $\lal,$ and the leftmost $A$-node in this sequence is called {\em the good addable $i$-node} of $\lal.$ Clearly there can be at most one of each. For instance, for the multipartition considered in Example \ref{ex:AR} the nodes $A_{-3,4}$ and $R_{9,1}$ are good $0$-nodes.     

Let $A \subset \Q(q)$ be the ring of rational functions without pole at $q=0.$ Let $\Lc[\sal] = \oplus_{\lal \in \Pi^l} A |\lal,\sal \ket,$ and let $\Bc[\sal]$ be the $\Q$-basis of $\Lc[\sal]/q \Lc[\sal]$ given by $\Bc[\sal] = \{ |\lal,\sal\ket \bmod q \Lc[\sal]\:|\: \lal \in \Pi^l\}.$ 

\begin{theorem}[\cite{J,Fo}] \label{t:CRYSTAL}
The pair $(\Lc[\sal],\Bc[\sal])$ is a lower crystal basis of $\FN$ at $q=0.$ Moreover, the crystal graph $\Bc[\sal]$ contains the arrow 
$$ |\lal,\sal \ket\bmod q \Lc[\sal] \overset{i}{\longrightarrow} |\mul,\sal \ket\bmod q \Lc[\sal]  $$
if and only if $\mul$ is obtained from $\lal$ by adding a good $i$-node.
\end{theorem}

\noindent Let $\Pi^l(\sal)$ be the subset of $\Pi^l$ such that $\Bc[\sal]^{\circ}=\{|\lal,\sal\ket \bmod q \Lc[\sal]\:|\: \lal \in \Pi^l(\sal)\}$ is the set of vertices in the connected component of $|\ol,\sal\ket \bmod q\Lc[\sal]$ in the crystal graph $\Bc[\sal]$ of $\FN.$ Then Theorem 3 of \cite{Ka1} implies that $\Bc[\sal]^{\circ}$ is isomorphic to the crystal graph of the irreducible submodule $\MN$ of $\FN.$   

Let us now briefly review the notion of the global crystal base of an irreducible module $\MN.$ First, recall 
the involution $x \mapsto \ov{x}$ of $\UNp$ defined as the unique algebra automorphism satisfying 
$$ \textstyle{\ov{q} = q^{-1},\qquad \ov{\tn_i} = (\tn_i)^{-1}, \qquad \ov{\en_i} = \en_i, \qquad \ov{\fn_i} = \fn_i.}$$
Now, each vector $v$ of $\MN$ can be written as $v = x|\ol,\sal\ket$ for some $x \in \UNp.$ Then we set $\ov{v} = \ov{x}|\ol,\sal\ket.$ Finally, denote by $U^-_{\Q}$ the $\Q[q,q^{-1}]$-subalgebra of $\UNp$ generated by the $q$-divided differences $\fn_i^k/[k]!,$ and let ${{\MN}}_{\Q} = U^-_{\Q}|\ol,\sal\ket.$ (Here, $[k]!$ denotes the $q$-factorial, that is $[k] = {(q^k - q^{-k})}/{(q-q^{-1})}$ and $[k]!  = [k][k-1]\cdots [1].$)

\begin{theorem}[\cite{Ka1}] \label{t:GLOB}
There exists a unique $\Q[q,q^{-1}]$-basis $\{ \Gc(\lal,\sal)\:|\: \lal \in \Pi^l(\sal)\}$ of ${{\MN}}_{\Q}$ such that 
\begin{align*}
&\ov{\Gc(\lal,\sal)} = \Gc(\lal,\sal) \tag*{(i)}, \\ 
&\Gc(\lal,\sal) \equiv |\lal,\sal \ket\bmod q \Lc[\sal]. \tag*{(ii)}
\end{align*}
\end{theorem}
\noindent The basis $\{ \Gc(\lal,\sal)\:|\: \lal \in \Pi^l(\sal)\}$ is called {\em the lower global crystal basis} of $\MN.$

\section{Canonical bases of wedge products}
\subsection{Affine Weyl group}
Let $\tf^* = \oplus_{i=1}^r \C \ep_i$ be the dual space of the Cartan subalgebra of 
$\gl_r.$ 
Let $\hat{\tf}^* = \tf^*\oplus\C \Lambda_0 \oplus \C \delta$ be the dual space of the  
Cartan subalgebra of $\widehat{\gl}_r.$ The space $\hat{\tf}^*$ is equipped with the  
bilinear symmetric form defined by 
$(\ep_i,\ep_j) = \delta_{ij},$ $(\ep_i,\Lambda_0)=(\ep_i,\delta)=(\delta,\delta)=(\Lambda_0,\Lambda_0)=0,$ $(\Lambda_0,\delta)=1.$
The systems of roots $R,$ of positive roots $R^+$ and simple roots $\Pi$ of type $A_{r-1}$ are the subsets of $\tf^*$ defined by   
\begin{alignat*}{3}
&R   &&=&& \:\{ \alpha_{ij}=\ep_i-\ep_j \:|\: i\neq j\},\\
&R^+ &&=&& \:\{ \alpha_{ij} \:|\: i < j\},\\
&\Pi &&=&& \:\{\alpha_1,\dots,\alpha_{r-1}\}, \quad (\alpha_i := \alpha_{ii+1}).
\end{alignat*}
The systems of roots $\affR,$ of positive roots $\affR^+$ 
and simple roots $\widehat{\Pi}$ of type $A_{r-1}^{(1)}$ are the subsets of $\hat{\tf}^*$ defined by
\begin{alignat*}{3}
&\affR   &&=&& \: \{\alpha + k\delta \:|\: \alpha \in R,\:k\in \Z\},\\
&\affR^+ &&=&& \: \{\alpha + k\delta \:|\: \alpha \in R^+,\:k \geq 0\}\sqcup\{ -\alpha + k\delta \:|\: \alpha \in R^+,\:k > 0\},\\ 
&\widehat{\Pi} &&=&& \: \{\alpha_0:= \delta -(\ep_1-\ep_r)\}\sqcup\Pi.
\end{alignat*} 
The Weyl group $W$ of $\gl_r$ is isomorphic to the symmetric group ${\mathfrak {S}}_r,$ and has a realization as the group generated by the reflections $s_{\alpha}(\xi) = \xi - (\alpha,\xi)\,\alpha,$  $(\alpha \in R)$  of  $\tf^*.$ Let $Q = \oplus_{i=1}^{r-1}\Z \alpha_i$ and $P=\oplus_{i=1}^{r}\Z \ep_i$ be, respectively, the root and the weight lattices of $\gl_r.$ They both are preserved by $W.$ 

The affine Weyl group is defined as the semi-direct product
\[
\affW = W \ltimes P
\]
with relations $ w t_{\eta} = t_{w(\eta)} w,$ where $w$ and $t_{\eta}$ are elements of $\affW$ that correspond to $w \in W,$ $\eta \in P.$ The group $\affW$ contains the Weyl group $\coxW = W \ltimes Q$ of type $A_{r-1}^{(1)}$ as a subgroup. The group $\affW$ acts on $\hat{\tf}^*$ by 
\begin{alignat*}{4}
&s_{\alpha}(\zeta) & &=   & & \zeta - (\alpha,\zeta)\,\alpha, & & \quad (\zeta \in \hat{\tf}^*,\: \alpha \in R),\\
&t_{\eta}(\zeta) & &= & &\zeta + (\delta,\zeta)\,\eta - ( (\eta,\zeta) +  {\textstyle\frac{1}{2}} (\eta,\eta)(\delta,\zeta))\,\delta, & &\quad (\zeta \in \hat{\tf}^*,\: \eta \in P).
\end{alignat*} 
This action preserves the root system $\affR,$ and the bilinear form on $\hat{\tf}^*$ is invariant with respect to this action. 

For an affine root $\hat{\alpha} = \alpha + k\delta$ $(\alpha \in R, k\in \Z),$ define the corresponding affine reflection as $s_{\hat{\alpha}} = t_{-k\alpha}\,s_{\alpha},$ and put $s_i = s_{\alpha_i}$ $(i=0,1,\dots,r-1),$ $\pi = t_{\ep_1}s_1\cdots s_{r-1}.$ The group $\affW$ is generated by $\pi,\pi^{-1},\,s_0,s_1,\dots,s_{r-1}$ and is defined by the relations 
\begin{align*}
&s_i s_{i+1} s_i = s_{i+1} s_i s_{i+1}, \\
&s_i s_j = s_j s_i \quad (i-j \neq \pm 1), \\
&s_i^2 = 1,\qquad \pi s_i = s_{i+1} \pi, 
\end{align*}
where the subscripts are understood to be modulo $r.$ In this presentation  $\coxW$ is the Coxeter subgroup generated by $s_0,s_1,\dots,s_{r-1},$ and $\affW \cong \Omega \ltimes \coxW$ where $\Omega \cong \Z$ is the subgroup of $\affW$ generated by  $\pi,\pi^{-1}.$   

For $w\in \affW$ let $S(w) = \affR^+\cap w^{-1}(\affR^-)$ where $\affR^-=\affR\setminus\affR^+$ is the set of negative roots. The length $l(w)$ of $w$ is defined as the number $\sharp S(w)$ of elements in $S(w).$  The length of $w$ is zero  if and only if $w \in \Omega.$ A partial order on $\affW$ is defined by $\pi^k w \preceq \pi^{k'}w'$ $(w,w' \in \coxW)$ if $k=k',$ and  $w \preceq w'$ in the Bruhat order of $\coxW.$

The following lemma follows immediately from the definition of $S(w).$ 
\begin{lemma} \label{l:LENGTH} \mbox{} \\ {\em (i)} For $w \in \affW,$ $S(w^{-1}) = -w(S(w)).$ \\ 
{\em (ii)} For $u,v \in \coxW,$ $S(u) = S(v)$ implies $u=v.$  \\
{\em (iii)} For $w\in W,$ $\lambda \in P$ 
\[
l(w\, t_{\lambda}) =\sum_{  \alpha \in R^+,\: w(\alpha) \in R^+} |(\lambda,\alpha)| + \sum_{  \alpha \in R^+,\:w(\alpha) \in R^-} |1+ (\lambda,\alpha)|. 
\]
\end{lemma}
\noindent A corollary to (iii) above is the equality 
$ l(t_{\lambda}) = l(t_{\mu}) $ for $\lambda, \mu \in P$ such that $\lambda = w(\mu)$ for some $w\in W.$

The following lemma is contained in \cite{AST} as Definition and Proposition 2.2.2.
\begin{lemma} \label{l:SHORTEST}
For $\lambda \in P,$ let $w$ be the shortest element of $W$ such that $w(\lambda) \in P^+.$ Then  
\[
S(w) = \{ \alpha \in R^+\:|\: (\lambda,\alpha) < 0 \}.
\]
\end{lemma}

\begin{proposition} \label{p:SPLIT}
For every  $x\in \affW$ there is a unique factorization of the form  $x= u\,t_{\lambda} v,$ where $u,v \in W,$ $\lambda \in P^+,$ and $S(v) = \{ \alpha \in R^+\:|\: (\lambda,v(\alpha)) < 0 \}.$ Moreover, $l(x) = l(u) + l(t_{\lambda}) - l(v).$
\end{proposition}

\proof Every $x \in \affW$ can be factorized as $x= w\,t_{\mu},$ where $w\in W,$ $\mu \in P.$ Let $v\in W$ be the shortest element such that $v(\mu)\in P^+.$ By Lemma \ref{l:SHORTEST}, $S(v)=\{ \alpha \in R^+\:|\: (\mu,\alpha) < 0 \}.$ The desired factorization is afforded by $x= u\,t_{\lambda} v,$ with  $u=wv^{-1},$ $\lambda = v(\mu).$ 

Assume $x= u_1\,t_{\lambda_1} v_1 = u_2\,t_{\lambda_2} v_2 ,$ where  
$u_i,v_i,\lambda_i$ satisfy the conditions listed in the statement of the proposition. 
%$u_i,v_i\in W,$ $\lambda_i \in P^+,$ and $S(v_i) = \{ \alpha \in R^+\:|\: (\lambda_i,v_i(\alpha)) < 0 \}.$ 
Put $\mu_i = v_i^{-1}(\lambda_i),$ so that $x= u_i v_i\, t_{\mu_i}.$ The presentation of $x$ in the form $w\,t_{\mu},$ $(w\in W,\,\mu \in P)$ is unique, hence $\mu_1=\mu_2,$ $u_1 v_1 = u_2 v_2.$ The equality $\mu_1=\mu_2$ implies $S(v_1) = S(v_2),$ whence $v_1 = v_2,$ and, therefore, $u_1=u_2.$  The factorization is unique.

It remains to show  the relation $l(x) = l(u) + l(t_{\lambda}) - l(v).$ The length formula of Lemma \ref{l:LENGTH} together with  $S(v)=\{ \alpha \in R^+\:|\: (\mu,\alpha) < 0 \}$  give
$$
l(x)=l(w\,t_{\mu})= 
      \sum_{\alpha \in R^+\setminus S(v),\: w(\alpha) \in R^- } 1 - 
      \sum_{\alpha \in S(v),\: w(\alpha) \in R^- } 1  + 
      \sum_{\alpha \in R^+\setminus S(v)}(\mu,\alpha) - 
      \sum_{\alpha \in S(v)}(\mu,\alpha).
$$
On the other hand 
\begin{align*}
&l(v) = \sum_{\alpha \in S(v)} 1 = \sum_{\alpha \in S(v),\: w(\alpha) \in R^+ } 1 + \sum_{\alpha \in S(v),\: w(\alpha) \in R^- } 1, \\
&l(u) = l(w v^{-1}) = \sum_{\alpha \in R^+\setminus S(v),\: w(\alpha) \in R^- } 1 + \sum_{\alpha \in S(v),\: w(\alpha) \in R^+ } 1, \\ 
&l(t_{\lambda}) = l(t_{\mu}) =  \sum_{\alpha \in R^+\setminus S(v)}(\mu,\alpha) - \sum_{\alpha \in S(v)}(\mu,\alpha).
\end{align*}
The relation $l(x) = l(u) + l(t_{\lambda}) - l(v)$ follows. \finpf

\subsubsection{A right action of $\affW$ on $P.$}
Let $n$ be a positive integer, and define a right action of $\affW$ on $P$ by 
\begin{equation}\label{e:RA}
\begin{aligned}
 \mbox{}&\zeta\!\cdot s_i & &= s_i(\zeta), & &\qquad (\zeta \in P,\; 1\leq i < r), \\
 \mbox{}&\zeta\!\cdot t_{\mu}&  &= \zeta + n \mu, & & \qquad (\zeta \in P,\; \mu \in P).
\end{aligned}
\end{equation}
In coordinates $(\zeta_1,\dots,\zeta_r)$ of $\zeta = \sum_{i=1}^r \zeta_i \ep_i
$ this action looks as follows 
\begin{alignat*}{3}
&(\zeta_1,\dots,\zeta_r) \cdot s_i & & = (\dots,\zeta_{i+1},\zeta_i,\dots), & &\qquad (1\leq i < r), \\
&(\zeta_1,\dots,\zeta_r)\cdot t_{\ep_i}& & = (\dots, \zeta_i + n , \dots ). & & \qquad (1\leq i \leq  r). 
\end{alignat*}
Hence 
\begin{alignat*}{2}
&(\zeta_1,\dots,\zeta_r) \cdot \pi & & = (\zeta_2,\dots,\zeta_r,\zeta_1+n),\\
&(\zeta_1,\dots,\zeta_r)\cdot s_0 & & =  (\zeta_r-n,\zeta_2,\dots,\zeta_{r-1},\zeta_1 + n).
\end{alignat*}
Define  $A^n \subset P$ by 
$$ 
A^n = \{ a=(a_1,\dots,a_r) \in P\:|\: 1 \leq a_1 \leq a_2 \leq \cdots \leq a_r \leq n\}.
$$ 
Then $A^n$ is a fundamental domain of the action given by (\ref{e:RA}).
For $a \in A^n,$ denote by $\adW$ the stabilizer of $a.$ The inequality $a_r-a_1 < n$ implies that $\adW \subset W.$ Let $\auaffW$ (resp.$\auW$)  be the set of minimal length representatives in the cosets $\adW\setminus\affW$ (resp.$\adW\setminus W$).  

\begin{lemma} \label{l:WaWa} Let $x = u\,t_{\lambda} v$ be the factorized presentation of $x \in \affW$ given by Proposition \ref{p:SPLIT}. Then $x\in \auaffW$ if and only if $u\in \auW.$ 
\end{lemma}
\proof For any $w \in W_a$ the factorized presentation of $wx$ is $ wx = (w u)\,t_{\lambda} v.$ It now follows from the length relation of  Proposition \ref{p:SPLIT} that $l(wx) \geq l(x) \Longleftrightarrow l(wu) \geq l(u),$ i.e. $x$ is the shortest element of its coset if and only if $u$ is the shortest element of {\em its} coset.  \finpf

\begin{lemma}[\cite{AST}] \label{l:R_a}
For $a \in A^n,$ let $R^+_a = \{\ep_i - \ep_j \in R^+ \: | \: a_i = a_j\}.$ Then  
$$\auW = \{ u \in W\:|\: S(u^{-1}) \subset R^+\setminus R^+_a \}.$$ 
\end{lemma}

\noindent For $w \in W$ let  a map $w : \{1,\dots,r\} \rightarrow  \{1,\dots,r\}$ be defined by $ \ep_{w^{-1}(i)} = w(\ep_i).$ Note that for $u,v \in W,$ $u(v(i)) = vu(i).$  

\begin{lemma} \label{l:S(u)}
Let $u \in \auW,$ and let $c=(c_1,\dots,c_r) = a\cdot u.$ Then 
$$  S(u) = \{ \ep_i - \ep_j \in R^+ \: | \: c_i > c_j\}. $$ 
\end{lemma}
\proof 
Observe that $c_i = a_{u^{-1}(i)}$ for all $i=1,\dots,r.$ Also, $\ep_i - \ep_j \in S(u)$ if and only if $i<j, u^{-1}(i) > u^{-1}(j).$    
 Since the sequence $a$ is non-decreasing, $i < j, a_{u^{-1}(i)} > a_{u^{-1}(j)}$ implies $i < j, {u^{-1}(i)} > {u^{-1}(j)},$ \ie $\ep_i - \ep_j \in S(u).$
Conversely, $\ep_i - \ep_j \in S(u)$ implies, by Lemma \ref{l:LENGTH} (i) and Lemma \ref{l:R_a},  that $u( \ep_j - \ep_i ) = \ep_{u^{-1}(j)}-\ep_{u^{-1}(i)}   \in R^+\setminus R^+_a.$ This gives $a_{u^{-1}(i)} > a_{u^{-1}(j)},$ and the  lemma follows. \finpf

\begin{proposition} \label{p:LENGTHR}
For $a\in A^n$ and $x \in \auaffW,$ let $ a\cdot x = (c_1 + n\mu_1,\dots,c_r + n\mu_r),$ where $c_i \in \{1,\dots,n\}$ and $\mu_i \in \Z.$  Then   
\begin{gather*} 
l(x) = \sharp\{i < j\:| \: c_i > c_j,\, \mu_i \geq \mu_j\} + \sharp\{i<j\: | \: c_i < c_j,\, \mu_i < \mu_j\} + \\ \sum_{i<j,\, \mu_i > \mu_j } (\mu_i - \mu_j) + \sum_{i<j,\, \mu_i < \mu_j } (\mu_j - \mu_i -1). 
\end{gather*} 
\end{proposition}
\proof Let $ x= u\,t_{\lambda} v$  be the factorized presentation of Proposition \ref{p:SPLIT}. The expression for $l(x)$ follows from $l(x) = l(u) + l(t_{\lambda}) - l(v),$ Lemma \ref{l:WaWa}, Lemma \ref{l:S(u)} and the length formula of Lemma \ref{l:LENGTH}(iii) \finpf

\begin{proposition} \label{p:RIGHTACTIONofW}
For $a\in A^n,$ let $x\in \auaffW,$ and let $\zeta = (\zeta_1,\dots,\zeta_r) = a\cdot x.$ Put $\zeta_0 = \zeta_r - n.$ Then for each $i=0,1,\dots,r-1$ one has the following complete set of alternatives: 
\begin{alignat}{3}
&\zeta_i = \zeta_{i+1} &  &\Longleftrightarrow  xs_i \not\in  \auaffW, & &\quad  \tag*{(i)}\\ 
&\zeta_i > \zeta_{i+1} &  &\Longleftrightarrow  xs_i \in  \auaffW, & &\quad l(xs_i) = l(x) -1, \tag*{(ii)}\\
&\zeta_i < \zeta_{i+1} &  &\Longleftrightarrow  xs_i \in  \auaffW, & &\quad l(xs_i) = l(x) + 1. \tag*{(iii)}
\end{alignat}
Moreover, in the case {\em (i)}, $ xs_i = s_jx $ where $s_j \in \adW.$ 
\end{proposition}
\proof  First, we show (i). Let $\zeta_i = \zeta_{i+1},$ then $\zeta\!\cdot s_i = \zeta.$ Assuming  $xs_i \in \auaffW$ the length formula of  Proposition \ref{p:LENGTHR} is applicable, and immediately  gives  $l(xs_i) = l(x),$ which is impossible. Hence, $\zeta_i = \zeta_{i+1}$ $\Longrightarrow$ $xs_i \not\in \auaffW.$ Now let $xs_i \not\in \auaffW.$ By \cite[Lemma 2.1 (iii)]{Deodhar}, in this case  $xs_i = s_jx$ for some $s_j \in \adW,$ which implies $\zeta\cdot s_i = \zeta,$ hence $\zeta_i = \zeta_{i+1}.$ Thus $\zeta_i = \zeta_{i+1}$ $\Longleftarrow$ $xs_i \not\in \auaffW.$ 

Let $\zeta_i > \zeta_{i+1}$ (resp. $\zeta_i < \zeta_{i+1}$). Then $xs_i \in  \auaffW,$ and one may use the length formula of Proposition \ref{p:LENGTHR} to show $l(xs_i) = l(x)-1$ (resp.  $l(xs_i) = l(x) +1$). Since (i) is already established, this proves  (ii) and (iii). 

Finally, \cite[Lemma 2.1 (iii)]{Deodhar} states that in the case (i), $ xs_i = s_jx $ where $s_j \in \adW.$    \finpf  

%\begin{remark} {\rm 
%The {\em left} action of $\affW$ on $P$ obtained from the right action (\ref{e%:RA}) by means of the antiautomorphism $w \mapsto w^{-1}$ is identical with the action denoted $\pi_{-n}$ in \cite{LT2}.  
%}
%\end{remark}

\subsubsection{A left  action of $\affW$ on $P.$}
Let $l$ be a positive integer, and define a left action of $\affW$ on $P$ by 
\begin{equation} \label{e:LA}
\begin{aligned}
\mbox{}& s_i\cdot \eta  &  &= s_i(\eta), & &\qquad (\eta \in P,\; 1\leq i < r), \\
\mbox{}&t_{\mu}\cdot \eta &  &= \eta + l \mu, & &\qquad (\eta \in P,\; \mu \in P).
\end{aligned}
\end{equation}
In coordinates $(\eta_1,\dots,\eta_r)$ of $\eta = \sum_{i=1}^r \eta_i \ep_i
$ this action looks as follows 
\begin{alignat*}{3}
&s_i\cdot (\eta_1,\dots,\eta_r)  & & = (\dots,\eta_{i+1},\eta_i,\dots), & &\qquad (1\leq i < r), \\
& t_{\ep_i}\cdot (\eta_1,\dots,\eta_r)& & = (\dots, \eta_i + l , \dots ). & & \qquad (1\leq i \leq  r). 
\end{alignat*}
Hence 
\begin{alignat*}{2}
&\pi\cdot (\eta_1,\dots,\eta_r)  & & = (\eta_r+l,\eta_1,\dots,\eta_{r-1}),\\
&s_0\cdot (\eta_1,\dots,\eta_r) & & =  (\eta_r+l,\eta_2,\dots,\eta_{r-1},\eta_1 -l).
\end{alignat*}
Define  $B^l \subset P$ by 
$$ 
B^l = \{ b=(b_1,\dots,b_r) \in P\:|\: l \geq b_1 \geq b_2 \geq \cdots \geq b_r \geq 1\}.
$$ 
Then $B^l$ is a fundamental domain of the action given by (\ref{e:LA}).
For $b \in B^l,$ denote by $\Wbd$ the stabilizer of $b.$ The inequality $b_1-b_r < l$ implies that $\Wbd \subset W.$ Let $\affWbu$ (resp. $\Wbu$) be the set of minimal length representatives in the cosets $\affW/\Wbd$ (resp.$ W/\Wbd$).

\begin{lemma}[\cite{AST}] \label{l:R_b}
For $b \in B^l,$ let $R^+_b = \{\ep_i - \ep_j \in R^+ \: | \: b_i = b_j\}.$ Then 
$$\Wbu = \{ v \in W\:|\: S(v) \subset R^+\setminus R^+_b \}.$$ 
\end{lemma}

%\begin{lemma} \label{l:vi<vj}
%Let $v \in \Wbu,$ $i<j,$ and $b_i=b_j.$ Then $v^{-1}(i) < v^{-1}(j).$ 
%\end{lemma}
%\proof $S(v) \subset R^+\setminus R^+_b$ implies $\ep_i - \ep_j \not\in S(u),$ \ie $v(\ep_i-\ep_j) = \ep_{v^{-1}(i)} - \ep_{v^{-1}(j)} \in R^+.$ \finpf

\begin{lemma} \label{l:S(v^{-1})}
Let $v \in \Wbu,$ and let $d=(d_1,\dots,d_r) = v\cdot b.$ Then 
$$  S(v^{-1}) = \{ \ep_i - \ep_j \in R^+\: | \: d_i < d_j\}. $$ 
\end{lemma}

\begin{proposition} \label{p:LENGTHL}
For $b\in B^l$ and  $x \in \affWbu,$ let $  x\cdot b = (d_1 + l\mu_1,\dots,d_r + l\mu_r),$ where $d_i \in \{1,\dots,l\}$ and $\mu_i \in \Z.$  Then   
\begin{gather*} 
l(x) = \sharp\{i < j\:| \: d_i < d_j,\, \mu_i \leq \mu_j\} + \sharp\{i<j\: | \: d_i > d_j,\, \mu_i > \mu_j\} + \\ \sum_{i<j,\, \mu_i > \mu_j } (\mu_i - \mu_j -1) + \sum_{i<j,\, \mu_i < \mu_j } (\mu_j - \mu_i). 
\end{gather*} 
\end{proposition}

\begin{proposition} \label{p:LEFTACTIONofW}
For $b\in B^l,$ let $x\in \affWbu,$ and let $\eta = (\eta_1,\dots,\eta_r) =  x\cdot b.$ Put $\eta_0 = \eta_r +l.$ Then for each $i=0,1,\dots,r-1$ one has the following complete set of alternatives: 
\begin{alignat}{3}
&\eta_i = \eta_{i+1} &  &\Longleftrightarrow  s_ix \not\in  \affWbu, & &\quad  \tag*{(i)}\\ 
&\eta_i > \eta_{i+1} &  &\Longleftrightarrow  s_ix \in  \affWbu, & &\quad l(s_ix) = l(x) +1, \tag*{(ii)}\\
&\eta_i < \eta_{i+1} &  &\Longleftrightarrow  s_ix \in  \affWbu, & &\quad l(s_ix) = l(x) - 1. \tag*{(iii)}
\end{alignat}
Moreover, in the case {\em (i)}, $ s_ix = xs_j $ where $s_j \in \Wbd.$ 
\end{proposition}

\noindent We omit proofs of Lemma \ref{l:S(v^{-1})} and of the last two propositions because they  are almost identical to the proofs of Lemma \ref{l:S(u)} and of Propositions \ref{p:LENGTHR}, \ref{p:RIGHTACTIONofW}. 

%\begin{remark} {\rm
%In \cite{LT2}, the left action of $\affW$ defined by (\ref{e:LA}) is  denoted  $\pi_l.$ 
%
%}
%\end{remark}

\subsection{Affine Hecke algebra}

The Hecke algebra $\affH$ of the Weyl group $\affW$ is the algebra over $\K =\Q(q)$ with basis $ T_x $ $(x \in \affW)$ and relations  
\begin{gather*}
T_x T_y  = T_{xy} \; \; \text{ whenever $\;l(xy) = l(x) + l(y),$} \\
( T_{s_i} - q^{-1}) (T_{s_i} + q) = 0 \; \; \text{ for all $\;i=0,1,\dots,r-1.$}
\end{gather*}
The subalgebra $H$ of $\affH$ generated by $T_w$ $(w \in W)$ is isomorphic to the Hecke algebra of the finite Weyl group $W.$  A system of generators of $\affH$ is afforded by elements $T_{\pi}, T_{\pi^{-1}}$ and $T_0,T_1,\dots,T_{r-1}$ where, for simplicity, we put $T_i := T_{s_i}.$ 

Another system of generators is obtained as follows. For $\lambda \in P,$ write  $\lambda = \mu - \nu $ where $\mu,\nu \in P^+,$ and define 
$$X^{\lambda} = T_{t_{\nu}}T_{t_{\mu}}^{-1}.$$ 
Note that Lemma \ref{l:LENGTH} (iii) implies  $l(t_{\mu} t_{\nu}) =  l(t_{\mu}) +   l(t_{\nu})$ for $\mu,\nu \in P^+.$ From this it follows that $X^{\lambda}$ does not depend on the choice of $\mu,\nu \in P^+,$ and  
$$  X^{\lambda} X^{\mu} = X^{\lambda + \mu}\;\;\text{ for $\lambda,\mu \in P.$} $$ 
A proof of the following result is contained, \eg in \cite{Ki}.
\begin{lemma} \label{l:Ki}
{\em (i)} The elements $X^{\lambda}$,$\lambda \in P,$ $T_1,\dots,T_{r-1}$ generate $\affH.$ \\
{\em (ii)} For $\lambda \in P,$ $i=1,\dots,r-1,$  
\begin{align*}
X^{\lambda} T_i = T_i X^{s_i(\lambda)} + (q-q^{-1}) \frac{ X^{s_i(\lambda)} - X^{\lambda}}{1 - X^{\alpha_i}}, \\ 
T_i X^{\lambda} = X^{s_i(\lambda)} T_i + (q-q^{-1}) \frac{ X^{s_i(\lambda)} - X^{\lambda}}{1 - X^{\alpha_i}}.
\end{align*}
\end{lemma}

Following \cite{So} let us now briefly recall the notions of the canonical bases and the Kazhdan-Lusztig polynomials of $\affW.$ 

First, we recall that there is a canonical involution $h \mapsto \ov{h}$ of $\affH$ defined as the unique algebra automorphism such that $\ov{T_x} = (T_{x^{-1}})^{-1}$ and $\ov{q} = q^{-1}.$ A proof of the following lemma is straightforward.  
\begin{lemma}   \label{l:INVO1}
For $u,v \in W$ and $\lambda \in P,$
$$ \ov{ T_u\, X^{\lambda} \,(T_{v^{-1}})^{-1} } = T_{u\omega}\, X^{\omega(\lambda)} \,(T_{(\omega v)^{-1}})^{-1} ,$$ 
where $\omega$ is the longest element of $W.$ 

\end{lemma}

Let $L^+$ (resp. $L^-$) be the lattice spanned over $\Z[q]$ (resp. $\Z[q^{-1}]$) by $T_x$ $(x \in \affW).$ The canonical bases $C_x',C_x$ $(x\in\affW)$ are the unique bases of $\affH$ with the properties   
\begin{align*}
\ov{C_x'} = {\textstyle{C_x'}}, & \quad \ov{C_x} = C_x, \\ 
C_x' \equiv T_x \bmod qL^+, & \quad C_x \equiv T_x \bmod q^{-1}L^-.
\end{align*}
Let 
$$ C_x' = \sum_y {\mathcal{P}}^+_{y,x} T_y,\qquad C_x = \sum_y {\mathcal{P}}^-_{y,x} T_y.$$  
The coefficients ${\mathcal{P}}^{\pm}_{y,x}$ are called the Kazhdan-Lusztig polynomials of $\affW,$ they are non-zero only if $y \preceq x,$ that is, only if $x=\pi^k \tilde{x}, y=\pi^k \tilde{y}$ for some $k\in\Z,\tilde{x},\tilde{y} \in \coxW$ such that $\tilde{y} \preceq \tilde{x}.$ 
In this case 
 $${\mathcal{P}}^+_{y,x} = q^{l(x)-l(y)} P_{\tilde{y},\tilde{x}},\qquad {\mathcal{P}}^-_{y,x} = (-q)^{l(y)-l(x)} \ov{P_{\tilde{y},\tilde{x}}},$$ 
where $P_{\tilde{y},\tilde{x}} \in \Z_{\geq 0}[q^{-2}]$ are the Kazhdan-Lusztig polynomials of the Coxeter group $\coxW$.  

\subsubsection{A right representation of $\affH.$}
For $a \in A^n,$ let $H_a$ be the parabolic subalgebra of $\affH$ generated by $T_w$ $(w \in W_a).$ Let $\K\Ia$ be a one-dimensional right representation of $H_a$ defined by  
$$ \Ia\cdot T_i = q^{-1}\,\Ia \qquad (s_i \in W_a). $$  
The induced right representation 
$$  \K\Ia\otimes_{H_a}\affH $$
of $\affH$ has as its basis $\Ia\otimes_{H_a} T_x$ $(x \in \auaffW).$ For $x\in \auaffW,$ define 
$$ (\zeta | := \Ia\otimes_{H_a} T_x, \qquad \text{where $\quad \zeta = a\cdot x,$ $\zeta \in P.$} $$  
Then $(\zeta |$ $( \zeta \in a\cdot \affW)$ is a basis of $  \K\Ia\otimes_{H_a}\affH, $ and $(\zeta |$ $( \zeta \in P)$ is a basis of $  \oplus_{a \in A^n} \K\Ia\otimes_{H_a}\affH. $ Proposition \ref{p:RIGHTACTIONofW} allows to describe the action of the affine Hecke algebra  in the basis $(\zeta |$ explicitly. We have    

\begin{align}
 & (\zeta |\cdot T_i = \begin{cases}  \;( \zeta\cdot s_i |  & \quad \text{if $\zeta_i < \zeta_{i+1},$} \\   
                                   \;q^{-1} ( \zeta |  & \quad \text{if $\zeta_i = \zeta_{i+1},$} \\ 
                                 \;( \zeta\cdot s_i |  - (q-q^{-1}) ( \zeta | & \quad \text{if $\zeta_i > \zeta_{i+1}$}
\end{cases}\qquad (0\leq i < r), \label{e:RAH}\\ 
& (\zeta |\cdot T_{\pi} = (\zeta\cdot \pi|.  \nonumber
\end{align}

Define a canonical involution $ v \mapsto \ov{v}$ of  $\K\Ia\otimes_{H_a}\affH $ by 
$$ \ov{\Ia\otimes_{H_a} h } = \textstyle{\Ia\otimes_{H_a} \ov{h }} \qquad (h \in \affH),$$ 
and two lattices by 
\[
 L_a^+ := \bigoplus_{ \zeta \in a\cdot\affW} \Z[q] \, (\zeta |,\qquad L_a^- := \bigoplus_{ \zeta \in a\cdot\affW} \Z[q^{-1}] \, (\zeta |.  
\]

\begin{theorem}[\cite{Deodhar}]
There are unique bases $C^{\pm}_{\zeta}$ $(\zeta \in a\cdot \affW)$ of $\K\Ia\otimes_{H_a}\affH, $ such that 
\begin{align}
& \ov{C_{\zeta}^{\pm} } = C_{\zeta}^{\pm}, \tag*{(i)}\\ 
& C_{\zeta}^{\pm}  \equiv (\zeta | \bmod q^{\pm 1} L^{\pm}_a. \tag*{(ii)}
\end{align}
Moreover, if $$ C_{\zeta}^{\pm} = \sum_{\eta} P^{\pm}_{\eta,\zeta}\,(\eta\, |,$$
then 
$$ P^+_{\eta,\zeta} = {\mathcal {P}}^+_{\omega_a y,\omega_a x},\qquad P^-_{\eta,\zeta} = \sum_{u \in W_a} q^{-l(u)}\: {\mathcal {P}}^-_{u y, x},$$
where $x,y$ are unique elements of $\auaffW$ such that $a\cdot x = \zeta,$ $a\cdot y = \eta,$ and $\omega_a$ is the longest element of $W_a.$ 
\end{theorem}

\noindent In the proof of Theorem \ref{t:KL} below we shall use the following relation \cite[formula (31)]{LT2}:
\begin{equation} \label{e:pex}
P^-_{\eta\cdot s_i,\zeta} = -q^{-1} P^-_{\eta,\zeta}\qquad \text{if} \; \zeta_i > \zeta_{i+1}, \: \eta_i > \eta_{i+1}.
\end{equation}
Here $i=0,1,\dots, r-1.$ 

\subsubsection{A left representation of $\affH.$}
For $b \in B^l,$ let $H_b$ be the parabolic subalgebra of $\affH$ generated by $T_w$ $(w \in W_b).$ Let $\K\Ib$ be a one-dimensional left representation of $H_b$ defined by  
$$  T_i\cdot\Ib = -q\,\Ib \qquad (s_i \in W_b). $$  
The induced left representation 
$$  \affH\otimes_{H_b}\K\Ib $$
of $\affH$ has as its basis $(T_{x^{-1}})^{-1} \otimes_{H_b} \Ib$ $(x \in \affWbu).$ For $x\in \affWbu,$ define 
$$ |\eta ) := (T_{x^{-1}})^{-1} \otimes_{H_b} \Ib, \qquad \text{where $\quad \eta = x\cdot b,$ $\eta \in P.$} $$  
Then $|\eta )$ $( \eta \in \affW\!\cdot b)$ is a basis of $  \affH\otimes_{H_b}\K\Ib, $ and $|\eta )$ $( \eta \in P)$ is a basis of $  \oplus_{b \in B^l} \affH\otimes_{H_b}\K\Ib. $ Proposition \ref{p:LEFTACTIONofW} allows to describe the action of the affine Hecke algebra  in the basis $|\eta )$ explicitly:    

\begin{align}
 & T_i\cdot |\eta ) = \begin{cases}  \;|s_i\cdot \eta )  & \quad \text{if $\eta_i < \eta_{i+1},$} \\   
                                   \;-q\, |\eta )  & \quad \text{if $\eta_i = \eta_{i+1},$} \\ 
                                 \;|s_i\cdot \eta )  - (q-q^{-1})\, | \eta ) & \quad \text{if $\eta_i > \eta_{i+1}$}
\end{cases}\qquad (0\leq i < r), \label{e:LAH}\\ 
& T_{\pi}\cdot|\eta )  = |\pi\cdot\eta).  \nonumber
\end{align}

\subsection{Wedge product} \label{ss:wedgeproduct}

For $a\in A^n$ and $b\in B^l$ define a vector space $\aWeb$ by 
$$ \aWeb := \K\Ia\otimes_{H_a}\affH\otimes_{H_b}\K\Ib.$$ 
Note that maps 
\begin{align}
&\K\Ia\otimes_{H_a}\affH \otimes_H H \otimes_{H_b} \K\Ib {\rightarrow} \aWeb  \; : \; 
 \Ia\otimes\hat{h}\otimes h\otimes\Ib \mapsto \Ia\otimes \hat{h} h \otimes\Ib, \label{e:I1}\\  
&\K\Ia\otimes_{H_a} H \otimes_H \affH \otimes_{H_b} \K\Ib {\rightarrow}  \aWeb \; : \;
 \Ia\otimes h\otimes\hat{h}\otimes\Ib \mapsto \Ia\otimes h \hat{h}\otimes \Ib\label{e:I2}
\end{align}
are isomorphisms of vector spaces. Let  
$$ \Lambda^r(a) := \bigoplus_{b \in B^l} \Lambda^r(a,b)\;\; (a\in A^n),\quad \Lambda^r(b) := \bigoplus_{a \in A^n} \Lambda^r(a,b)\;\; (b\in B^l),\quad \Lambda^r:= \bigoplus_{b\in B^l} \Lambda^r(b).$$
Let  $v_1,\dots,v_n$ (resp. $\vd_1,\dots,\vd_l$) be a basis of $\K^n$ (resp. $\K^l$). With a sequence  
$$
v_{c_1}X^{\mu_1}\vd_{d_1}, \dots,  v_{c_r}X^{\mu_r}\vd_{d_r}, \quad \text{where $\quad v_{c_i}X^{\mu_i}\vd_{d_i} \in (\K^n\otimes\K^l)[X,X^{-1}],$} 
$$
 we associate unique $a \in A^n,$ $b\in B^l,$ and unique $u \in \auW,$ $v\in \Wbu,$ such that   
$$ c=(c_1,\dots,c_r) = a\cdot u,\quad d=(d_1,\dots,d_r) = v\cdot b,$$
and define the following vector of $\Lambda^r(a,b)$ (here $\omega_b$ is the longest element of $\Wbd$):
\begin{equation}
 (v_{c_1}X^{\mu_1}\vd_{d_1})\wedge  \cdots  \wedge (v_{c_r}X^{\mu_r}\vd_{d_r}):= (-q^{-1})^{l(\omega_b)}\, \Ia\otimes_{H_a}T_u X^{\mu} (T_{v^{-1}})^{-1} \otimes_{H_b} \Ib. \label{e:wedge}
\end{equation} 
Note that using the isomorphism (\ref{e:I1}) to identify the vector spaces we have 
\begin{equation}
 (v_{c_1}X^{\mu_1}\vd_{d_1})\wedge  \cdots  \wedge (v_{c_r}X^{\mu_r}\vd_{d_r}) = (-q^{-1})^{l(\omega_b)} \, ( c|\cdot X^{\mu}\otimes_H |d), \label{e:w1}
\end{equation} 
and using the isomorphism (\ref{e:I2}) we have 
\begin{equation} (v_{c_1}X^{\mu_1}\vd_{d_1})\wedge  \cdots  \wedge (v_{c_r}X^{\mu_r}\vd_{d_r}) = (-q^{-1})^{l(\omega_b)} \, ( c|\otimes_H X^{\mu}\cdot |d).\label{e:w2} 
\end{equation}

We shall call a vector of the form (\ref{e:wedge}) {\em a wedge}, and the vector space $\Lambda^r$ {\em the wedge product}. In what follows it will often be convenient to use a slightly different indexation of wedges: in notations of (\ref{e:wedge}) put $u_{k_i} := v_{c_i}X^{\mu_i}\vd_{d_i},$ where $k_i := c_i+n(d_i-1)-nl\mu_i.$ Since integers $c_i$ (resp. $d_i$) range from $1$ to $n$ (resp. from $1$ to $l$), a wedge (hence $c,d,\mu,a,b,u,v$)  is completely determined by the sequence $\kb = (k_1,k_2,\dots,k_r).$ To emphasize this we write 
\begin{equation} \label{e:x(k)}
c=c(\kb),\; d=d(\kb),\; \mu=\mu(\kb),\; a=a(\kb) ,\; b=b(\kb),\; u=u(\kb),\;v=v(\kb).
\end{equation}
  Denote the left-hand side of (\ref{e:wedge}) by     
$$ u_{\kb} = u_{k_1}\wedge \cdots \wedge u_{k_r}.$$ 
Then  $u_{\kb}$ $( \kb \in \P:=\Z^r )$ is a spanning set of $\Lambda^r.$ However, the vectors of this set are not linearly independent. Indeed, using \eg (\ref{e:w1}) it follows that there are relations among these vectors coming from  
\begin{equation} \label{e:T=T}
 (c|\cdot X^{\mu} T_i \otimes_H | d) = (c|\cdot X^{\mu} \otimes_H T_i\cdot | d)\qquad (i=1,2,\dots,r-1).
\end{equation}   
The exchange formula for $X^{\mu}$ and $T_i$ of  Lemma  \ref{l:Ki} (ii), and the formulas for the action of $T_i$ on $(c|$ and $|d)$ given, respectively, by (\ref{e:RAH}) and (\ref{e:LAH}), allow to compute the relations among the wedges explicitly. Note that the relations for general $r$ follow from those for $r=2.$

Let us call a wedge $u_{\kb}$ {\em ordered} if $\kb \in \P^{++}:=\{ \kb \in \P\:|\: k_1 > k_2 > \cdots > k_r \}.$  

\begin{proposition} \mbox{} \label{p:RULES}\\
{\em (i)} Let $r=2.$ For integers $k_1,k_2$ such that $k_1 \leq k_2,$ let $c_i \in \{1,\dots,n\},$ $d_i \in \{1,\dots,l\},$ $\mu_i \in \Z$ be the unique numbers satisfying $k_i = c_i + n(d_i - 1) - nl\mu_i.$ Let $\gamma$ {\em (}resp. $\delta${\em )} be the residue of $c_2-c_1$ {\em (}resp. $n(d_2-d_1)${\em )} modulo $nl.$ Then    
\begin{alignat}{2}
 u_{k_1}\wedge u_{k_2} = & - u_{k_2}\wedge u_{k_1},  \qquad & & \text{{\em if} $\gamma=0,\delta=0$},   \tag*{(R1)}
\end{alignat}
\begin{alignat}{2}
 u_{k_1}\wedge u_{k_2} =  & - q^{-1} u_{k_2}\wedge u_{k_1} + & & \tag*{(R2)} \\
  & + (q^{-2}-1)\sum_{m\geq 0} q^{-2m} u_{k_2 - \gamma -nlm} \wedge u_{k_1 + \gamma + nlm}\; -  & & \nonumber\\
& - \; (q^{-2}-1)\sum_{m\geq 1} q^{-2m+1} u_{k_2 - nlm} \wedge u_{k_1 + nlm}, & & \text{{\em if} $\gamma >0,\delta=0$}, \nonumber  
\end{alignat}
\begin{alignat}{2}
 u_{k_1}\wedge u_{k_2} =  & q u_{k_2}\wedge u_{k_1} +  & & \tag*{(R3)} \\
  & + (q^2-1)\sum_{m\geq 0} q^{2m} u_{k_2 - \delta -nlm} \wedge u_{k_1 + \delta + nlm}\;+  & & \nonumber\\
& +  \;(q^2-1)\sum_{m\geq 1} q^{2m-1} u_{k_2 - nlm} \wedge u_{k_1 + nlm}, & & \text{{\em if} $\gamma =0,\delta >0$},\nonumber \end{alignat}
\begin{alignat}{2} 
u_{k_1}\wedge u_{k_2} = &   u_{k_2}\wedge u_{k_1} + & & \tag*{(R4)}\\
&+(q-q^{-1})\sum_{m\geq 0}\frac{(q^{2m+1}+q^{-2m-1})}{(q+q^{-1})} \; 
u_{k_2 - \delta -nlm} \wedge u_{k_1 + \delta + nlm} + & & \nonumber \\
& +(q-q^{-1})\sum_{m\geq 0}\frac{(q^{2m+1}+q^{-2m-1})}{(q+q^{-1})} \;
u_{k_2 - \gamma -nlm} \wedge u_{k_1 + \gamma + nlm} + & & \nonumber \\   
&+(q-q^{-1})\sum_{m\geq 1}\frac{(q^{2m}-q^{-2m})}{(q+q^{-1})} \; 
u_{k_2 +nl - \gamma - \delta -nlm} \wedge u_{k_1 -nl +\gamma + \delta + nlm} +\!\!\!\!\!\!\!\!\!\!\!\!\!\negthickspace \negthickspace \negthickspace \!\!\!\!\!\!\!  & & \nonumber \\   
&+(q-q^{-1})\sum_{m\geq 1}\frac{(q^{2m}-q^{-2m})}{(q+q^{-1})}  \; 
u_{k_2  -nlm} \wedge u_{k_1 + nlm}, & & \text{{\em if} $\gamma >0,\delta >0$}, \nonumber 
\end{alignat}
where summations continue as long as  wedges  appearing under the sums  remain ordered. \\
{\em (ii)} Let $r>2,$ then the relations of {\em (i)} hold in every pair of adjacent factors of $u_{k_1}\wedge u_{k_2} \wedge \cdots \wedge u_{k_r}.$
\end{proposition}
\noindent It follows from this proposition, that ordered wedges span $\Lambda^r.$ The relations (R1--R4) then can be thought of as ordering rules that allow to straighten an arbitrary wedge as a linear combination of ordered wedges.  

\begin{remark}{\rm  (i)
In order to compute the ordering rules of the wedges one can use, instead of the isomorphism (\ref{e:I1}) and relations (\ref{e:T=T}), the isomorphism (\ref{e:I2}) and relations $(c|\cdot T_i \otimes_H X^{\mu}\cdot | d) = (c|\otimes_H T_iX^{\mu}\cdot | d).$ The result is easily seen to be the same.\\
(ii) The ordering rules given in Proposition \ref{p:RULES} differ from those used in \cite{TU,U}. This difference is due to a different definition of wedges adopted here. In the present notations, the wedge of \cite{U} is  $(-1)^{l(v)}u_{\kb},$ where $v$ is the same as in (\ref{e:wedge}). 
}
\end{remark}
The next lemma follows easily from Proposition \ref{p:RULES}.
\begin{lemma}  \label{l:ZERO}
 Let $k\geq m,$ then 
\begin{align}  &u_m\wedge u_{k}\wedge u_{k-1}\wedge \cdots \wedge u_m =0,  \tag*{(i)} \\
               &u_{k}\wedge u_{k-1}\wedge \cdots \wedge u_m\wedge u_k =0.\tag*{(ii)}
\end{align}
\end{lemma}
\medskip 

\noindent Now our aim is to show that ordered wedges form a basis of the wedge product. To this end for $b \in B^l$ and $\zeta \in P$ we define $[\zeta ]_b \in \Lambda^r(b)$ by 
$$ [\zeta ]_b := (-q^{-1})^{l(\omega_b)}\,(\zeta |\otimes_{H_b}\Ib.$$
Then  (\ref{e:RAH}) implies  
\begin{equation} \label{e:ze}
[\zeta]_b  = \begin{cases} 0 & \text{if $\zeta_i=\zeta_{i+1},$ $b_i=b_{i+1},$} \\
 -q^{-1} [\zeta\cdot s_i]_b & \text{if $\zeta_i<\zeta_{i+1},$ $b_i=b_{i+1}.$} 
\end{cases}
\end{equation}
From this it follows that $ [\zeta]_b $ $(\zeta \in P^{++}_b),$ where $P^{++}_b:= \{ \zeta \in P \: |\: (\alpha_i,\zeta)>0 \;\;\text{ if $ \: b_i=b_{i+1}$}\},$ span $\Lambda^r(b).$ For $b$ such that $P^{++}_b = P^{++}$ a proof of the following lemma is given in \cite{LT2}, a proof for general $b$ is  completely similar.   
\begin{lemma} \label{l:BASIS}
For each $b \in B^l$ the set $\{ [\zeta]_b \:|\:\zeta \in P^{++}_b \}$ is a basis of $\Lambda^r(b).$  
\end{lemma}

For $\kb \in \P,$ we define $\zeta(\kb) \in P$ by 
$$ \zeta(\kb):= v(\kb)^{-1}\cdot( c(\kb) - n \mu(\kb)),$$ 
where $v(\kb),c(\kb)$ and $\mu(\kb)$ are defined in (\ref{e:x(k)}).
\begin{proposition} \label{p:CON}
Suppose  $\kb \in \P^{++}.$ Then $\zeta(\kb) \in P^{++}_{b(\kb)}$ and $u_{\kb} = [\zeta(\kb)]_{b(\kb)}.$
%Let $\kb \in \P^{++},$ then $u_{\kb} = [\zeta]_b,$ where $\zeta = \zeta(\kb)$ %and $b=b(\kb).$ Moreover, $\zeta \in P^{++}_b.$ 
Conversely, for $b \in B^l$ and $\zeta \in P^{++}_b,$ there is $\kb \in \P^{++}$ such that $b=b(\kb),$ $\zeta = \zeta(\kb).$  
\end{proposition}
\proof First, let us show that $\kb \in \P^{++}$ implies $u_{\kb} = [\zeta(\kb)]_{b(\kb)}.$ We put $\zeta:=\zeta(\kb),$ and use notations of (\ref{e:wedge}) and (\ref{e:x(k)}). Let $\lambda:=-\mu.$ Observe that $\kb \in P^{++}$  implies 
\begin{align}
 &\lambda \in P^+, \label{e:CON1} \\ 
 & i<j,\: d_i<d_j  \Longrightarrow   \lambda_i > \lambda_j,  \label{e:CON2} \\
 & i<j, \: d_i = d_j  \Longrightarrow  c_i + n\lambda_i > c_j + n\lambda_j. \label{e:CON3}
\end{align}
\mbox{}From the dominance of $\lambda$ it follows that 
$
X^{\mu} = T_{t_{\lambda}}.
$
 Using $v\in W^b,$ Lemma \ref{l:S(v^{-1})} and (\ref{e:CON2}) we get $\alpha \in S(v^{-1}) \Longrightarrow (\alpha , \lambda) > 0.$ Since $S(v^{-1}) = - v(S(v)),$ this gives $\alpha \in S(v) \Longrightarrow (v(\alpha),\lambda) < 0,$ \ie $S(v) \subset \{\alpha \in R^+\: |\: (\alpha,\nu) < 0\},$ where we put $\nu:=v^{-1}(\lambda).$ Let $w\in W$ be the shortest element such that $w(\nu) \in P^+,$ clearly, $w(\nu) = \lambda.$ By Lemma \ref{l:SHORTEST}, $S(w) = \{ \alpha \in R^+\:|\: (\alpha,\nu)\},$ hence $S(v) \subset S(w).$ On the other hand, the length of $w$ does not exceed the length of $v.$ This is possible only if $S(v) = S(w),$ hence $v=w.$ Thus, $S(v) = \{\alpha \in R^+\: |\: (v(\alpha),\lambda) < 0\}.$ 
Now, Proposition \ref{p:SPLIT} implies $l(u \,t_{\lambda} v)= l(u) + l(t_{\lambda}) - l(v).$ This, $l(u) + l(t_{\lambda}) = l(u \,t_{\lambda})$ and  $X^{\mu} = T_{t_{\lambda}}$ give
$$ T_u \,X^{\mu} (T_{v^{-1}})^{-1} = T_{u t_{\lambda} v}.$$  
Since $u \in \auW,$ by Lemma \ref{l:WaWa}, we have $u\, t_{\lambda} v \in \auaffW,$ hence 
$$ 
 u_{\kb} =  [ a\cdot u\, t_{\lambda} v ]_b. 
$$ 
It remains to observe  that 
$$ a\cdot u\, t_{\lambda} v = (c + n\lambda )\cdot v = v^{-1} \cdot (c + n\lambda ) = \zeta.$$ 

Next, we show that $\kb \in \P^{++}$ implies $\zeta \in P^{++}_b.$ 
It follows from $v\in \Wbu$ and Lemma \ref{l:R_b} that for $i<j,$ 
$$
 b_i=b_j \Rightarrow v^{-1}(i) < v^{-1}(j),\;\;\text{equivalently,}\; \; \,d_{v^{-1}(i)} = d_{v^{-1}(j)}  \Rightarrow v^{-1}(i) < v^{-1}(j).
$$ 
Therefore, using $c_{v^{-1}(i)} + n \lambda_{v^{-1}(i)} = \zeta_i,$ and (\ref{e:CON3}) we get  $\zeta \in P^{++}_b.$ 

 Finally, $u_{\kb}$ with $\kb \in \P^{++}$ span $\Lambda^r.$ Hence, for $b \in B^l$ and $\zeta \in P^{++}_b,$ we have 
$$[\zeta]_b = \sum_{\lb \in \P^{++}} e_{\lb} u_{\lb} \qquad  (e_{\lb} \in \K). 
$$ 
But $\{[\zeta]_b\:|\: b \in B^l,\: \: \zeta \in P^{++}_b\}$ is a basis of $\Lambda^r.$ Therefore writing $u_{\lb}$ as $[\zeta(\lb)]_{b(\lb)},$  we get $[\zeta]_b = u_{\kb},$ where $\kb \in \P^{++}$ is such that  $b = b(\kb)$ and  $\zeta = \zeta(\kb).$  \finpf  

\noindent This proposition and Lemma \ref{l:BASIS} immediately imply 
\begin{proposition}
$\{ u_{\kb}\:|\: \kb \in \P^{++}\}$ is a basis of $\Lambda^r.$ 
\end{proposition}

\subsection{Canonical bases of the wedge product} \label{ss:CBfin}
Define an involution $x \mapsto \ov{x}$ of $\Lambda^r$ by 
$$ \textstyle{\ov{\Ia\otimes_{H_a}h \otimes_{H_b}\Ib} = {\Ia\otimes_{H_a}\ov{h}\otimes_{H_b}\Ib}}, \qquad \textstyle{\ov{q} = q^{-1}} \qquad ( h \in \affH ). 
$$ 

\begin{lemma} \label{l:INVO2}
Let $u\in \auW,$ $v\in \Wbu,$ $\mu \in P.$ Then $\omega_a u \omega \in \auW,$ 
$ \omega v \omega_b \in \Wbu,$ and 
$$ \textstyle{\ov{T_u X^{\mu} (T_{v^{-1}})^{-1}}} = \textstyle{T_{\omega_a} T_{\omega_a u \omega} X^{\omega(\mu)} (T_{(\omega v \omega_b)^{-1}})^{-1} (T_{\omega_b})^{-1}}.$$
\end{lemma}
\proof Since for $w\in W$ we have $l(w \omega) = l(\omega) - l(w),$ $u\omega$ is the longest element of the coset $W_a u\omega.$ Hence $\omega_a  u\omega$ is the shortest element of this coset, \ie $\omega_a  u\omega \in \auW$. Therefore $l(u\omega) = l(\omega) + l(\omega_a u\omega)$ and    
$$ \textstyle{T_{u\omega} = T_{\omega_a} T_{\omega_a u \omega}.}$$ 
In a completely similar fashion we get $ \omega v \omega_b \in \Wbu$ and 
$$\textstyle{(T_{(\omega v)^{-1}})^{-1} = (T_{(\omega v \omega_b)^{-1}})^{-1} (T_{\omega_b})^{-1}}.$$
 Lemma \ref{l:INVO1} implies  the remaining statement. \finpf

\begin{proposition} \label{p:INVO}
For $u_{k_1}\wedge u_{k_2}\wedge \cdots \wedge u_{k_r} \in \Lambda^r(a,b)$ we have   
$$ \ov{u_{k_1}\wedge u_{k_2}\wedge  \cdots \wedge u_{k_r}} = (-q)^{l(\omega_b)} q^{-l(\omega_a)} {u_{k_r}\wedge \cdots \wedge u_{k_2}\wedge u_{k_1}}.$$ 
\end{proposition}
\proof Using (\ref{e:w1}) in the right-hand side of (\ref{e:wedge}) we have, by Lemma \ref{l:INVO2}, 
$$ \textstyle{ \ov{\Ia\otimes_{H_a} T_u X^{\mu} (T_{v^{-1}})^{-1} \otimes_{H_b} \Ib}} = q^{-l(\omega_a)} (-q)^{-l(\omega_b)}(a\cdot u\omega | \cdot X^{\omega(\mu)} \otimes_H |\omega v\cdot b).$$ 
The result follows. \finpf

\begin{remark} \label{r:leng} {\rm 
It is easy to see that in notations of (\ref{e:wedge}) we have 
$$ l(\omega_a) = \sharp\{i<j\:|\: c_i=c_j\},\qquad l(\omega_b) = \sharp\{i<j\:|\: d_i=d_j\}.$$
}
\end{remark}

\noindent For $\kb \in \P^{++}$ put 
$$ \ov{u_{\kb}} = \sum_{\lb \in \P^{++}} R_{\kb,\lb}(q)\: u_{\lb}.$$
By Proposition \ref{p:RULES}, the entries of matrix  $\|R_{\kb,\lb}(q)\|$ are  Laurent polynomials in $q$ with integral coefficients, and by Remark \ref{r:leng} we have $R_{\kb,\kb}(q) =1.$ 

We define a partial order on $\P^{++}$ by 
$$ \kb \geq \lb \quad \text{iff} \quad \left|\begin{array}{l l} \sum_{i=1}^j k_i \geq \sum_{i=1}^j l_i & \text{for all $j=1,\dots,r,$} \\ \sum_{i=1}^r k_i = \sum_{i=1}^r l_i, &  \\ k_r \leq l_r. &      \end{array} \right. 
$$
The ordering rules of Proposition \ref{p:RULES} imply that matrix $\|R_{\kb,\lb}(q)\|$ is lower triangular with respect to this order. That is $R_{\kb,\lb}(q)$ is not zero only if $\kb \geq \lb.$  

Define two lattices of $\Lambda^r$ by 
$$ {\mathcal{ L}}^+ := \bigoplus_{\kb \in \P^{++}} \Q[q] \, u_{\kb},\qquad  {\mathcal{ L}}^- := \bigoplus_{\kb \in \P^{++}} \Q[q^{-1}]\, u_{\kb}.$$ 
The unitriangularity of  $\|R_{\kb,\lb}(q)\|$ implies, by the standard argument going back to Kazhdan and Lusztig, that    
\begin{theorem}
There are unique bases $\{ G^+_{\kb}\:| \: \kb \in \P^{++}\},$ $\{ G^-_{\kb}\:| \: \kb \in \P^{++}\}$ of $\Lambda^r$ such that 
\begin{alignat*}{2}
&\textstyle{\ov{G^+_{\kb}} = G^+_{\kb},}&\qquad &\textstyle{\ov{G^-_{\kb}} = G^-_{\kb},} \tag*{(i)} \\ 
 &G^+_{\kb} \equiv u_{\kb}\bmod q{\mathcal{ L}}^+,&\qquad &G^-_{\kb} \equiv u_{\kb}\bmod q^{-1}{\mathcal{ L}}^-.  \tag*{(ii)}
\end{alignat*}
\end{theorem}
\noindent Set  
$$ G^+_{\kb} = \sum_{\lb \in \P^{++}} \dplus_{\kb,\lb}(q)\: u_{\lb},\qquad G^-_{\kb} = \sum_{\lb \in \P^{++}} \dminus_{\kb,\lb}(q)\: u_{\lb}.$$ 
It is clear that $\dplus_{\kb,\lb}(q)$ or $\dminus_{\kb,\lb}(q)$ is non-zero  only if $u_{\kb}$ and $u_{\lb}$ belong to the same subspace $\Lambda^r(a,b).$ That is, in notations of (\ref{e:x(k)}), only if $a(\kb) = a(\lb)$ and $b(\kb) = b(\lb).$ 

\begin{theorem} \label{t:KL}
For $\kb,\lb \in \P^{++}$ such that $a(\kb) = a(\lb),$ $b(\kb) = b(\lb),$ put $\xi=\zeta(\kb),$ $\eta=\zeta(\lb)$ and $a=a(\kb),$ $b=b(\kb).$ Then   
\begin{align}
&\dminus_{\kb,\lb}(q) = P^{-}_{\eta,\,\xi},\tag*{(i)} \\
&\dplus_{\kb,\lb}(q) = \sum_{v \in W_b} (-q)^{l(v)} P^+_{\eta\cdot\omega_bv,\,\xi\cdot\omega_b}, \tag*{(ii)}
\end{align}
where $P^{-}_{\eta,\xi},$ $P^{+}_{\eta,\xi}$ are the parabolic Kazhdan-Lusztig polynomials associated with $\Ia\otimes_{H_a}\affH.$
\end{theorem}

\proof Put $D_{\xi} = C^-_{\xi}\otimes \Ib.$ Then $\ov{D_{\xi}} = D_{\xi}.$ Using (\ref{e:pex}) and  (\ref{e:ze}) we obtain 
$$ 
D_{\xi} = (-q)^{l(\omega_b)} \sum_{\eta \in P} \, P^-_{\eta,\xi} \:[\eta]_b = z_b \: \sum_{\eta \in P^{++}_b} \,P^-_{\eta,\xi} \: [\eta]_b, 
$$ 
where  
$ z_b  = (-q)^{l(\omega_b)} \sum_{v \in W_b} q^{-2 l(v)} .$ Since $\ov{z_b} = z_b,$ we have $D_{\xi} = z_b\, G^-_{\kb}$ and (i) follows.      

Next let  $E_{\xi} = C^-_{\xi\cdot\omega_b}\otimes \Ib.$ Then $\ov{E_{\xi}} = E_{\xi},$ and 
$$ 
E_{\xi} = \sum_{\eta \in P^{++}_b} \sum_{v \in W_b} (-q)^{l(\omega_b) - l(v)}  P^+_{\eta\cdot v,\xi\cdot \omega_b} \:[\eta]_b = \sum_{\eta \in P^{++}_b} \sum_{v \in W_b} (-q)^{l(v)}  P^+_{\eta\cdot \omega_b v,\xi\cdot \omega_b} \:[\eta]_b.
$$ 
Hence $E_{\xi} = G^+_{\kb},$ which implies (ii). \finpf 

\begin{remark} {\rm 
For $\xi \in a\cdot\affW,$ let $x(\xi)$ be the unique element of $\auaffW$ such that $ \xi = a\cdot x(\xi).$ It follows from Proposition \ref{p:RIGHTACTIONofW} that for $\xi \in P^{++}_b$ we have $ x(\xi\cdot v) = x(\xi)v$ for all $v\in W_b.$ Hence, one can rewrite the formulas of the theorem in terms of (ordinary) Kazhdan-Lusztig polynomials for $\affH$ as 
\begin{align*}
&\dminus_{\kb,\lb}(q) = \sum_{u \in W_a} q^{-l(u)} {\mathcal {P}}^{-}_{u\, x(\eta),\,x(\xi)},\tag*{(i)} \\
&\dplus_{\kb,\lb}(q) = \sum_{v \in W_b} (-q)^{l(v)} {\mathcal {P}}^+_{\omega_a\, x(\eta)\,\omega_b \,v,\,\omega_a \,x(\xi)\,\omega_b}. \tag*{(ii)}
\end{align*}

}
\end{remark}

\subsection{Actions of quantum affine algebras on the wedge product}

It is easy to verify that formulas 
\begin{alignat*}{2}
&\en_i(v_c X^{m})  =  \delta_{i+1\equiv c\bmod n} \;v_{c-1} X^{m + \delta_{i0}},  &
&\qquad   \tn_i(v_c X^m )  =  q^{\delta_{i\equiv c\bmod n} - \delta_{i+1\equiv c\bmod n}} \; v_c X^m,  \\  
&\fn_i(v_c X^{m})  =  \delta_{i\equiv c\bmod n} \;v_{c+1} X^{m - \delta_{i0}},  &
&\qquad   \dn(v_c X^m )  = m v_c X^m. 
\end{alignat*}
define a level zero action of $\UN$ on $\K^n[X,X^{-1}].$ Here it is understood that $v_0 = v_{n},$  $v_{n+1} = v_1.$ Also, for a statement $S$ we put $\delta_S=1$ (resp. $0$) if $S$ is true (resp. false). 
Using the coproduct 
\begin{alignat*}{2}
&\Delta(\en_i)  = \en_i\otimes {\tn_i}^{-1} + 1\otimes \en_i,&
&\qquad \Delta(\tn_i) = \tn_i\otimes \tn_i,  \\
&\Delta(\fn_i) = \fn_i \otimes 1 + \tn_i\otimes \fn_i, &
&\qquad \Delta(\dn) = \dn\otimes 1 + 1 \otimes \dn, 
\end{alignat*}
we extend this action on the tensor product $(\K^n[X,X^{-1}])^{\otimes r}.$ Next, using the isomorphism 
$$ 
(\K^n[X,X^{-1}])^{\otimes r} \;\overset{\sim}{\longrightarrow}\; \bigoplus_{a\in A^n} \K\Ia\otimes_{H_a} \affH,\quad v_{c_1}X^{\mu_1}\otimes \cdots  \otimes v_{c_r}X^{\mu_r} \;\mapsto \; (c|\cdot X^{\mu}, 
$$ 
where $c = (c_1,\dots,c_r),$ $\mu = (\mu_1,\dots,\mu_r),$ to identify the involved  vector spaces, we obtain an action of   $\UN$ on the $\affH$--module $M_R=\oplus_{a\in A^n} \K\Ia\otimes_{H_a} \affH.$ 
%Let $\UNp$ be the subalgebra of $\UN$ generated by $\en_i,$ $\fn_i,$ $\tn_i$ and ${\tn_i}^{-1}.$ 
The following proposition is due to \cite{GRV}, it is easily verified by reducing to the case $r=2$ (cf. proof of Proposition 7.1 in \cite{LT2}).     
\begin{proposition} \label{p:UNCOMMH}
%The actions of $\UNp$ and $\affH$ on $M_R$ commute. 
The actions of $\UN$ and $H\subset \affH$ on $M_R$ commute.   
\end{proposition}
\noindent By the isomorphism (\ref{e:I1}), for any $b \in B^l,$ the commutativity of  $\UN$ and $H$ allows to restrict the action of $\UN$ on $\Lambda^r(b),$ which gives then an action of $\UN$ on $\Lambda^r.$ In terms of the wedge vectors  $u_{\kb} = u_{k_1}\wedge \cdots \wedge u_{k_r}$ $(\kb = (k_1,\dots,k_r) \in \P),$  this action is  written as   
\begin{align}
&\en_i(u_{\kb}) = \textstyle{\sum_{j=1}^r u_{k_1}\wedge \cdots \wedge u_{k_{j-1}}\wedge \en_i(u_{k_j})\wedge {\tn_i}^{-1}(u_{k_{j+1}})\wedge \cdots \wedge {\tn_i}^{-1}(u_{k_r})},     \label{e:en}\\ 
&\fn_i(u_{\kb}) = \textstyle{\sum_{j=1}^r \tn_i(u_{k_1})\wedge \cdots \wedge \tn_i(u_{k_{j-1}})\wedge \fn_i(u_{k_j})\wedge u_{k_{j+1}}\wedge \cdots \wedge u_{k_r}}, \label{e:fn}\\                  
&\tn_i(u_{\kb}) =  \textstyle{ \tn_i(u_{k_1})\wedge \cdots \wedge \tn_i(u_{k_r})}, \label{e:tn} \\
&\dn(u_{\kb}) =   \textstyle{ \sum_{j=1}^r u_{k_1}\wedge \cdots \wedge u_{k_{j-1}}\wedge \dn(u_{k_j})\wedge u_{k_{j+1}}\wedge \cdots \wedge u_{k_r}}.\label{e:dn}
\end{align}
Here we put $e_i(u_{k_j})= e_i(v_{c_j}X^{\mu_j}\vd_{d_j}):=e_i(v_{c_j}X^{\mu_j})\vd_{d_j}$ and similarly for the rest of the generators.   

\bigskip 

In a completely similar fashion we define on the wedge product an action of the quantum affine algebra $\UL,$ where $p:=-q^{-1}.$ In order to distinguish between this action and the action of $\UN$ we put dots over the generators of $\UL.$ Define a level zero action of $\UL$ on $\K^l[X,X^{-1}]$ by    
\begin{alignat*}{2}
&\el_i(X^{m}\vd_d)  =  \delta_{i+1\equiv d\bmod l} \;X^{m + \delta_{i0}}\vd_{d-1},  &
&\qquad   \tl_i( X^m \vd_d)  =  p^{\delta_{i\equiv d\bmod l} - \delta_{i+1\equiv d\bmod l}} \; X^m\vd_d,  \\  
&\fl_i( X^{m}\vd_d)  =  \delta_{i\equiv d\bmod l} \; X^{m - \delta_{i0}}\vd_{d+1} ,  &
&\qquad   \dl(X^m  \vd_d  )  = m X^m \vd_d. 
\end{alignat*}
Here it is understood that $\vd_0 = \vd_{l},$  $\vd_{l+1} = \vd_1.$ 
Using the coproduct 
\begin{alignat*}{2}
&\Delta(\el_i)  = \el_i\otimes {\tl_i}^{-1} + 1\otimes \el_i,&
&\qquad \Delta(\tl_i) = \tl_i\otimes \tl_i,  \\
&\Delta(\fl_i) = \fl_i \otimes 1 + \tl_i\otimes \fl_i, &
&\qquad \Delta(\dl) = \dl\otimes 1 + 1 \otimes \dl, 
\end{alignat*}
we extend this action on the tensor product $(\K^l[X,X^{-1}])^{\otimes r}.$ Next, using the isomorphism 
$$ 
(\K^l[X,X^{-1}])^{\otimes r} \;\overset{\sim}{\longrightarrow}\; \bigoplus_{b\in B^l} \affH\otimes_{H_b} \K\Ib,\quad X^{\mu_1}\vd_{d_1}\otimes \cdots  \otimes X^{\mu_r}\vd_{d_r} \;\mapsto \; X^{\mu}\cdot |d)\, (-q^{-1})^{l(\omega_b)}, 
$$ 
where $d = (d_1,\dots,d_r),$ $\mu = (\mu_1,\dots,\mu_r),$ and $b$ is the unique point of $B^l$ in the orbit $W\cdot d,$ to identify the involved  vector spaces, we obtain an action of   $\UL$ on the $\affH$--module $M_L=\oplus_{b\in B^l} \affH\otimes_{H_b} \K\Ib.$ 
%Let $\ULp$ be the subalgebra of $\UL$ generated by $\el_i,$ $\fl_i,$ $\tl_i$ and ${\tl_i}^{-1}.$ 
The following result is  an analogue of Proposition \ref{p:UNCOMMH} and is verified similarly.
\begin{proposition} \label{p:ULCOMMH}
%The actions of $\ULp$ and $\affH$ commute. 
The actions of $\UL$ and $H\subset \affH$ on $M_L$ commute.   
\end{proposition}
\noindent By the isomorphism (\ref{e:I2}), for any $a \in A^n,$ the commutativity of  $\UL$ and $H$ allows to restrict the action of $\UL$ on $\Lambda^r(a),$ which gives then an action of $\UL$ on $\Lambda^r.$ In terms of the wedge vectors  $u_{\kb} = u_{k_1}\wedge \cdots \wedge u_{k_r}$ $(\kb = (k_1,\dots,k_r) \in \P),$  this action is  written as   
\begin{align}
&\el_i(u_{\kb}) = \textstyle{\sum_{j=1}^r u_{k_1}\wedge \cdots \wedge u_{k_{j-1}}\wedge \el_i(u_{k_j})\wedge {\tl_i}^{-1}(u_{k_{j+1}})\wedge \cdots \wedge {\tl_i}^{-1}(u_{k_r})},   \label{e:el}   \\ 
&\fl_i(u_{\kb}) = \textstyle{\sum_{j=1}^r \tl_i(u_{k_1})\wedge \cdots \wedge \tl_i(u_{k_{j-1}})\wedge \fl_i(u_{k_j})\wedge u_{k_{j+1}}\wedge \cdots \wedge u_{k_r}}, \label{e:fl}\\                  
&\tl_i(u_{\kb}) =  \textstyle{ \tl_i(u_{k_1})\wedge \cdots \wedge \tl_i(u_{k_r})}, \label{e:tl}\\
&\dl(u_{\kb}) =   \textstyle{ \sum_{j=1}^r u_{k_1}\wedge \cdots \wedge u_{k_{j-1}}\wedge \dl(u_{k_j})\wedge u_{k_{j+1}}\wedge \cdots \wedge u_{k_r}}.\label{e:dl}
\end{align}
Here we put $\el_i(u_{k_j})= \el_i(v_{c_j}X^{\mu_j}\vd_{d_j}):=v_{c_j}\el_i(X^{\mu_j}\vd_{d_j})$ and similarly for the rest of the generators.   

\bigskip 

A well-known result of Bernstein (see \eg \cite{Ki}) says that the centre $Z(\affH)$ of $\affH$ is generated by symmetric Laurent polynomials in $X_i:= X^{\ep_i}.$ It follows, by using either (\ref{e:I1}) or (\ref{e:I2}), that $Z(\affH)$ acts on the wedge product $\Lambda^r.$ This action may be computed in terms of the wedge vectors $u_{\kb} = u_{k_1}\wedge \cdots \wedge u_{k_r}$ $(\kb = (k_1,\dots,k_r) \in \P)$ by using either (\ref{e:w1}) or (\ref{e:w2}). In particular, for  $B_m := \sum_{i=1}^r X_i^m$ $(m\in\Z^*),$ we get
\begin{equation} 
B_m(u_{\kb}) = \textstyle{ \sum_{j=1}^r u_{k_1}\wedge \cdots \wedge u_{k_{j-1}}\wedge u_{k_j - nlm}\wedge u_{k_{j+1}}\wedge \cdots \wedge u_{k_r}}.\label{e:Bo}
\end{equation}

\begin{proposition} \label{p:COMMUTE}
The actions of $\UNp,$ $\ULp$ and $Z(\affH)$ on $\Lambda^r$ are pairwise mutually commutative.
\end{proposition}
\proof The commutativity of $\UNp$ (resp. $\ULp$) with $Z(\affH)$ is immediate by Proposition \ref{p:UNCOMMH} (resp. Proposition \ref{p:ULCOMMH}). The commutativity of $\UNp$ with $\ULp$ follows from (\ref{e:en} -- \ref{e:tn}) and (\ref{e:el} -- \ref{e:tl}).\finpf

\noindent The actions of $\UN,$ $\UL$ and $Z(\affH)$ on $\Lambda^r$ are compatible with the bar involution in the sense that 
\begin{proposition} \label{p:COMP}
For $v \in \Lambda^r$
\begin{alignat*}{5}
&\ov{\en_i(v)} = \en_i(\ov{v}),\quad & &\ov{\fn_i(v)} = \fn_i(\ov{v}),\quad & & \ov{\tn_i(v)} = {\tn_i}^{-1}(\ov{v}),\quad & & \ov{\dn(v)} = {\dn}(\ov{v}),\qquad & & (i=0,1,\dots,n-1),     \\
&  \ov{\el_j(v)} = \el_j(\ov{v}),\quad & &\ov{\fl_j(v)} = \fl_j(\ov{v}),\quad & &  \ov{\tl_j(v)} = {\tl_j}^{-1}(\ov{v}),\quad & & \ov{\dl(v)} = {\dl}(\ov{v}),\qquad & & (j=0,1,\dots,l-1),     \\
&\ov{B_m(v)} = B_m(\ov{v}),& &    &       &   & &    \qquad & &(m\in\Z^*).
\end{alignat*}
\end{proposition}
\proof It is enough to show this for $v = u_{\kb}$ $(\kb =(k_1,\dots,k_r) \in \P).$ As in Section \ref{ss:wedgeproduct}  we put $k_i = c_i + n(d_i-1)-nl\mu_i$ with $c_i \in \{1,\dots,n\},$ $d_i \in \{1,\dots,l\},$ $\mu_i \in \Z.$ By Proposition \ref{p:INVO} and Remark \ref{r:leng} we have    
$$ \ov{u_{\kb}} = q^{-\kappa(c)} (-q)^{\kappa(d)} u_{k_r}\wedge \cdots \wedge u_{k_1} = p^{-\kappa(d)} (-p)^{\kappa(c)} u_{k_r}\wedge \cdots \wedge u_{k_1}.$$ 
where $\kappa(c)$ (resp. $\kappa(d)$) is the number of pairs $(i,j)$ such that $c_i=c_j$ (resp. $d_i=d_j$).
Using  (\ref{e:fl}) we have 
$$ %\begin{align*}
\ov{\fl_0(u_{\kb})} = {(-p)^{\kappa(c)} \sum_{j,\, d_j =l} p^{ -\kappa(d'_{(j)}) + \sum_{i<j} \delta_{d_i=1}-\delta_{d_i=l}} u_{k_r}\wedge \cdots \wedge (v_{c_j}X^{\mu_j-1}\vd_1) \wedge \cdots \wedge u_{k_1}},
$$ 
and
$$
{\fl_0(\ov{u_{\kb}})} = {(-p)^{\kappa(c)} \sum_{j,\, d_j =l} p^{ -\kappa(d)+ \sum_{i>j} \delta_{d_i=l}-\delta_{d_i=1}} u_{k_r}\wedge \cdots \wedge (v_{c_j}X^{\mu_j-1}\vd_1) \wedge  \cdots \wedge u_{k_1}},
$$ %\end{align*}
where  
$$ \kappa(d'_{(j)}) = \kappa(d+\ep_j-l\ep_j) = \kappa(d) + \sum_{i\neq j} \delta_{d_i=1} - \delta_{d_i=l}.$$ 
This  gives $\ov{\fl_0(u_{\kb})} = {\fl_0(\ov{u_{\kb}})}.$ The relations for the rest of the generators of $\UL$ and for the generators of $\UN$ are verified similarly. Finally, $\ov{B_m(v)} = B_m(\ov{v})$ follows from Lemma \ref{l:INVO1}. \finpf

\section{Canonical bases of the $q$-deformed Fock spaces}
\subsection{Semi-infinite wedge product} \label{ss:SIW}
For each integer $s,$ we define a vector space $\SIW$ as the inductive limit $\underset{\rightarrow}{\lim}\Lambda^r$ where maps $\Lambda^r \rightarrow \Lambda^t$ $(t>r)$ are given by 
$$ v \mapsto v\wedge u_{s-r}\wedge u_{s-r-1}\wedge \cdots \wedge u_{s-t+1}.$$ 
The vector space $\SIW$ will be called the semi-infinite wedge product of charge $s.$ For $v \in \Lambda^r$ we shall use the semi-infinite expression 
$$ v\wedge u_{s-r}\wedge u_{s-r-1}\wedge u_{s-r-2}\wedge \cdots $$ 
to denote the image of $v$ with respect to the canonical map from $\Lambda^r$ to  $\underset{\rightarrow}{\lim}\Lambda^r.$ Let $\PP(s)$ be the set of semi-infinite sequences $\kb=(k_1,k_2,\dots\:)$ $\in $ $\Z^{\infty}$ such that $k_i = s-i+1$ for all $i\gg 1.$ By definition, $\SIW$ is spanned by $u_{\kb} := u_{k_1}\wedge u_{k_2} \wedge \cdots $ with $\kb \in \PP(s).$ We shall call a semi-infinite wedge $u_{\kb}$ {\em ordered} if $\kb \in \PP^{++}(s):= \{ (k_1,k_2,\dots\:) \in \PP(s)\:|\: k_1>k_2>\cdots \;\}.$ 
\begin{proposition}[\cite{TU}]
Ordered wedges form a basis of $\SIW.$ 
\end{proposition}

In what follows we shall use, besides the indexation by $\PP^{++}(s),$  three other indexations of the basis formed by ordered wedges. The first of these is the obvious indexation by the set $\Pi$ of all  partitions. Namely, with $\kb = (k_1,k_2,\dots\:) \in \PP^{++}(s)$ we associate a partition $\lambda =(\lambda_1,\lambda_2,\dots\:)$ by taking $\lambda_i = k_i - s + i -1$ for all $i \in \N,$ and we put $|\lambda,s\rangle := u_{\kb}.$   

Another indexation is by the set of pairs $(\lal,\sal)$ where $\lal = (\lambda^{(1)},\dots,\lambda^{(l)})$ is an  $l$-multipartition, and $\sal =(s_1,\dots,s_l)$ is a sequence of integers summing up to $s.$ For any $d\in \{1,\dots,l\}$ the partition $\lambda^{(d)}$ and the integer $s_d$ are defined as follows. For each $i\in \N^{*}$ write $k_i = c_i + n(d_i-1) -nl m_i$ where $c_i \in \{1,\dots,n\},$ $d_i \in \{1,\dots,l\}$ and $m_i \in \Z.$ Then, let $k_1^{(d)}$ be equal $c_i - n m_i$ where $i$ is the smallest such that $d_i =d,$ let $k_2^{(d)}$ be equal 
$c_j - n m_j$ where $j$ is the second smallest such that $d_j =d,$ and so on. In this way we obtain a strictly decreasing semi-infinite sequence $(k_1^{(d)},k_2^{(d)},\dots\:)$ such that $k_i^{(d)} = s_d - i+1$ $(i\gg 1)$ for some unique integer $s_d.$ Now we define the  partition $\lambda^{(d)} = (\lambda^{(d)}_1,\lambda^{(d)}_2,\dots\:)$ by $\lambda^{(d)}_i = k_i^{(d)} - s_d + i - 1.$ It is easy to check that the  $s_d$ obtained in this way satisfy $s_1+\cdots + s_l =s.$

In a completely similar fashion we associate with each $\kb \in \PP^{++}(s)$ a pair $(\lan,\san)$ where $\lan$ is an $n$-multipartition, and $\san$ is a sequence of $n$ integers summing up to $s.$ The $(\lan,\san)$ is obtained by the same procedure as the $(\lal,\sal)$ reversing everywhere  the roles of $n$ and $l,$ and the roles of $c_i$ and $d_i.$

For any $s\in \Z,$ let $\Z^l(s)$ be the set of $l$-tuples of integers summing to $s.$ Define the map $\tau_l^s : \Pi \rightarrow \Pi^l\times \Z^l(s)$ (respectively, the map $\tau_n^s : \Pi \rightarrow \Pi^n\times \Z^n(s)$) by  
\begin{equation}
 \tau_l^s \: : \: \lambda \mapsto (\lal,\sal) ,\quad 
\text{(respectively, by)}\quad 
 \tau_n^s \: : \: \lambda \mapsto (\lan,\san) \label{e:tau}.
\end{equation}
It is not difficult to see that for each $s,$ the maps $\tau_l^s$ and $\tau_n^s$ are bijections. Hence, setting $|\lal,\sal\rangle := |\lambda,s\rangle$ (resp. $|\lan,\san\rangle := |\lambda,s\rangle$) if $(\lal,\sal) = \tau_l^s(\lambda)$ (resp. if  $(\lan,\san) = \tau_n^s(\lambda)),$ we obtain that $B(s):= \{ |\lal,\sal\rangle \: | \: \lal \in \Pi^l,\, \sal \in \Z^l(s)\}$ $=$ $\{ |\lan,\san\rangle \: | \: \lan \in \Pi^n,\, \san \in \Z^n(s)\}$ $=$ $\{ |\lambda,s\rangle \: | \: \lambda\in \Pi\}$  is a basis of $\SIW.$ 

\begin{remark}\label{r:QUO}{\rm  (i)
The pair $(\lan = (\lambda^{(1)},\dots,\lambda^{(n)}), \san=(s_1,\dots,s_n))$ can be read off the  diagram of the partition $\lambda$ in the following manner. For $r\in \{0,1,\dots,n-1\}$ let $R_r$ (resp. $C_r$) be the set of all rows (resp. columns) of $\lambda$  that have a node of residue $r$ as their rightmost (resp. bottom) node. Then for each $c \in \{1,\dots,n\},$ the  diagram of $\lambda^{(c)}$ is embedded into the  diagram of $\lambda$ by  
$$ \lambda^{(c)} = \lambda\cap C_c\cap R_{c-1},$$ 
where we put $C_n = C_0.$ Hence, $\lan$ is the $n$-quotient of $\lambda.$ On the other hand, the sequence $\san$ is obtained as follows. Let $N_r(\lambda |s,n)$ be the number of addable nodes of residue $r$ minus the number of removable nodes of residue $r.$ Then, for each $c\in \{1,\dots,n-1\}$ we have 
$s_c - s_{c+1} = N_c(\lambda |s,n).$ These equalities, together with $\sum_{a=1}^n s_a = s,$ determine $\san$ completely. It follows that $\san$ is a particular label for the $n$-core (cf. \cite{Mac}) of $\lambda.$ 
\medskip     
\\ \noindent (ii) In a similar manner one can describe the pair $(\lal, \sal)$ corresponding to $\lambda$ and $s.$ In this case we first associate with $\lambda$ and $s$ another partition which we denote by $\sigma^s(\lambda).$ 
Let $\kb = (k_1,k_2,\dots\:)$ be the element of $\PP^{++}(s)$ defined by $(\lambda,s).$ For each $i\in \N$ we  put $k_i = e_i - nlm_i$ where $e_i \in \{1,\dots,nl\},$ $m_i \in \Z.$ 
Next, for every $e \in \{1,\dots, nl\}$ we define $(k^{(e)}_1,k_2^{(e)},\dots \:)$ to be the semi-infinite strictly decreasing  sequence $( 1 - m_i\: | \: e_i = e).$ This sequence stabilizes, for $i\gg 1,$  to  $t_e - i + 1$ where $t_e$ is a uniquely defined integer. Then, we define a partition $\mu^{(e)} = (\mu^{(e)}_1,\mu^{(e)}_2,\dots\:)$ by $\mu^{(e)}_i = k^{(e)}_i - t_e + i -1.$ In this way we get an $nl$-multipartition $(\mu^{(1)},\dots, \mu^{(nl)})$ and a sequence of $nl$ integers $(t_1,\dots, t_{nl})$ summing to $s.$ Of course, these are nothing but the $nl$-quotient and the $nl$-core of $\lambda.$ Let $\sigma$ be the permutation of $\{1,\dots,nl\}$ defined by  
$$ 
\sigma : d + l(c-1) \mapsto c + n(d-1)\qquad ( c \in \{1,\dots,n\},\; d \in \{1,\dots,l\}).
$$ 
Now, $\sigma^s(\lambda)$ is defined to be the unique partition with the $nl$-quotient $(\mu^{(\sigma(1))},\dots,\mu^{(\sigma(nl))}),$ and the $nl$-core $(t_{\sigma(1)},\dots,t_{\sigma(nl)}).$ Finally, $\lal$ and $\sal$ can be read off the  diagram of $\sigma^s(\lambda)$ in exactly the same way as $(\lan,\san)$ are read off the  diagram of $\lambda,$ i.e. as explained in (i) above, but reversing everywhere the roles of $n$ and $l.$   
}
\end{remark}

\begin{example} {\rm The following table illustrates, for  partitions $\lambda$ of size $7,$ the three indexations defined above. Here $n=3,$ $l=2,$ $s=0.$ 
\nopagebreak 
$$
\begin{array}{|l| l | l |  l | l | l |}  \hline 
   \lambda &  \sigma^s(\lambda) & \lal & \sal & \lan & \san   \\  \hline 
  ( 7)  &  (7) & ( ( 3) ,\varnothing )  & ( 1,-1)  & ( ( 2) ,\varnothing ,\varnothing )  & ( 1,0,  -1) \\
  ( 6,1)  & (6,1^2) & ( \varnothing ,( 3,1) )  & ( 0, 0)  & ( \varnothing ,\varnothing ,( 1) )  & ( 0,-1, 1)  \\
  ( 5,2)  & (4,2,1)& ( \varnothing ,( 3) )  & ( 1, -1)   &   ( \varnothing ,( 2) ,\varnothing )  & ( 1,0, -1)  \\
  ( 5,1^2)  & (4,1^4)&  ( \varnothing ,( 2,1^2) )  & ( 0, 0)  & ( \varnothing ,( 1) ,\varnothing )  & ( -1,1, 0)  \\
  ( 4,3)  & (3^2,1) & ( ( 1) ,( 2) )  & ( 1, -1)  &  ( ( 1) ,( 1) ,\varnothing )  & ( 1,0,-1)  \\
  ( 4,2,1)  &  (2^2,1^3) & ( \varnothing ,( 2,1) )  & ( 1,-1) &  ( ( 1^2) ,\varnothing ,\varnothing )  & ( 1,0, -1) \\  
  ( 4,1^3)  & (2,1) & ( \varnothing ,\varnothing )  & ( -1,1)  &  ( ( 1) ,\varnothing ,( 1) )  & ( 1,0, -1) \\
  ( 3^2,1)  & (5,4,1^2) & ( ( 1^2) ,( 1) )  & ( 2,-2) & ( \varnothing ,\varnothing ,( 1) )  & ( -1,1, 0)   \\
  ( 3,2^2)  & (5,2^3,1) & ( ( 1) ,( 2) )  & ( 2,-2)  & ( ( 1) ,\varnothing ,\varnothing )  & ( 0,-1, 1) \\ 
  ( 3,2,1^2)  & (5,2) & ( ( 2,1) ,\varnothing )  & ( 1,-1) &  ( \varnothing ,\varnothing ,( 2) )  & ( 1,0, -1) \\
  ( 3,1^4)  & (5,1^3) & ( ( 3,1) ,\varnothing )  & ( 0, 0)  &  ( \varnothing ,( 1) ,\varnothing )  & ( 0,-1, 1) \\ 
  ( 2^3,1)  & (3,2^2) & ( ( 1^2) ,( 1) )  & ( 1,-1)  &  ( \varnothing ,( 1) ,( 1) )  & ( 1,0,-1)  \\
  ( 2^2,1^3)  & (3,2,1^2) & ( ( 1^3) ,\varnothing )  & ( 1,-1) &  ( \varnothing ,( 1^2) ,\varnothing )  & ( 1,0, -1)  \\
  ( 2,1^5)  & (3,1^5) & ( ( 2,1^2) , \varnothing )  & ( 0, 0) &  ( ( 1) ,\varnothing , \varnothing )  & ( -1,1, 0) \\ 
  ( 1^7)  & (1^7) & ( \varnothing ,( 1^3) )  & ( 1,-1)  &  ( \varnothing ,\varnothing ,( 1^2) )  & ( 1,0,-1) \\  \hline   
\end{array}  
$$
}
\nopagebreak
\finex \end{example}

\subsection{Actions of quantum affine algebras on $\SIW.$} \label{ss:ACTIONS}
Taking the action (\ref{e:en} -- \ref{e:dn}) of $\UN,$ and the action  (\ref{e:el} -- \ref{e:dl}) of $\UL$ as the input, we shall define actions of $\UN$ and $\UL$ on the semi-infinite wedge product $\SIW.$ First, we assign to each vector $u$ of $\Lambda^{*+\frac{\infty}{2}}= \oplus_{s\in\Z}\SIW$ a weight $\wt(u)$ of $\asl_n,$ and a weight $\wtd(v)$ of $\asl_l.$ We shall write $\Wt(u)$ for the sum of $\wt(u)$ and $\wtd(u).$ For $s\in \Z,$ set $|s\rangle := u_s\wedge u_{s-1}\wedge \cdots\:.$ Let $\wt(|0\rangle) := l\Lambda_0,$ and  $\wtd(|0\rangle) := n\dLambda_0.$ Here, and in what follows we put dots over the fundamental weights, fundamental roots and the null root of $\asl_l$ in order to distinguish them from those of $\asl_n.$ For each non-zero integer $s$ let   
$$
\Wt(|s\rangle):= \Wt(|0\rangle) + \begin{cases} -\Wt(u_0\wedge u_{-1} \wedge \cdots \wedge u_{s+1}) & \text{ if $s<0,$ }\\ \;\;\:\Wt(u_s\wedge u_{s-1}\wedge \cdots \wedge u_1)& \text{ if $s>0.$}  \end{cases}
$$ 
Here $\Wt(u_0\wedge u_{-1} \wedge \cdots \wedge u_{s+1})$ and $\Wt(u_s\wedge u_{s-1}\wedge \cdots \wedge u_1)$ are defined by (\ref{e:tn}), (\ref{e:dn})  and (\ref{e:tl}), (\ref{e:dl}).
Then for $r\leq t,$ and $v\in \Lambda^r,$ the expression
$$ \Wt(v\wedge u_{s-r}\wedge u_{s-r-1}\wedge \cdots \wedge u_{s-t+1}) + \Wt(|s-t\rangle)$$   
is independent of the choice of $t.$ Hence, the assignment 
$$ \Wt(v\wedge |s-r ) : =  \Wt(v) + \Wt(|s-r\rangle)$$
gives a well-defined weight for the vector $v\wedge |s-r\rangle$ of $\SIW.$ Thus we have obtained a weight decomposition of $\SIW.$ It is straightforward to verify that in terms of the basis $\{ |\lal,\sal\rangle \: | \: \lal \in \Pi^l,\, \sal \in \Z^l(s)\}$ $=$ $\{ |\lan,\san\rangle \: | \: \lan \in \Pi^n,\, \san \in \Z^n(s)\}$ this decomposition looks as follows 
\begin{alignat}{2}
&\wt(|\lal,\sal\rangle)&= & -\Delta(\sal,n)\,\delta + \Lambda_{s_1}+\cdots+\Lambda_{s_l}-\sum_{i=0}^{n-1} M_i(\lal|\sal,n) \alpha_i, \label{e:carl}\\
&\wtd(|\lal,\sal\rangle)&= &  -(\Delta(\sal,n) + M_0(\lal|\sal,n))\,\ddelta + \label{e:cardl}\\
 & & & (n-s_1+s_l)\dLambda_0 + (s_1-s_2)\dLambda_1+\cdots +(s_{l-1}-s_l)\dLambda_{l-1},  \nonumber\\
&\wtd(|\lan,\san\rangle)&= & -\Delta(\san,l)\,\ddelta + \dLambda_{s_1}+\cdots+\dLambda_{s_n}-\sum_{i=0}^{l-1} M_i(\lan|\san,l) \dalpha_i, \label{e:cardn}\\
&\wt(|\lan,\san\rangle)&= &  -(\Delta(\san,l) + M_0(\lan |\san,l))\,\delta + \label{e:carn}\\
 & & & (l-s_1+s_n)\Lambda_0 + (s_1-s_2)\Lambda_1+\cdots +(s_{n-1}-s_n)\Lambda_{n-1}.  \nonumber
\end{alignat}
Now we define actions of the Cartan parts of $\UN$ and $\UL$ by  
\begin{alignat}{3}
&\tn_i\, u = q^{\langle \wt(u), h_i\rangle} u,\qquad & & \dn\,u = \langle \wt(u), \dn\rangle u && \qquad (i=0,1,\dots,n-1),\\
&\tl_j\, u = p^{\langle \wtd(u), \dot{h}_j\rangle} u,\qquad & & \dl\,u = \langle \wtd(u), \dl\rangle u && \qquad (j=0,1,\dots,l-1).
\end{alignat}
Next, we define on $\SIW$ actions of the raising generators of $\UN$ and $\UL.$ Let $v\in \Lambda^r$ for some $r\in \N.$ Lemma \ref{l:ZERO} (i) and (\ref{e:en}), (\ref{e:el}) imply that expressions 
\begin{alignat*}{2}
&\en_i(v\wedge u_{s-r}\wedge u_{s-r-1}\wedge \cdots \wedge u_{s-t+1})\wedge (\tn_i)^{-1}|s-t\rangle \qquad & & (i=0,\dots,n-1), \\
&\el_j(v\wedge u_{s-r}\wedge u_{s-r-1}\wedge \cdots \wedge u_{s-t+1})\wedge (\tl_j)^{-1}|s-t\rangle \qquad & & (j=0,\dots,l-1). 
\end{alignat*}
are independent of $t$ for $t\geq r.$ Hence, the assignments  
\begin{alignat*}{2}
&\en_i(v\wedge|s-r\rangle) := \en_i(v)\wedge (\tn_i)^{-1}|s-r\rangle  \qquad & & (i=0,\dots,n-1), \\
&\el_j(v\wedge|s-r\rangle) := \el_j(v)\wedge (\tl_j)^{-1}|s-r\rangle  \qquad & & (j=0,\dots,l-1)
\end{alignat*}
determine well-defined endomorphisms of $\SIW.$ 

Finally, we  define actions of the lowering generators of $\UN$ and $\UL$ on $\SIW.$ Using now Lemma \ref{l:ZERO} (ii) and  (\ref{e:fn}), (\ref{e:fl}) one can check, that expressions 
\begin{alignat*}{2}
&\fn_i(v\wedge u_{s-r}\wedge u_{s-r-1}\wedge \cdots \wedge u_{s-t+1})\wedge u_{s-t}\wedge u_{s-t+1}\wedge \cdots \wedge u_{s-m+1}  \quad & & (i=0,\dots,n-1), \\
&\fl_j(v\wedge u_{s-r}\wedge u_{s-r-1}\wedge \cdots \wedge u_{s-t+1})\wedge u_{s-t}\wedge u_{s-t+1}\wedge \cdots \wedge u_{s-m+1}  \quad & & (j=0,\dots,l-1) 
\end{alignat*}
are independent of $t\leq m$ provided $t>r$ and $t$ is sufficiently large. Hence, we obtain well-defined endomorphisms of $\SIW$ by setting 
\begin{alignat*}{2}
&\fn_i(v\wedge|s-r\rangle) := \fn_i(v\wedge u_{s-r}\wedge u_{s-r-1}\wedge \cdots \wedge u_{s-t+1})\wedge |s-t\rangle   \qquad & & (i=0,\dots,n-1), \\
&\fl_j(v\wedge|s-r\rangle) := \fl_j(v\wedge u_{s-r}\wedge u_{s-r-1}\wedge \cdots \wedge u_{s-t+1})\wedge |s-t\rangle   \qquad & & (j=0,\dots,l-1)
\end{alignat*}
where $t$ is arbitrary such that $t\gg r.$

The actions of $\fn_i$ and $\en_i$ defined above are easily  described in terms of the basis $B(s)$ = $\{ |\lal,\sal\rangle \: | \: \lal \in \Pi^l,\, \sal \in \Z^l(s)\}.$ Using notations of Section  \ref{ss:DEF} we have 

\begin{alignat}{2} 
&\fn_i |\lal,\sal\rangle  = \sum_{{\mathrm {res}}_n(\mul/\lal) = i} q^{\; N_i^{>}(\lal,\,\mul |\sal,n)} \: & & |\mul,\sal\rangle, \label{e:fnmult}\\
&\en_i |\mul,\sal\rangle  = \sum_{{\mathrm {res}}_n(\mul/\lal) = i} q^{- N_i^{<}(\lal,\,\mul |\sal,n)} \: & & |\lal,\sal\rangle.\label{e:enmult}
\end{alignat}
These actions, and the actions of the generators $\tn_i,$ $\dn$ defined by (\ref{e:carl}), are identical with those defined on the combinatorial Fock space $\FN.$ It follows that $\en_i,\fn_i,\tn_i$ and $\dn$ satisfy the defining relations of $\UN.$ Moreover, we see that for each $\sal = (s_1,\dots,s_l) \in \Z^l,$ the Fock space $\FN$  is realized inside $\SIW$ with $s= s_1+\cdots+s_l$ as the subspace spanned by $\{ |\lal,\sal\rangle\:|\:\lal \in \Pi^l\}.$ Note that by (\ref{e:cardl}), the  Fock space $\FN$ is the set of all vectors $u \in \SIW$ such that $\wtd(u)$ is congruent to $(n-s_1+s_l)\dLambda_0 + (s_1-s_2)\dLambda_1+\cdots +(s_{l-1}-s_l)\dLambda_{l-1}$ modulo $\Z\ddelta.$      

Similarly, one can describe the actions of $\fl_j$ and $\el_j$ in terms of the basis $B(s).$ Now the formulas acquire a simple form if we use the indexation of $B(s)$ by $(\lan,\san)$ $(\lan \in \Pi^n, \san \in \Z^n(s)).$ We have   
\begin{alignat}{2}
&\fl_j |\lan,\san\rangle  = \sum_{{\mathrm {res}}_l(\mun/\lan) = j} p^{\; N_j^{>}(\lan,\,\mun |\san,l)} \: & & |\mun,\san\rangle, \label{e:flmult}\\
&\el_j |\mun,\san\rangle  = \sum_{{\mathrm {res}}_l(\mun/\lan) = j} p^{- N_j^{<}(\lan,\,\mun |\san,l)} \: & & |\lan,\san\rangle.\label{e:elmult}
\end{alignat}
These actions, and the actions of the generators $\tl_j,$ $\dl$ defined by (\ref{e:cardn}), are identical with those defined on the combinatorial Fock space $\FL.$ It follows that $\el_j,\fl_j,\tl_j$ and $\dl$ satisfy the defining relations of $\UL.$ Moreover, for each $\san = (s_1,\dots,s_n) \in \Z^n,$ the Fock space $\FL$  is realized inside $\SIW$ with $s= s_1+\cdots+s_n$ as the subspace spanned by $\{ |\lan,\san\rangle\:|\:\lan \in \Pi^n\}.$ Note that by (\ref{e:carn}), the  Fock space $\FL$ is the set of all vectors $u \in \SIW$ such that $\wt(u)$ is congruent to $(l-s_1+s_n)\Lambda_0 + (s_1-s_2)\Lambda_1+\cdots +(s_{n-1}-s_n)\Lambda_{n-1}$ modulo $\Z\delta.$ 

\subsection{Action of the Heisenberg algebra on $\SIW.$}
Let $v\in \Lambda^r$ for some $r\in \N.$ Lemma \ref{l:ZERO} (i)  and (\ref{e:Bo}) imply that for $m>0$ the expression 
$$
B_m(v\wedge u_{s-r}\wedge u_{s-r-1}\wedge \cdots \wedge u_{s-t+1})\wedge |s-t\rangle 
$$
is  independent of $t$ for $t\geq r.$ Hence, the assignment  
$$
B_m(v\wedge|s-r\rangle) := B_m(v)\wedge |s-r\rangle \qquad (m>0)
$$ 
determines a  well-defined endomorphism of $\SIW.$ 
Using now Lemma \ref{l:ZERO} (ii) one can check that for $m<0$ the expression
$$
B_m(v\wedge u_{s-r}\wedge u_{s-r-1}\wedge \cdots \wedge u_{s-t+1})\wedge u_{s-t}\wedge u_{s-t+1}\wedge \cdots \wedge u_{s-k+1} 
$$
is  independent of $t\leq k$ provided $t>r$ and $t$ is sufficiently large. Hence, we obtain a well-defined endomorphism of $\SIW$ by setting 
$$
B_m(v\wedge|s-r\rangle) := B_m(v\wedge u_{s-r}\wedge u_{s-r-1}\wedge \cdots \wedge u_{s-t+1})\wedge |s-t\rangle   \qquad (m<0)
$$
where $t$ is arbitrary such that $t\gg r.$ It is clear that $B_m$ $(m\in \N)$ has weight $m\delta  + m\ddelta$ with respect to the weight decomposition of $\SIW$ defined in the previous section. This implies that the subspaces $\FN$ and $\FL$ are preserved by $B_m.$  The next proposition shows that $B_m$ generate a Heisenberg algebra $\He.$

\begin{proposition} 
There are non-zero  $\gamma_m \in \K$ {\em (}independent of $s${\em)} such that 
$$ [B_m,B_{m'}] = \delta_{m+m',0} \gamma_m. $$ 
\end{proposition}
\proof Since for each $r$ the actions of $B_m$ on $\Lambda^r$ commute,  for any $v \in \Lambda^r$ we have 
\begin{equation} 
[B_m,B_{m'}](v\wedge |s-r\rangle) = v\wedge[B_m,B_{m'}]|s-r\rangle. \label{e:bbb1}
\end{equation}
It is, therefore, enough to show that the statement of the proposition holds when $[B_m,B_{m'}]$ is applied to $|s\rangle$ with arbitrary $s\in \Z.$ Let us first assume that $m+m'>0.$ Then it is easy to see from (\ref{e:carl} -- \ref{e:carn}) that $\Wt(|s\rangle ) + (m+m')(\delta + \ddelta)$ is not a weight of $\SIW.$ Hence, $[B_m,B_{m'}]|s\rangle =0$ in this case. Next, let $m+m'<0.$ Write $[B_m,B_{m'}]|s\rangle$ as a linear combination of ordered wedges
$$
[B_m,B_{m'}]|s\rangle = \sum_{\nu} c_{\nu} u_{k_1^{\nu}}\wedge u_{k_2^{\nu}}\wedge \cdots  
$$ 
where $c_{\nu}$ are non-zero coefficients. Since $\Wt(|s\rangle) + (m+m')(\delta + \ddelta)$ is distinct from the weight of $|s\rangle,$ we have $k_1^{\nu} > s$ for all $\nu.$ For any $t>0,$ (\ref{e:bbb1}) gives   
$$ 
[B_m,B_{m'}]|s\rangle  = u_s\wedge u_{s-1} \wedge \cdots \wedge u_{s-nlt+1}\wedge [B_m,B_{m'}]|s-nlt\rangle. 
$$ 
From the structure of the ordering rules of Proposition \ref{p:RULES} it follows that 
$$
[B_m,B_{m'}]|s-nlt\rangle = \sum_{\nu} c_{\nu} u_{k_1^{\nu}-nlt}\wedge u_{k_2^{\nu}-nlt}\wedge \cdots\quad.  
$$ 
Choosing  $t$ sufficiently large, so that $s \geq k_1^{\nu} - nlt$ holds for all $\nu,$ and taking into account the inequality $k_1^{\nu} - nlt > s - nlt,$ we obtain, by Lemma \ref{l:ZERO} (ii), that $u_s\wedge u_{s-1} \wedge \cdots \wedge u_{s-nlt+1}\wedge u_{k_1^{\nu}} \wedge u_{k_2^{\nu}} \wedge \cdots $ vanishes for all $\nu.$ Hence  $[B_m,B_{m'}]|s\rangle $ is zero. 

Finally, let $m+m'=0.$ Then the weight of $[B_m,B_{m'}]|s\rangle$ equals $\Wt(|s\rangle).$ It is easy to see that  the vector $|s\rangle$ coincides with $|\ol,\sal\rangle $ for a certain $\sal \in \Z^l(s).$  But from (\ref{e:carl} -- \ref{e:carn}) it is clear that the weight subspace of $\FN$ with the weight $\Wt(|\ol,\sal\rangle)$ is one-dimensional. Since $B_m$ preserve  $\FN,$ it follows that $[B_m,B_{m'}]|s\rangle = \gamma_m |s\rangle$ for some $\gamma_m \in \K.$ Using $[B_m,B_{m'}]|s\rangle  = u_s\wedge [B_m,B_{m'}]|s-1\rangle$ shows that $\gamma_m$ does not depend on $s.$ Specializing to $q=1,$ we obtain  $\gamma_m|_{q=1} = mnl.$ Hence, $\gamma_m \neq 0.$ \finpf   

\begin{proposition} \label{p:gamma}
For $m >0$ we have 
$$ \gamma_m = m\frac{1-q^{-2mn}}{1-q^{-2m}} \frac{1-q^{2ml}}{1-q^{2m}}. $$ 
\end{proposition}
\noindent For $l=1$ a proof of this proposition is given in \cite{KMS}. We give a proof for all $l\in\N$ in Section \ref{sss:gamma}.

Following the same reasoning that  was used above to show that $B_m$ commutes with $B_{m'}$ unless $m+m'=0,$ we obtain  
\begin{proposition}
The actions of $\UNp,$ $\ULp$ and $\He$ on $\SIW$ are pairwise mutually commutative. 
\end{proposition}

Taking into account the weight decomposition of $\SIW$ given by (\ref{e:carl} -- \ref{e:carn}), and  the fact that weight of $B_m$ equals $m(\delta + \ddelta),$ we see that vectors $|\ol,\sal\rangle $ $(\sal \in \Z^l(s))$ are  singular vectors of $\He,$ \ie are annihilated by $B_m$ with positive $m.$ Clearly, these  vectors are also singular vectors for $\UN.$ Likewise, the vectors $|\on,\san\rangle $ $(\san \in \Z^n(s))$ are singular both for $\He$ and $\UL.$ Define sets $A_l^n(s)$ and $A_n^l(s)$ by 
\begin{alignat*}{2}
&A_l^n(s) := \{\sal = (s_1,\dots,s_l) \in \Z^l(s)\:& & |\: s_1\geq s_2\geq \cdots \geq s_l,\: s_1 - s_l \leq n\},\\
&A_n^l(s) := \{\san = (s_1,\dots,s_n) \in \Z^n(s)\:& & |\: s_1\geq s_2\geq \cdots \geq s_n,\: s_1 - s_n \leq l\}.
\end{alignat*}
Definitions of $|\lal,\sal\rangle$ and $|\lan,\san\rangle$ given in Section \ref{ss:SIW} imply 
\begin{equation}  \label{e:twosets}
\{ |\ol,\sal \rangle \: | \: \sal \in A_l^n(s)\}  = \{ |\on,\san\rangle \: | \: \san \in A_n^l(s)\}. 
\end{equation} 
Hence, $|\ol,\sal \rangle$ $(\sal \in A_l^n(s))$ or, equivalently, $|\on,\san \rangle$ $(\san \in A_n^l(s))$ are the only vectors of the basis $B(s)$ that are simultaneously singular for $\UN,$ $\UL$ and $\He.$ The equality (\ref{e:twosets}) shows that we have a bijection   $A_l^n(s) \rightarrow A_n^l(s)$ such that 
$$  \sal \mapsto \san \qquad \text{if and only if}\qquad  |\ol,\sal \rangle = |\on,\san\rangle .$$ This bijection is completely determined by comparing the weights of $ |\ol,\sal \rangle$ and $|\on,\san\rangle$ according to (\ref{e:carl} -- \ref{e:carn}). Namely, $(t_1,\dots,t_n) \in A_n^l(s)$ is the image of $(s_1,\dots,s_l) \in A_l^n(s)$ if and only if    
$$ \Lambda_{s_1}+\cdots + \Lambda_{s_l} = (l-t_1+t_n)\Lambda_0 + (t_1-t_2)\Lambda_1 + \cdots +(t_{n-1} - t_n)\Lambda_{n-1},$$ 
or, equivalently, if and only if 
$$ \dLambda_{t_1}+\cdots + \dLambda_{t_n} = (n-s_1+s_l)\dLambda_0 + (s_1-s_2)\dLambda_1 + \cdots +(s_{l-1} - s_l)\dLambda_{l-1}.$$ 
\begin{example}{\rm 
Let $n=5,$ $l=2$ and $s=11.$ Then $A_2^5(11)$ contains three elements: $(6,5),(7,4)$ and $(8,3).$ On the other hand, $A_5^2(11)$ is formed by  elements $$(3,2,2,2,2),(3,3,2,2,1),(3,3,3,1,1).$$ The bijective correspondence between $A_2^5(11)$ and $A_5^2(11)$ is given by  
\begin{align*}
& |\varnothing_2,(6,5)\rangle =  |\varnothing_5,(3,2,2,2,2)\rangle, && -(\delta+\ddelta) + \Lambda_0 + \Lambda_1 + 4\dLambda_0+\dLambda_1,\\ 
&|\varnothing_2,(7,4)\rangle =  |\varnothing_5,(3,3,2,2,1)\rangle, && -2(\delta+\ddelta) + \Lambda_2 + \Lambda_4 + 2\dLambda_0  + 3\dLambda_1,  \\
&|\varnothing_2,(8,3)\rangle =  |\varnothing_5,(3,3,3,1,1)\rangle,  && -3(\delta+\ddelta) + 2\Lambda_3 + 5\dLambda_1.  
\end{align*}
Here we listed weights of the corresponding vectors in the right column. 
}
\finex \end{example}
The next theorem shows that $\{ |\ol,\sal \rangle \: | \: \sal \in A_l^n(s)\}  = \{ |\on,\san\rangle \: | \: \san \in A_n^l(s)\}$ is the complete set of singular vectors in $\SIW.$ A proof follows immediately from \cite[Theorem 1.6]{F1} (see also \cite[Theorem 3.2]{F2}).
\begin{theorem} \label{t:TD}
$$ \SIW = \bigoplus_{ \sal \in A_l^n(s)} \textstyle{\UNp\cdot\He\cdot\ULp \, |\ol,\sal\rangle },$$
equivalently, 
$$ \SIW = \bigoplus_{ \san \in A_n^l(s)} \textstyle{\UNp\cdot\He\cdot\ULp \, |\on,\san\rangle }.$$

\end{theorem}

\begin{corollary} \label{c:ts} \mbox{} \\
{\em (i)} For ${\boldsymbol {t}}_l = (t_1,\dots,t_l) \in A_l^n(s),$ let $\sal = (s_1,\dots,s_l)$ be any element of $\Z^l(s)$ such that $\Lambda_{s_1}+\cdots + \Lambda_{s_l} = \Lambda_{t_1}+\cdots + \Lambda_{t_l}.$ Then $|\ol,{\boldsymbol {t}}_l\ket \in \ULp |\ol,{\boldsymbol {s}}_l\ket.$  

\noindent {\em (ii)} For ${\boldsymbol {t}}_n = (t_1,\dots,t_n) \in A_n^l(s),$ let $\san = (s_1,\dots,s_n)$ be any element of $\Z^n(s)$ such that $\dLambda_{s_1}+\cdots + \dLambda_{s_n} = \dLambda_{t_1}+\cdots + \dLambda_{t_n}.$ Then $|\on,{\boldsymbol {t}}_n\ket \in \UNp |\on,{\boldsymbol {s}}_n\ket.$
\end{corollary}

\proof Proofs of (i) and (ii) being almost identical, we show (i) only. Since $|\ol,{\boldsymbol {s}}_l\ket$ is a singular vector for $\UN\cdot \He,$ we have, by Theorem \ref{t:TD},   
$$ |\ol,{\boldsymbol {s}}_l\ket  \in \bigoplus_{ {\boldsymbol {r}}_l \in A_l^n(s)} \ULp\,|\ol,{\boldsymbol {r}}_l\ket.$$ 
Observe  that for two distinct elements  ${\boldsymbol {r}}_l$ and ${\boldsymbol {t}}_l $  of $A_l^n(s),$ we have $\Lambda_{r_1}+\cdots + \Lambda_{r_l} \neq \Lambda_{t_1}+\cdots + \Lambda_{t_l}.$ Therefore, comparing  $\asl_n$-weights, we have
$$ |\ol,{\boldsymbol {s}}_l\ket  \in  \ULp\,|\ol,{\boldsymbol {t}}_l\ket.$$ 
Since  $\ULp\,|\ol,{\boldsymbol {t}}_l\ket$ is an irreducible representation of $\ULp,$ the claim follows.
\finpf

\subsection{Canonical bases of the $q$-deformed Fock spaces}
Fix an arbitrary integer $s$ and define a gradation of the semi-infinite wedge product $\SIW$  by setting $\deg |\lambda,s\rangle = |\lambda|.$  
\begin{lemma} \label{l:STAB}
Let $\kb = (k_1,k_2,\dots) \in \PP^{++}(s).$ Then for any $t,r \in \N$ such that $t > r \geq \deg u_{\kb}$ we have 
$$ \ov{u_{k_1}\wedge u_{k_2}\wedge \cdots \wedge u_{k_r}}\wedge u_{k_{r+1}}\wedge \cdots \wedge u_{k_t} =\ov{u_{k_1}\wedge u_{k_2}\wedge \cdots \wedge u_{k_r} \wedge u_{k_{r+1}}\wedge \cdots \wedge u_{k_t}}. $$   
\end{lemma}
\noindent For $l=1$ a proof of this lemma is given in \cite[Proof of Lemma 7.7]{LT2}, for arbitrary $l$ a proof is virtually identical and will be omitted. 

From this lemma it follows that for $\kb =(k_1,k_2,\dots) \in \PP^{++}(s),$ the assignment 
\begin{equation}\label{e:siInvo}
\ov{u_{\kb}} := \ov{u_{k_1}\wedge u_{k_2}\wedge \cdots \wedge u_{k_r}}\wedge u_{k_{r+1}}\wedge u_{k_{r+2}}\wedge\cdots  \qquad (r \geq \deg u_{\kb})
\end{equation}
determines a well-defined semi-linear involution $u \mapsto \ov{u}$ of $\SIW.$ It is easily seen from the weight decomposition of $\SIW$ defined  in Section \ref{ss:ACTIONS} that $\Wt(\ov{u}) = \Wt(u)$ for any weight vector $u$ of $\SIW.$ Hence, for $\sal \in \Z^l(s)$ (resp. $\san \in \Z^n(s)$), the Fock space $\FN$ (resp. $\FL$) is invariant with respect to the bar-involution. In particular, therefore, we have   
$$ \ov{|\lal,\sal\rangle } = \sum_{\mul \in \Pi^l} R_{\lal,\mul}(\sal | q) \, |\mul, \sal \rangle, $$
where $R_{\lal,\mul}(\sal | q)$ is a Laurent polynomial in $q$ with integral coefficients. From (\ref{e:carl}) and the fact that the involution preserves the weight subspaces of $\SIW,$ it follows that $R_{\lal,\mul}(\sal | q)$ is non-zero only if $|\lal| = |\mul|.$     

For $\lal,\mul \in \Pi^l$ and $\sal \in \Z^l(s),$ let $\kb = (k_1,k_2,\dots),$ and $\lb = (l_1,l_2,\dots )$ be the unique elements of $\PP^{++}(s)$ such that $|\lal,\sal \rangle = u_{\kb},$ $|\mul,\sal\rangle = u_{\lb}.$ Then (\ref{e:siInvo}) implies    
\begin{equation} \label{e:R=R}
R_{\lal,\mul}(\sal | q) = R_{(\kb)_r , (\lb)_r}(q),
\end{equation}
where $(\kb)_r := (k_1,k_2,\dots,k_r),$ $(\lb)_r := (l_1,l_2,\dots,l_r)$ and $r$ is an arbitrary integer satisfying $r \geq \deg u_{\kb},$ $\deg u_{\lb}.$ Here, the coefficient $R_{(\kb)_r , (\lb)_r}(q)$ is defined in Section \ref{ss:CBfin}. The unitriangularity of the matrix $\| R_{\kb,\lb}(q) \|$ $(\kb,\lb \in \P^{++})$ described in that section immediately leads to    
\begin{proposition} For $\lal,\mul \in \Pi^l$ and $\sal \in \Z^l(s),$ the coefficient $R_{\lal,\mul}(\sal | q)$ is zero unless  the partition $\lambda = (\tau_l^s)^{-1}(\lal,\sal)$ is greater or equals the partition $\mu =  (\tau_l^s)^{-1}(\mul,\sal)$  with respect to the dominance order on partitions. Moreover, $R_{\lal,\lal}(\sal | q) =1.$

\end{proposition}

The unitriangularity of the involution matrix $\|R_{\lal,\mul}(\sal | q)\, \|$ allows to define canonical bases $\{ \Gc^+(\lal,\sal)\: | \: \lal \in \Pi^l\},$ $\{ \Gc^-(\lal,\sal)\: | \: \lal \in \Pi^l\}$ of the Fock space $\FN$ for arbitrary $\sal \in \Z^l(s).$ These bases are characterized by   
\begin{align}
&\textstyle{\ov{\Gc^+(\lal,\sal)} = \Gc^+(\lal,\sal),}\qquad \textstyle{\ov{\Gc^-(\lal,\sal)} = \Gc^-(\lal,\sal)}, \tag*{(i)} \\ 
&\Gc^+(\lal,\sal) \equiv |\lal,\sal\rangle \bmod q \Lc^+(s),\qquad \Gc^-(\lal,\sal) \equiv |\lal,\sal\rangle \bmod q^{-1} \Lc^-(s), \tag*{(ii)}
\end{align}
where $\Lc^+(s)$ (resp. $\Lc^-(s)$) is the $\Q[q]$-lattice   (resp. $\Q[q^{-1}]$-lattice) of $\SIW$ generated by the basis 
$$B(s) = \{ |\lambda,s\rangle \: | \: \lambda \in \Pi\}=\{ |\lal,\sal\rangle \: |\: \lal \in \Pi^l,\, \sal \in \Z^l(s)\}= \{ |\lan,\san\rangle \: |\: \lan \in \Pi^n,\, \san \in \Z^n(s)\}.$$ Set     
$$ \Gc^+(\lal,\sal) = \sum_{\mul \in \Pi^l} \dplus_{\lal,\mul}(\sal|q)\, |\mul,\sal\rangle,\qquad \Gc^-(\lal,\sal) = \sum_{\mul \in \Pi^l} \dminus_{\lal,\mul}(\sal|q)\, |\mul,\sal\rangle.$$ Then, keeping notations as in (\ref{e:R=R}) we have 
$$ \dplus_{\lal,\mul}(\sal|q) = \dplus_{(\kb)_r,(\lb)_r}(q),\qquad \dminus_{\lal,\mul}(\sal|q) = \dminus_{(\kb)_r,(\lb)_r}(q),$$
where the matrices $\|\Delta^{\pm}_{\kb,\lb}(q)\|$ $(\kb,\lb \in \P^{++})$ are defined in Section \ref{ss:CBfin}. Hence, Theorem \ref{t:KL} shows that $\Delta^{\pm}_{\lal,\mul}(\sal|q)$ are  parabolic Kazhdan-Lusztig polynomials. Note that $R_{\lal,\mul}(\sal | q) \neq 0$ only if $|\lal|=|\mul|$ implies  $\Delta^{\pm}_{\lal,\mul}(\sal|q)\neq 0$  only if $|\lal|=|\mul|.$ For each non-negative integer $k$ let us put 
$$ \|\Delta^{\pm}_{\lal,\mul}(\sal|q)\|_k = \|\Delta^{\pm}_{\lal,\mul}(\sal|q)\|\quad (|\lal| = |\mul| =k). $$ 
A proof of the next proposition in the special case $l=1$ is given in \cite{LT2}, a proof of the general case is similar and will be omitted.  
\begin{proposition}  \label{p:BAR}
For each $u \in \SIW$ one has 
\begin{alignat*}{3}
&\ov{\en_i u} = \en_i\ov{u},\qquad & &\ov{\fn_i u} =  \fn_i\ov{u}\qquad & &(i\in \{0,1,\dots,n-1\}), \\ 
&\ov{\el_j u} = \el_j\ov{u},\qquad & &\ov{\fl_j u} = \fl_j\ov{u}\qquad & &(j\in \{0,1,\dots,l-1\}), \\ 
&\ov{B_{-m} u} =  B_{-m} \ov{u},\qquad & &\ov{B_m u} = q^{2m(n-l)} B_m\ov{u} \qquad & &( m >0).
\end{alignat*}
\end{proposition}

Let us now describe how the canonical bases relate to the global crystal bases of Kashiwara.  As in Section \ref{ss:CRYSTAL}, let $\Lc[\sal] = \oplus_{\lal \in \Pi^l} A \:|\lal,\sal\ket$ be the lower crystal lattice of $\FN$ at $q=0.$ Proposition \ref{p:BAR} then implies that $\ov{\Lc[\sal]} = \oplus_{\lal \in \Pi^l} \ov{A}\:\ov{|\lal,\sal\ket},$ where $\ov{A} \:\subset \Q(q)$ is the subring of rational functions regular at $q=\infty,$ is a lower crystal lattice of $\FN$ at $q=\infty$ (cf. \cite{Ka1,Ka2}). Let ${\UN}_{\Q}$ be the $\Q[q,q^{-1}]$-subalgebra of $\UN$ generated by the $q$-divided differences $\en^{(m)}_i, \fn^{(m)}_i$ and 
$$ \prod_{k=1}^m \frac{ q^{1-k} \tn_i -  (\tn_i)^{-1}q^{k-1}}{q^k - q^{-k}}$$ 
with $m\in \N.$ One can show \cite[Lemma 2.7]{AM} that 
$$ {\FN}_{\Q} = \oplus_{\lal \in \Pi^l} \Q[q,q^{-1}] \:|\lal,\sal\ket$$
is invariant with respect to the action of ${\UN}_{\Q}$ on $\FN.$ 

The existence and uniqueness of the basis $\{ \Gc^+(\lal,\sal)\:|\: \lal \in \Pi^l\}$ can be reformulated, by using the unitriangularity of the bar-involution, as the existence of an isomorphism    
$$ {\FN}_{\Q}\cap \Lc[\sal] \cap \ov{\Lc[\sal]} \;\overset{\sim}{\rightarrow} \; \Lc[\sal]/q\Lc[\sal] $$ 
such that the preimage of $|\lal,\sal \ket \bmod q\Lc[\sal]$ is $\Gc^+(\lal,\sal).$ In terminology of \cite{Ka1,Ka2}, therefore,  
$$\{ \Gc^+(\lal,\sal)\:|\: \lal \in \Pi^l\}$$
is a lower global crystal basis of the integrable $\UN$-module $\FN.$ 

Now let us use the indexation of the basis $B(s)$ by pairs $(\lan,\san)$ with $\lan \in \Pi^n$ and $\san \in \Z^n(s).$ Certainly, we may label by these pairs the canonical bases as well, so that 
$$ \{ \Gc^{\pm}(\lan,\san) \: | \: \lan \in \Pi^n,\, \san \in \Z^n(s)\} = \{ \Gc^{\pm}(\lal,\sal) \: | \: \lal \in \Pi^l,\, \sal \in \Z^l(s)\} $$ 
and $\Gc^{\pm}(\lan,\san)$ are congruent to $|\lan,\san\ket$ modulo $q^{\pm 1}\Lc^{\pm}(s).$ Comparing formulas (\ref{e:flmult}, \ref{e:elmult}) with (\ref{e:fnmult}, \ref{e:enmult}), and taking into account Theorem \ref{t:CRYSTAL}, we see that 
$\Lc[\san] = \oplus_{\lan \in \Pi^n} \ov{A} \:|\lan,\san\ket$ is a lower crystal lattice of the $\UL$-module $\FL$ at $p:=-q^{-1}=0.$ Then, by Proposition \ref{p:BAR} again,  $\ov{\Lc[\san]} = \oplus_{\lan \in \Pi^n} {A} \:\ov{|\lan,\san\ket}$ is a lower crystal lattice of $\FL$ at $p=\infty,$ and the existence and uniqueness of the basis $\{ \Gc^{-}(\lan,\san) \: | \: \lan \in \Pi^n\}$ imply that there is an isomorphism  
$$ {\FL}_{\Q}\cap \Lc[\san] \cap \ov{\Lc[\san]} \;\overset{\sim}{\rightarrow} \; \Lc[\san]/p\Lc[\san] $$ 
taking $\Gc^{-}(\lan,\san)$ to $|\lan,\san \ket \bmod p\Lc[\san].$ Therefore,  
$$\{ \Gc^-(\lan,\san)\:|\: \lan \in \Pi^n\}$$
is a lower global crystal basis of the integrable $\UL$-module $\FL.$

Let us now comment on how the canonical basis $\{ \Gc^+(\lal,\sal)\:|\: \lal \in \Pi^l\}$  is related to the lower global crystal basis $\{ \Gc(\lal,\sal)\:|\: \lal \in \Pi^l(\sal)\}$ of the irreducible $\UN$-submodule $\MN$ generated by $|\ol,\sal\ket$ (cf. Theorem \ref{t:GLOB}). 
Using Lemma 2.7 in \cite{AM}, one can show, that the rational form ${\MN}_{\Q}$ of $\MN$ belongs to ${\FN}_{\Q}.$ From the definition of $\Gc(\lal,\sal)$ it now follows that $\Gc(\lal,\sal)$ belongs to $\Lc^+(s),$ hence  has the same congruence property with respect to $\Lc^+(s)$ as does $\Gc^+(\lal,\sal).$ On the other hand, by Proposition \ref{p:BAR}, the restriction of the bar-involution on $\MN$ coincides with the involution of $\MN$ defined in Section \ref{ss:CRYSTAL}. 
By the uniquennes of  $\Gc^+(\lal,\sal)$ it now follows  that
  $\Gc^+(\lal,\sal) = \Gc(\lal,\sal)$
for all $\lal \in \Pi^l(\sal).$ Hence  
$$ \{ \Gc^+(\lal,\sal)\:|\: \lal \in \Pi^l(\sal)\} $$
is the lower global crystal basis of the irreducible $\UN$-submodule $\MN.$

For $\san \in \Z^n,$ let $\ML$ be the irreducible $\UL$-submodule of $\FL$ generated by the highest weight vector $|\on,\san\ket.$ By the same argument as above  we conclude  that 
$$ \{ \Gc^-(\lan,\san)\:|\: \lan \in \Pi^n(\san)\} $$
is the lower global crystal basis of $\ML.$

%Thus, the picture is completely symmetric under the exchange  $n \leftrightarrow l$ and $q \leftrightarrow p.$ 

Note that the involution matrix $R_{\lal,\mul}(\sal|q),$ because of its unitriangularity, can be computed by using the ordering rules of Proposition \ref{p:RULES}. We have, therefore, an algorithm for computation of the transition matrices $\| \Delta^{\pm}_{\lal,\mul}(\sal|q)\|.$ By the deep result of \cite{A}  the coefficients  $\Delta^{+}_{\lal,\mul}(\sal|1)$ for $\lal \in \Pi^l(\sal), \mul \in \Pi^l$ are identified with  the decomposition numbers  of   Specht modules for an Ariki-Koike algebra, hence are non-negative integers. 
Tables of  the transition matrices  suggest that for all  $\lal, \mul \in \Pi^l$ the entries $\Delta^{+}_{\lal,\mul}(\sal|q)$ are in $\Z_{\geq 0}[q]$ (and those of $\| \Delta^{-}_{\lal,\mul}(\sal|q)\|$  are in $\Z_{\geq 0}[p]$).

\section{An inversion formula for Kazhdan-Lusztig polynomials}
The aim of this section is to prove Theorem \ref{t:MAIN} which gives an inversion formula relating the matrices $\|\Delta^{+}_{\lal,\mul}(\sal|q)\|$ with $\|\Delta^{-}_{\lal,\mul}(\sal|q)\|.$ In the case $l=1$ this formula has already been proved by Leclerc and Thibon in \cite{LT2}.

\subsection{Some properties of the Heisenberg algebra action on $\SIW$}

\begin{definition}\label{d:dom}{\rm 
Let $m \in \Z_{\geq 0}.$ We shall say that a pair $({\boldsymbol{\lambda}}_r =(\lambda^{(1)},\dots,\lambda^{(r)}),$ ${\boldsymbol{s}}_r = (s_1,\dots,s_r))$ $\in$ $\Pi^r\times \Z^r$ is {\em $m$-dominant} if for all $a=1,2,\dots,r-1$ we have the inequalities 
$$ s_a - s_{a+1} \geq m + |{\boldsymbol{\lambda}}_r|, $$
}
were  $|{\boldsymbol{\lambda}}_r| = |\lambda^{(1)}| + \cdots + |\lambda^{(r)}|.$
\end{definition}

\noindent Also, we shall say that a basis vector $|\lal,\sal\rangle$ (resp. $|\lan,\san\rangle$) is $m$-dominant  if the pair $(\lal,\sal)$ (resp. $(\lan,\san)$) is.  To explain the reason for introducing this definition, we need to prepare some notations. Let $n \in \N,$ $l=1.$ Let $x$ be a linear operator on $\SIW$ acting on the elements of the basis $B(s)$ by 
$$ x\,|\lambda ,s \rangle = \sum_{\mu \in \Pi} x_{\lambda,\mu}(s) \, | \mu,s\rangle,$$ 
where $x_{\lambda,\mu}(s)$ are coefficients in $\K.$ Let now $n \in \N,$ $l\in \N.$ For each $b=1,2,\dots,l$ we define an endomorphism $x^{(b)}[n,1]$ of $\SIW$ by  
$$ 
x^{(b)}[n,1]\,|\lal,\sal\rangle = \sum_{\mu\in \Pi} x_{\lambda^{(b)},\mu}(s_b) \, | (\lambda^{(1)},\dots,\lambda^{(b-1)},\mu,\lambda^{(b+1)},\dots,\lambda^{(l)}), \sal \rangle. $$  
Similarly, let $n=1,$ $l\in \N.$ For an endomorphism  $y$ of $\SIW$ we introduce the corresponding matrix elements $y_{\lambda,\mu}(s)$ on the basis $B(s)$ by  
$$ y\,|\lambda ,s \rangle = \sum_{\mu \in \Pi} y_{\lambda,\mu}(s) \, | \mu,s\rangle.$$ 
Again, for arbitrary $n\in \N,$ $l\in \N$ we define, for each $a=1,2,\dots,n,$ an endomorphism $y^{(a)}[1,l]$ of $\SIW$ by  
$$ 
y^{(a)}[1,l]\,|\lan,\san\rangle = \sum_{\mu\in \Pi} y_{\lambda^{(a)},\mu}(s_a) \, | (\lambda^{(1)},\dots,\lambda^{(a-1)},\mu,\lambda^{(a+1)},\dots,\lambda^{(n)}), \san \rangle. $$   

\begin{example} \label{ex:SPLIT} {\rm 
For $n=2,$ $l=1,$ and $s\in 2\Z$ one finds, using the ordering rules of Proposition \ref{p:RULES}, that 
\begin{eqnarray*}
 B_{-2}|\varnothing,s\rangle &= &|(4),s\rangle - q^{-1} |(3,1),s\rangle + (q^{-2} - 1) |(2^2),s\rangle + \\
       &  &+q^{-1} |(2,1^2),s\rangle - q^{-2} | (1^4),s\rangle.
\end{eqnarray*}
Hence, for $n=2,$ $l=2,$ taking $\sal = (s_1,s_2)$ such that $s_1,s_2 \in 2\Z,$ we have 
\begin{eqnarray*}
B_{-2}^{(1)}[2,1]\,|\ol,\sal\rangle &= &|((4),\varnothing),\sal \rangle - q^{-1} |((3,1),\varnothing),\sal\rangle + (q^{-2} - 1) |((2^2),\varnothing),\sal\rangle + \\
       &  &+q^{-1} |((2,1^2),\varnothing),\sal\rangle - q^{-2} | ((1^4),\varnothing),\sal\rangle, \\  
B_{-2}^{(2)}[2,1]\,|\ol,\sal\rangle &= &|(\varnothing,(4)),\sal \rangle - q^{-1} |(\varnothing,(3,1)),\sal\rangle + (q^{-2} - 1) |(\varnothing,(2^2)),\sal\rangle + \\
       &  &+q^{-1} |(\varnothing,(2,1^2)),\sal\rangle - q^{-2} | (\varnothing,(1^4)),\sal\rangle.
\end{eqnarray*}
}
\finex \end{example}
\begin{proposition} \label{p:MAIN} \mbox{} \\
{\em (i)} Let $(\lal,\sal) \in \Pi^l\times \Z^l$ be $n m$-dominant for some $m\in \N.$ Then 
$$ B_{-m}\,|\lal,\sal\rangle  = \sum_{b=1}^l q^{(b-1)m} B_{-m}^{(b)}[n,1]\,|\lal,\sal\rangle.$$ 
{\em (ii)} Let $(\lal,\sal) \in \Pi^l\times \Z^l$ be $0$-dominant. Then, for any $m\in \N,$ 
$$ B_{m}\,|\lal,\sal\rangle  = \sum_{b=1}^l q^{(b-1)m} B_{m}^{(b)}[n,1]\,|\lal,\sal\rangle.$$ 
{\em (i')} Let $(\lan,\san) \in \Pi^n\times \Z^n$ be $l m$-dominant for some $m\in \N.$ Then 
$$ B_{-m}\,|\lan,\san\rangle  = \sum_{a=1}^n p^{(a-1)m} B_{-m}^{(a)}[1,l]\,|\lan,\san\rangle.$$ 
{\em (ii')} Let $(\lan,\san) \in \Pi^n\times \Z^n$ be $0$-dominant. Then, for any $m\in \N,$ 
$$ B_{m}\,|\lan,\san\rangle  = \sum_{a=1}^n p^{(a-1)m} B_{m}^{(a)}[1,l]\,|\lan,\san\rangle.$$ 
\end{proposition}
\noindent A proof of this proposition is given in Section \ref{ss:PROOF}. 
\begin{example}{\rm
To illustrate the above proposition, take $n=2,$ $l=2,$ and $\sal = (2,-2).$ A straightforward computation using Proposition \ref{p:RULES} gives  
\begin{eqnarray*}
B_{-2}\,|\ol,\sal\rangle &= &|((4),\varnothing),\sal \rangle - q^{-1} |((3,1),\varnothing),\sal\rangle + (q^{-2} - 1) |((2^2),\varnothing),\sal\rangle + \\
       &  &+q^{-1} |((2,1^2),\varnothing),\sal\rangle - q^{-2} | ((1^4),\varnothing),\sal\rangle +  \\  
 & & +q^2 |(\varnothing,(4)),\sal \rangle - q |(\varnothing,(3,1)),\sal\rangle + (1- q^2) |(\varnothing,(2^2)),\sal\rangle + \\
       &  &+q |(\varnothing,(2,1^2)),\sal\rangle - | (\varnothing,(1^4)),\sal\rangle.
\end{eqnarray*}
The pair $(\ol,\sal)$ is $4$-dominant,  taking into account formulas of Example \ref{ex:SPLIT}, we see that the relation 
$$ 
B_{-2}\,|\ol,\sal\rangle = B_{-2}^{(1)}[2,1]\,|\ol,\sal\rangle + q^2  B_{-2}^{(2)}[2,1]\,|\ol,\sal\rangle
$$
is indeed satisfied.
}
\finex \end{example}

\begin{remark}{\rm 
Simple decompositions for the actions of the bosons described in Proposition \ref{p:MAIN} fail to hold, in general, when bosons are applied to  vectors that are not dominant. For example, let $n=l=2,$ and let $\sal=(0,0).$ Then the pair $(\ol,\sal)$ is not $m$-dominant for any $m\in \N.$ In this case an explicit computation yields
\begin{eqnarray*}
B_{-2}\,|\ol,\sal\rangle &=&q |((4),\varnothing),\sal \rangle - |((3,1),\varnothing),\sal\rangle + (1-q^2) |((2^2),\varnothing),\sal\rangle + \\
       &  &+q |((2,1^2),\varnothing),\sal\rangle - | ((1^4),\varnothing),\sal\rangle +  \\  
 & & + |(\varnothing,(4)),\sal \rangle - q^{-1} |(\varnothing,(3,1)),\sal\rangle + (q^{-2}-1) |(\varnothing,(2^2)),\sal\rangle + \\
       &  &+ |(\varnothing,(2,1^2)),\sal\rangle - q^{-1}| (\varnothing,(1^4)),\sal\rangle + \\
& & + (q^2 -1) |((2),(2)),\sal\rangle + (q^{-1} - q) |((1),(2,1)),\sal \rangle + \\ && +(q^{-1} - q) |((2,1),(1)),\sal \rangle +  (1-q^{-2}) |((1^2),(1^2)),\sal \rangle . 
\end{eqnarray*}
}
\end{remark}

For $m\in \Z,$ let $e_m$ and $h_m$ be, respectively, the elementary symmetric function and the complete symmetric function (cf. \cite{Mac}). In terms of the power-sum basis of the ring of symmetric functions one has
$$ e_m =e_m(p_1,p_2,\dots) = \!\!\!\sum_{\nu \in \Pi, |\nu| = m} a_{m,\nu}\, p_{\nu},\quad h_m =h_m(p_1,p_2,\dots) = \!\!\!\sum_{\nu \in \Pi, |\nu| = m} b_{m,\nu}\, p_{\nu},$$ 
where, as usual, for a partition $\nu = (\nu_1,\nu_2,\dots),$ we put $p_{\nu} = p_{\nu_1}p_{\nu_2}\cdots.$ It will be understood that $e_0 = h_0 =1,$ and $e_m = h_m = 0$ for $m<0.$ Now, for all integer $m$ we define 
\begin{align*}
&\Ep_m:=e_m(B_1,B_2,\dots),& &\quad \Em_{m}:=e_m(B_{-1},B_{-2},\dots), \\
&\Hp_m:=h_m(B_1,B_2,\dots),& &\quad \Hm_{m}:=h_m(B_{-1},B_{-2},\dots).
\end{align*}
\begin{corollary} \label{c:MAIN} \mbox{} \\
{\em (i)} Let $(\lal,\sal) \in \Pi^l\times \Z^l$ be $n m$-dominant for some $m\in \N.$ Then 
\begin{align*}
&\Em_{m}\,|\lal,\sal\rangle = \sum_{m_1+\cdots + m_l = m}\; \prod_{b=1}^l q^{(b-1)m_b} \Em_{m_b}^{(b)}[n,1]\,|\lal,\sal\rangle,  \\
&\Hm_{m}\,|\lal,\sal\rangle  = \sum_{m_1+\cdots + m_l = m}\; \prod_{b=1}^l q^{(b-1)m_b} \Hm_{m_b}^{(b)}[n,1]\,|\lal,\sal\rangle. 
\end{align*}
{\em (ii)} Let $(\lal,\sal) \in \Pi^l\times \Z^l$ be $0$-dominant. Then, for any $m\in \N,$ 
\begin{align*}
&\Ep_{m}\,|\lal,\sal\rangle = \sum_{m_1+\cdots + m_l = m}\; \prod_{b=1}^l q^{(b-1)m_b} \Ep_{m_b}^{(b)}[n,1]\,|\lal,\sal\rangle,  \\
&\Hp_{m}\,|\lal,\sal\rangle  = \sum_{m_1+\cdots + m_l = m}\; \prod_{b=1}^l q^{(b-1)m_b} \Hp_{m_b}^{(b)}[n,1]\,|\lal,\sal\rangle. 
\end{align*}
{\em (i')} Let $(\lan,\san) \in \Pi^n\times \Z^n$ be $l m$-dominant for some $m\in \N.$ Then 
\begin{align*}
&\Em_{m}\,|\lan,\san\rangle = \sum_{m_1+\cdots + m_n = m}\; \prod_{a=1}^n p^{(a-1)m_a} \Em_{m_a}^{(a)}[1,l]\,|\lan,\san\rangle,  \\
&\Hm_{m}\,|\lan,\san\rangle  = \sum_{m_1+\cdots + m_n = m}\; \prod_{a=1}^n p^{(a-1)m_a} \Hm_{m_a}^{(a)}[1,l]\,|\lan,\san\rangle. 
\end{align*}
{\em (ii')} Let $(\lan,\san) \in \Pi^n\times \Z^n$ be $0$-dominant. Then, for any $m\in \N,$ 
\begin{align*}
&\Ep_{m}\,|\lan,\san\rangle = \sum_{m_1+\cdots + m_n = m}\; \prod_{a=1}^n p^{(a-1)m_a} \Ep_{m_a}^{(a)}[1,l]\,|\lan,\san\rangle,  \\
&\Hp_{m}\,|\lan,\san\rangle  = \sum_{m_1+\cdots + m_n = m}\; \prod_{a=1}^n p^{(a-1)m_a} \Hp_{m_a}^{(a)}[1,l]\,|\lan,\san\rangle. 
\end{align*}
\end{corollary}
\proof It follows from (\ref{e:carl} -- \ref{e:carn}) that for each $m \in \Z^*,$ a vector $B_{m} |\lal,\sal \rangle$ (resp. a vector $B_{m} |\lan,\san \rangle$) is a linear combination of $|\mul,\sal \rangle$ (resp. $|\mun,\san \rangle$)  with $|\mul| = |\lal| -nm$ (resp. with  $|\mun| = |\lan| -lm$). Also, if $\Delta$ is the comultiplication on the ring of symmetric functions defined by $\Delta p_m = p_m\otimes 1 + 1\otimes p_m,$ then $\Delta e_m = \sum_{r+s=m} e_r\otimes e_s,$ $\Delta h_m = \sum_{r+s=m} h_r\otimes h_s$ (cf. \cite{Mac}). These facts and Proposition \ref{p:MAIN} imply the claims. \finpf   

\subsubsection{Proof of Proposition \ref{p:gamma}.} \label{sss:gamma}

To emphasize the dependence on $n$ and $l$ let us put, in notations of Proposition \ref{p:gamma}, $\gamma_m[n,l] := \gamma_m.$ First, let $n=l=1.$ In this case the ordering rules of Proposition \ref{p:RULES} reduce to $u_{k_1}\wedge u_{k_2} = - u_{k_2}\wedge u_{k_1}$ for all $k_1,k_2 \in \Z.$ This makes it easy to verify that $\gamma_m[1,1] = m.$ Next, let $n > 1,$ $l=1.$ It is clear that there  is $\san \in \Z^n$ such that the pair $(\on,\san)$ is $m$-dominant. Hence, applying Proposition \ref{p:MAIN} (i'), we obtain   
$$ \gamma_m[n,1] |\on,\san\rangle = B_m B_{-m} |\on,\san\rangle = B_m \sum_{a=1}^{n} p^{(a-1)m} B_{-m}^{(a)}[1,1] |\on,\san\rangle.$$ 
But $B_{-m} |\on,\san\rangle $ is a linear combination of $|\lan,\san \rangle$ with $|\lan| = m,$ hence a linear combination of $0$-dominant vectors. Therefore one may apply Proposition \ref{p:MAIN} (ii') and get 
$$
B_m \sum_{a=1}^{n} p^{(a-1)m} B_{-m}^{(a)}[1,1] |\on,\san\rangle = \sum_{a=1}^{n} p^{2(a-1)m} B_m^{(a)}[1,1] B_{-m}^{(a)}[1,1] |\on,\san\rangle .
$$ 
This implies 
$$ \gamma_m[n,1] |\on,\san\rangle = \left(\sum_{a=1}^n p^{2(a-1)m} \gamma_m[1,1]\right) |\on,\san\rangle.$$ Thus 
$$\gamma_m[n,1] =  m \frac{1-p^{2mn}}{1-p^{2m}}.$$  
Finally, let $n \geq 1,$ $l>1.$ Obviously there is $\sal \in \Z^l$ such that the pair $(\ol,\sal)$ is $nm$-dominant. Therefore, by Proposition \ref{p:MAIN} (i), we have 
$$ \gamma_m[n,l] |\ol,\sal\rangle = B_m B_{-m} |\ol,\sal\rangle = B_m \sum_{b=1}^{n} q^{(b-1)m} B_{-m}^{(b)}[n,1] |\ol,\sal\rangle.$$ 
Again, it is clear that $B_{-m} |\ol,\sal\rangle$ is a linear combination of $0$-dominant vectors, hence using Proposition \ref{p:MAIN} (ii), we obtain    
$$ \gamma_m[n,l] |\ol,\sal\rangle = \left(\sum_{b=1}^l q^{2(b-1)m} \gamma_m[n,1]\right) |\ol,\sal\rangle.$$
It follows that 
$$ \gamma_m[n,l] = m \frac{1-p^{2mn}}{1-p^{2m}} \frac{1-q^{2ml}}{1-q^{2m}}.$$
Recalling that $p:=-q^{-1},$ we get the desired result. \finpf

\subsection{A scalar product of $\SIW.$}
For each $s\in \Z$ we define on the semi-infinite wedge product $\SIW$ a $\K$-bilinear scalar product by $\bra b , b'\ket = \delta_{b,b'},$ where $b$ and $b'$ are any two elements of the basis 
$$B(s) = \{ |\lambda,s\ket\:|\: \lambda \in \Pi\} = \{ |\lal,\sal \ket \:|\: \lal \in \Pi^l,\, \sal \in \Z^l(s)\}=\{ |\lan,\san \ket \:|\: \lan \in \Pi^n,\, \san \in \Z^n(s)\}.$$
It is clear that this scalar product is symmetric, and that for two weight vectors $u$ and $v,$ $\bra u , v\ket$ is non-zero only if $\Wt(u) = \Wt(v).$  

\begin{proposition} \label{p:eSCALARP}
For $u,v$ $\in$ $\SIW$ one has
\begin{alignat*}{3}
&\bra \en_i u , v\ket = \bra u , q^{-1} (\tn_i)^{-1} \fn_i v \ket,\quad & &\bra \fn_i u , v\ket = \bra u , q^{-1} \tn_i \en_i v \ket  \qquad & & (i=0,1,\dots,n-1), \\
&\bra \el_j u , v\ket = \bra u , p^{-1} (\tl_j)^{-1} \fl_j v \ket,\quad & & \bra \fl_j u , v\ket = \bra u , p^{-1} \tl_j \el_j v \ket  \qquad &  &(j=0,1,\dots,l-1).
\end{alignat*}
\end{proposition}
\proof  This immediately follows from (\ref{e:fnmult} -- \ref{e:elmult}). \finpf 

\begin{proposition}\label{p:bSCALARP}
For $m\in \N,$ and $u,v \in$ $\in$ $\SIW$ one has 
$$ \bra B_{-m}u , v\ket = \bra u , B_m v \ket .$$
\end{proposition}
\noindent To prove this proposition we use the following lemma.
\begin{lemma} \label{l:bSCALARP} \mbox{} \\
{\em (i)} Assume that the statement of {\em Proposition \ref{p:bSCALARP}} is valid for some $n\in \N$ and $l=1.$ Then it is also valid
for the same $n$ and all $l \in \N_{>1}.$ \\
{\em (ii)} Assume that the statement of {\em Proposition \ref{p:bSCALARP}} is valid for some $l\in \N$ and $n=1.$ Then it is also valid
for the same $l$ and all $n \in \N_{>1}.$
\end{lemma}

\proof  
Proofs of (i) and (ii) being similar, we give a proof of (i) only. Let us keep notations as in Proposition \ref{p:bSCALARP}. Using  Theorem \ref{t:TD} we assume without loss of generality that 
$$ \textstyle{u = \sum_k x_k y_k |\ol, {\boldsymbol{t}}_l \ket,} $$
where $x_k$ is an element of $\UN^-\cdot\, \He^-,$ $y_k$ is an element of $\UL^-,$ and ${\boldsymbol{t}}_l = (t_1,\dots,t_l)$ is an element of $ A_l^n(s).$ By Corollary \ref{c:ts}, for any $\sal =(s_1,\dots,s_l) \in \Z^l(s),$ such that $\Lambda_{s_1}+\cdots + \Lambda_{s_l} = \Lambda_{t_1}+\cdots + \Lambda_{t_l}$ we have $|\ol,{\boldsymbol{t}}_l\ket = Y(\sal)|\ol,\sal \ket$ for some $Y(\sal) \in \ULp.$ Hence 
$$ \textstyle{\bra B_{m} u , v\ket = \sum_k  \bra B_{-m} x_k |\ol,\sal \ket , Y_k^*(\sal) v\ket,}$$
where $Y_k^*(\sal)$ is the adjoint of $Y_k(\sal):= y_k Y(\sal).$ Note that by Proposition \ref{p:eSCALARP}, $Y_k^*(\sal)\in \ULp.$ We may and do assume that for each $k$ the element $x_k\in \UN^-\cdot \He^-$  has a  definite weight 
$ \wt(x_k) = -(r_{k,0}\alpha_0 + \cdots + r_{k,n-1}\alpha_{n-1}),$ where $r_{k,i}$ are some non-negative integers. For all $k$ we have $x_k |\ol,\sal \ket \in \FN,$ and using (\ref{e:carl}) we see that $x_k |\ol,\sal \ket$  is a linear combination of $|\lal,\sal \ket$ with $|\lal| = r_{k}:= r_{k,0} + \cdots + r_{k,n-1}.$ Let now  $r:= \max_k\{r_k\},$ and choose $\sal$ so that  
\begin{equation} 
 s_b - s_{b+1} \geq r + nm\qquad (b=1,\dots,l-1). \label{e:dom2}
\end{equation}
 Then for each $k$ the vector $x_k |\ol,\sal \ket$ is a linear combination of $nm$-dominant vectors. Hence we may apply Proposition \ref{p:MAIN} (i) and get 
$$ \bra B_{m} u , v\ket = \sum_k  \sum_{b=1}^l q^{(b-1)m}\bra B_{-m}^{(b)}[n,1] x_k |\ol,\sal \ket , Y_k^*(\sal) v\ket.$$
Now using the assumption in the statement (i) of the lemma we obtain
$$ \bra B_{m} u , v\ket = \sum_k  \sum_{b=1}^l  q^{(b-1)m} \bra x_k |\ol,\sal \ket , B_{m}^{(b)}[n,1]Y_k^*(\sal) v\ket.$$
For each $k$ the scalar product $\bra B_{-m} x_k |\ol,\sal \ket , Y_k^*(\sal) v\ket$ is non-zero only if $Y_k^*(\sal) v$ $\in$ $\FN,$ and 
$$ \wt(Y_k^*(\sal) v) = \wt(B_{-m} x_k |\ol,\sal \ket) = \wt|\ol,\sal \ket - \sum_{i=0}^{n-1} (r_{k,i} + m) \alpha_i.$$ It follows that $Y_k^*(\sal) v$ is a linear combination of $|\mul,\sal\ket$ with $|\mul| = r_k + nm,$ whence, by (\ref{e:dom2}), it is a linear combination of $0$-dominant vectors. Therefore, we may apply Proposition \ref{p:MAIN} (ii), and obtain  
$$\bra B_{-m}u , v\ket  = \sum_k \bra  x_k |\on,\san \rangle , B_m Y_k(\san)^* v\ket = \bra u , B_m v \ket.$$
\finpf 

\noindent{\em Proof of Proposition \ref{p:bSCALARP}.}
In view of Lemma \ref{l:bSCALARP} it is sufficient to show that the statement of the proposition is valid for $n=l=1.$ 
But in  this case the relation $\bra B_{-m}u , v\ket = \bra u , B_m v \ket $ is just a restatement of the fact that the endomorphism $m \partial/\partial p_m$ of the ring of symmetric functions is adjoint to the multiplication by $p_m$ with respect to the scalar product orthonormalizing the basis of Schur functions. \finpf

\noindent Proposition \ref{p:eSCALARP} implies that for $m \in \Z_{\geq 0}$ and $u,v$ $\in$ $\SIW$  one has 
\begin{equation}
\bra \Em_m u ,v \ket  =  \bra u , \Ep_m v\ket,\qquad \bra \Hm_m u ,v \ket  =  \bra u , \Hp_m v\ket. \label{e:ESCALARP}
\end{equation}

\subsection{A symmetry of the bar-involution.}
Define a semi-linear involution $u \mapsto u'$ of $\Lambda^{*+\frac{\infty}{2}} = \oplus_{s\in\Z}\SIW$ by 
$$ |\lambda , s\ket' = |\lambda' , -s\ket,\qquad q' = q^{-1}. $$
Here $\lambda'$ stands for the conjugate partition of $\lambda.$ The description of the indexations of $|\lambda , s\ket$ given in Remark \ref{r:QUO} implies that for an $l$-multipartition $\lal =(\lambda^{(1)},\dots,\lambda^{(l)})$ and $\sal = (s_1,\dots,s_l)$ $\in$ $\Z^l$ we have     
$$ |\lal, \sal \ket' = |\lal',\sal'\ket,$$ 
where $\lal' =({\lambda^{(l)}}',\dots,{\lambda^{(1)}}')$ and $\sal' = (-s_l,\dots,-s_1).$  Likewise, for an $n$-multipartition $\lan$ and $\san \in \Z^n$ we have  
$$ |\lan, \san \ket' = |\lan',\san'\ket.$$

\begin{proposition} \label{p:ePRIME}
For $u \in \SIW$ $(s\in \Z)$ we have 
\begin{alignat*}{3}
&(\en_i u)' = q^{-1}\tn_{-i} \en_{-i} u',\quad & & (\fn_i u)' = q^{-1}(\tn_{-i})^{-1} \fn_{-i} u'\qquad & &(i=0,1,\dots,n-1),\\ 
&(\el_j u)' = p^{-1}\tl_{-j} \el_{-j} u',\quad  & &(\fl_j u)' = p^{-1}(\tl_{-j})^{-1} \fl_{-j} u'\qquad & &(j=0,1,\dots,l-1).
\end{alignat*}
\end{proposition}
\proof This follows from (\ref{e:fnmult} -- \ref{e:elmult}). \finpf

\begin{proposition} \label{p:bPRIME}
For $u \in \SIW$ $(s\in \Z)$ and $m\in \Z_{\geq 0}$ we have 
\begin{align}
&(\Em_m u)' = (-q)^{m(n-1)} (-p)^{m(l-1)} \Hm_m u', \tag*{(i)}\\
&(\Ep_m u)' = (-q)^{m(n-1)} (-p)^{m(l-1)} \Hp_m u'.\tag*{(ii)}
\end{align}
\end{proposition}
\noindent To prove this proposition we use the following lemma.
\begin{lemma} \label{l:bPRIME} \mbox{} \\
{\em (i)} Assume that the statement {\em (i)} of {\em Proposition \ref{p:bPRIME}} is valid for some $n\in \N$ and $l=1.$ Then it is also valid
for the same $n$ and all $l \in \N_{>1}.$ \\
{\em (ii)} Assume that the statement {\em (i)} of {\em Proposition \ref{p:bPRIME}} is valid for some $l\in \N$ and $n=1.$ Then it is also valid
for the same $l$ and all $n \in \N_{>1}.$
\end{lemma}

\proof Since proofs of (i) and (ii) are virtually identical, we give a proof of (i) only. Let us keep notations as in Proposition \ref{p:bPRIME}. Taking into account Theorem \ref{t:TD} we may, and do, assume without loss of generality that 
$$ \textstyle{u = \sum_k x_k y_k |\ol ,\tal \ket,}$$  
where $\tal =(t_1,\dots,t_l)$ is an element of $A_l^n(s),$ $x_k$ is an element of $\UN^-\cdot \,\He^-,$ and $y_k$ is an element of $\UL^-.$ Choose any sequence $\sal =(s_1,\dots,s_l) \in \Z^l(s)$  such that the relation $\Lambda_{s_1} + \cdots + \Lambda_{s_l} =  \Lambda_{t_1} + \cdots + \Lambda_{t_l}$ is satisfied.  Corollary \ref{c:ts} (i) implies that  there is $Y(\sal) \in \ULp$ such that $|\ol ,\tal \ket = Y(\sal)|\ol ,\sal \ket.$ Setting $Y_k(\sal)=y_k Y(\sal)$ we have   
$$ \textstyle{\Em_m u = \sum_k Y_k(\sal) \Em_m x_k |\ol ,\sal \ket.}$$ 
It is clear that we may assume all $x_k$ to be weight vectors of $\UN^-\cdot \,\He^-.$ Then $\wt(x_k) = - (r_{k,0}\alpha_0 + \cdots + r_{k,n-1}\alpha_{n-1})$ for some non-negative integers $r_{k,i}.$ Moreover, $x_k \in \UN\cdot \,\He$ implies  that $x_k |\ol ,\sal \ket$ belongs to $\FN.$  Hence, $x_k |\ol ,\sal \ket$ is a linear combination of vectors $|\lal,\sal\ket$ $(\lal \in \Pi^l),$ and from  formula (\ref{e:carl}) for the weight of $|\lal,\sal\ket$ we see that for all these vectors  $|\lal| = r_k:=r_{k,0}+\cdots + r_{k,n-1}.$ 
Now let $r:=\max_k\{ r_k \},$ and choose the sequence $\sal$ so that the inequalities 
\begin{equation}
s_{b} - s_{b+1} \geq r + nm   \label{e:dom4}
\end{equation}
are satisfied for all $b=1,\dots,l-1.$ Then for each $k$ the vector $x_k |\ol ,\sal \ket$ is a linear combination of $nm$-dominant $|\lal,\sal\ket,$ whereupon Corollary \ref{c:MAIN} (i) gives  
$$ \Em_m u = \sum_k Y_k(\sal) \sum_{m_1+\cdots+m_l=m} \;\prod_{b=1}^l q^{(b-1)m_b}\Em_{m_b}^{(b)}[n,1]\: x_k |\ol,\sal\ket.$$ 
Now we use the assumption in the statement (i) of the lemma, and obtain
$$ (\Em_m u)' = q^{-m(l-1)} (-q)^{m(n-1)} \sum_k Y_k(\sal)' \sum_{m_1+\cdots+m_l=m} \; \prod_{b=1}^l q^{(b-1)m_b}\Hm_{m_b}^{(b)}[n,1]\: (x_k |\ol,\sal\ket)',$$ 
where $Y_k(\sal)'$ is the element of $\ULp$ defined by $Y_k(\sal)'v' = (Y_k(\sal) v)'$ (cf. Proposition \ref{p:ePRIME}). Observe next, that if a vector $|\lal,\sal\ket$ is $nm$-dominant,  then so is $|\lal,\sal\ket'.$ Hence, from (\ref{e:dom4}) it follows that for each $k$  the vector $(x_k |\ol,\sal\ket)'$ is a linear combination of $nm$-dominant $|\lal,\sal'\ket$. Therefore, we may again apply Corollary \ref{c:MAIN} (i), and obtain
$$ (\Em_m u)' = q^{-m(l-1)} (-q)^{m(n-1)} \sum_k Y_k(\sal)' \Hm_m (x_k |\ol ,\sal \ket)' = (-p)^{m(l-1)} (-q)^{m(n-1)} \Hm_m u'.$$ 
Thus (i) is proved.  \finpf

\noindent {\em Proof of Proposition \ref{p:bPRIME}.}
By Lemma \ref{l:bPRIME}, the statement (i) of the proposition will be proved once it is shown to hold for $n=l=1.$ But in this case (i) is just a restatement of the relation $\omega(e_m) = h_m$ for the standard involution of the ring of symmetric functions defined, in terms of the Schur functions, by $\omega(s_{\lambda}) = s_{\lambda'}.$

It remains to observe that, the scalar product being non-degenerate, the relation  (ii) follows from (i), relations (\ref{e:ESCALARP}) and the easily checked formula  $\bra u', v \ket = \ov{ \bra u , v' \ket }.$
\finpf

\begin{proposition} \label{p:SYMMETRY}
For $u,v$ $\in$ $\SIW$ $(s\in\Z)$ one has 
$$ \textstyle{\bra \ov{u}, v \ket  = \bra u' , \ov{ v'}\ket.} $$  
\end{proposition}
\proof Using the decomposition of $\SIW$ described in Theorem \ref{t:TD} we define a gradation of $\SIW$ in the following way. We set the degrees of all the singular vectors $|\ol,\sal\ket $ $(\sal \in A_l^n(s))$ to be zero, and we require that with respect to our gradation the operators $\fn_i,$ $\fl_j,$ and $B_{-m}$ $(m \in \N)$ be homogeneous of degrees $1,$ $1,$ and $m$ respectively. Then, the operators $\en_i,$ $\el_j,$ are homogeneous of degree $-1,$ and the operators $\Em_m,$ $\Hm_m,$ $\Ep_m,$ $\Hp_m$ are homogeneous of respective degrees $m,$ $m,$ $-m,$ $-m.$    

Now we show the claim by induction. In degree zero we have  
\begin{align*}
&\bra \ov{|\ol,\sal\ket }, |\ol,\tal\ket   \ket = \bra {|\ol,\sal\ket }, |\ol,\tal\ket   \ket  = \delta_{ \sal,\tal },\\
&\bra {|\ol,\sal'\ket }, \ov{|\ol,\tal'\ket}   \ket = \bra {|\ol,\sal'\ket }, |\ol,\tal'\ket   \ket  = \delta_{ \sal,\tal }
\end{align*}
for all $\sal,$ $\tal$ $\in$ $A_l^n(s).$ Hence the claim holds in this case.

Assume that the proposition is proved for all $u,v$ with degrees $< k.$ By Theorem \ref{t:TD}, to prove the proposition for all $u,v,$ it is enough to show that  
\begin{alignat}{3}
&\bra \ov{ (\fn_i u)},v \ket &  &=&\; &\bra  (\fn_i u)',\ov{v'} \ket,   \label{e:1} \\
&\bra \ov{ (\fl_j u)},v \ket &  &=&\; &\bra  (\fl_j u)',\ov{v'} \ket,   \label{e:2} \\
&\bra \ov{ (\Em_m w)}, v \ket & &=&\; &\bra  (\Em_m w)', \ov{v'} \ket  \label{e:3}
\end{alignat}
for $u,v,w$ with degrees $k-1,k,k-m$ respectively.

Let us show (\ref{e:1}). We have 
$$\bra \ov{ (\fn_i u)},v \ket  = \bra\fn_i \ov{u},v \ket = \bra\ov{u},q^{-1}\tn_i \en_i v\ket = 
\bra u',\ov{(q^{-1}\tn_i \en_i v)'}\ket. $$ 
Here the first equality follows from Proposition \ref{p:BAR}, the second -- from Proposition \ref{p:eSCALARP} and the third -- from the induction assumption. Further,
$$ \bra u',\ov{(q^{-1}\tn_i \en_i v)'}\ket = \bra u', \ov{ \en_{-i} v'}\ket = \bra u',  \en_{-i} \ov{ v'}\ket  = \bra q^{-1} (\tn_{-i})^{-1} \fn_{-i} u', \ov{v'} \ket = \bra (\fn_i u)' , \ov{ v'}\ket .$$ 
Here we used Proposition \ref{p:ePRIME}, Proposition \ref{p:BAR} and Proposition \ref{p:eSCALARP}. Thus (\ref{e:1}) is established. A proof of (\ref{e:2}) is similar. 

Finally, let us show (\ref{e:3}). We have 
$$ \bra \ov{(\Em_m w)}, v \ket = \bra \Em_m \ov{w} , v \ket = \bra \ov{w} , \Ep_m v \ket = \bra w', \ov{(\Ep_m v)'}\ket. $$ 
Here the first equality comes from Proposition \ref{p:BAR}, the second -- from (\ref{e:ESCALARP}), and the third -- from the induction assumption. Continuing we have 
\begin{gather*} \bra w', \ov{(\Ep_m v)'}\ket = \bra w', \ov{(-q)^{m(n-1)} (-p)^{m(l-1)} \Hp_m v '} \ket = \\ =\bra w', (-q)^{m(n-1)} (-p)^{m(l-1)} \Hp_m \ov{v'}\ket = \bra (-q)^{m(n-1)} (-p)^{m(l-1)} \Hm_m  w' , \ov{v'}\ket = \bra (\Em_m w)', \ov{v'} \ket . 
\end{gather*}
Here we used Proposition \ref{p:bPRIME}, Proposition \ref{p:BAR} and  (\ref{e:ESCALARP}). Thus (\ref{e:3}) is proved. \finpf

For $\sal = (s_1,\dots,s_l) \in \Z^l$ let $\{ \Gc^*(\lal,\sal)\:|\: \lal \in \Pi^l\}$ be the basis of $\FN$ dual to  $\{ \Gc^+(\lal,\sal)\:|\: \lal \in \Pi^l\}$ with respect to the scalar product $\bra u , v \ket.$ Write  
$$ \Gc^*(\lal,\sal) = \sum_{\mul \in \Pi^l} \Delta^*_{\lal,\mul}(\sal | q) \, |\mul,\sal \ket .$$
Since the basis $\{ |\lal,\sal \ket \:|\: \lal \in \Pi^l\}$ of $\FN$ is orthonormal relative to the scalar product, the matrix $\|\Delta^*_{\lal,\mul}(\sal | q)\|$ is the transposed inverse of the matrix $\|\Delta^+_{\lal,\mul}(\sal | q)\|.$

\begin{proposition}
For $\sal \in \Z^l,$ $\lal \in \Pi^l$ one has 
$$ \Gc^*(\lal,\sal)'  = \Gc^-(\lal',\sal'). $$
\end{proposition}
\proof 
Since the matrix $\|\Delta^+_{\lal,\mul}(\sal | q)\|$ is unitriangular with off-diagonal entries in $q\Z[q],$ the same is true for its transposed inverse  $\|\Delta^*_{\lal,\mul}(\sal | q)\|.$ It follows that 
$$ \Gc^*(\lal,\sal)' \equiv |\lal',\sal'\ket \bmod q^{-1} \Lc^-(s) $$
where $s = s_1+ \cdots + s_l.$ Hence, $ \Gc^*(\lal,\sal)' $ has the required congruence property relative to the basis  $\{|\lal,\sal'\ket \}.$ 

It remains to show that  $ \Gc^*(\lal,\sal)' $ is invariant with respect to the bar-involution. Since $ \Gc^*(\lal,\sal) $ is dual to $ \Gc^+(\lal,\sal) $ this is equivalent to     
$$ \bra \ov{ \Gc^*(\lal,\sal)'}, \Gc^+(\mul,\sal)' \ket = \delta_{\lal,\mul}.$$ 
Using Proposition \ref{p:SYMMETRY} we obtain 
$$ \bra \ov{ \Gc^*(\lal,\sal)'}, \Gc^+(\mul,\sal)' \ket = \bra { \Gc^*(\lal,\sal)}, \ov{\Gc^+(\mul,\sal)} \ket  = \bra { \Gc^*(\lal,\sal)}, {\Gc^+(\mul,\sal)} \ket = \delta_{\lal,\mul}.$$ \finpf 

\noindent This proposition immediately implies 
\begin{theorem}[Inversion formula] \label{t:MAIN}
For $\sal \in \Z^l,$ $\lal,\mul \in \Pi^l$ one has 
$$ \sum_{\nul \in \Pi^l} \textstyle{\Delta^-_{\lal',\nul'}(\sal'|q^{-1})\: \Delta^+_{\mul,\nul}(\sal|q) = \delta_{\lal,\mul}.}$$ 
\end{theorem}

\begin{example}{\rm 
To illustrate the theorem, let $n=3,l=2$ and $\sal = (1,-1).$ Then the matrix $\|\Delta^+_{\lal,\mul}(\sal|q)\|_3$  consists of the two following blocks:
\[
\begin{array}{l c c c c c c c c}
 ((3),\vr)   & 1 & 0 & 0 & 0 & 0 & 0 & 0 & 0 \cr 
 (\vr,(3))   & 0 & 1 & 0 & 0 & 0 & 0 & 0 & 0 \cr 
 ((1),(2))   & q & q & 1 & 0 & 0 & 0 & 0 & 0 \cr 
 (\vr,(2,1)) & 0& {q^2} & q & 1 & 0 & 0 & 0 & 0 \cr 
 ((2,1),\vr) & q & 0 & 0 & 0 & 1 & 0 & 0 & 0 \cr 
((1^2),(1))  & {q^2} & 0 & q & 0 & q & 1 & 0 & 0 \cr 
((1^3),\vr)  &  0 & 0 & 0 & 0 & {q^2} & q & 1 & 0 \cr 
(\vr, (1^3)) &  0 & 0 & {q^2} & q & 0 & q & 0 & 1  
\end{array}
\qquad
\begin{array}{l c c}
((2),(1))   & 1 & 0 \cr 
((1),(1^2)) & 0 & 1 
\end{array}
\]
Here coefficients in each column are those of the  expansion of an element of the canonical basis on the basis $\{ |\lal\,\sal \ket \}.$ For example, $\Gc^+((\vr,(3)),\sal) = |(\vr,(3)),\sal\ket + q |((1),(2)),\sal\ket + q^2 |(\vr,(2,1)),\sal\ket.$   

On the other hand, the matrix $\|\Delta^-_{\lal,\mul}(\sal|q)\|_3$  consists of the following two blocks (we put $p=-q^{-1}$): 
\[
\begin{array}{l c c c c c c c c}
((3),\vr)    &1 & 0 & 0 & 0 & 0 & 0 & 0 & 0 \cr 
(\vr,(3))    &0 & 1 & 0 & 0 & 0 & 0 & 0 & 0 \cr 
((1),(2))    &p & p & 1 & 0 & 0 & 0 & 0 & 0 \cr 
(\vr,(2,1))  &{p^2} & 0 & p & 1 & 0 & 0 & 0 & 0 \cr 
((2,1),\vr)  &p & 0 & 0 & 0 & 1 & 0 & 0 & 0 \cr 
((1^2),(1))  &{p^2} & {p^2} & p & 0 & p & 1 & 0 & 0 \cr 
((1^3),\vr)  &0 & {p^3} & {p^2} & 0 & 0 & p & 1 & 0 \cr 
(\vr, (1^3)) &{p^3} & 0 & {p^2} & p & {p^2} & p & 0 & 1    
\end{array}
\qquad 
\begin{array}{l c c}
((2),(1))   & 1 & 0 \cr 
((1),(1^2)) & 0 & 1 
\end{array}
\]
It is easy to verify that $\|\Delta^-_{\lal',\mul'}(\sal|q^{-1})\|_3$ is the transposed inverse of $\|\Delta^+_{\lal,\mul}(\sal|q)\|_3$ (Note that in the present case $\sal' = \sal$).
}
\finex \end{example}

\subsection{Proof of Proposition \ref{p:MAIN}.} \label{ss:PROOF}

Let us prove part (i) of Proposition \ref{p:MAIN}. First, we  introduce some notations. 
Let $s$ be an integer. For any pair $(\lal,\sal)$ $\in$ $\Pi^l\times\Z^l(s)$ where 
$$\lal  =(\lambda^{(1)},\dots,\lambda^{(l)}),\qquad \sal=(s_1,\dots,s_l),$$
 let $\kb = (k_1,k_2,\dots )$ be the unique sequence from $\PP^{++}(s)$ (cf. Section \ref{ss:SIW})  such that 
$$ |\lal,\sal\ket = u_{\kb}.$$ 
As in Section \ref{ss:SIW} we write $k_i = a_i + n(b_i - 1) - nlm_i$ where $a_i \in \{1,\dots,n\},$ $b_i \in \{1,\dots,l\}$ and $m_i \in \Z.$  

For any natural number  $r$ put
$(\kb)_r = (k_1,k_2,\dots,k_r).$ Then $(\kb)_r \in \P^{++},$ and $u_{\kb} = u_{(\kb)_r}\w u_{k_{r+1}}\w u_{k_{r+2}}\w \cdots,$ where $u_{(\kb)_r} = u_{k_1}\wedge \cdots \wedge u_{k_r}$ is an element of $\Lambda^r.$   

We define  $(\kb)_r^+ = (k_1^+,\dots,k_r^+)$ to be the  unique  permutation of the sequence $(\kb)_r$ characterized by the following two conditions 
\begin{alignat*}{2}
&b^+_i \leq b_j^+  &&\text{ for all $i<j,$} \\ 
&k_i^+ > k_j^+  &&\text{ for all $i<j$ such that $b_i^+ = b_j^+.$}   
\end{alignat*}
Here  we  put $k_i^+ = a^+_i + n(b^+_i-1) - nl m_i^+$ where $a_i^+ \in \{1,\dots,n\},$ $b_i^+ \in \{1,\dots,l\}$ and $m_i^+ \in \Z.$

\begin{example} \label{ex:kr}{\rm 
Let $n=2, l=3$ and $s=-2.$ Take  $\lal = ((1^2),(1),(2))$ and $\sal = (5,0,-7).$ Let $r=25.$ In this case 
$$ (\kb)_r = \scriptstyle{(14,13,7,3,2,1,-3,-4,-5,
   -8,-9,-10,-11,-13,-14,-15,-16,-17,
   -20,-21,-22,-23,-24,-25,-26),} 
$$
and 
$$
(\kb)_r^+= \scriptstyle{ (14, 13, 7, 2, 1, -4, -5, -10, -11, -16, -17, -22, -23, 3, -3, -8, 
     -9, -14,  -15, -20, -21, -26, -13, -24, -25).} 
$$
}
\finex  \end{example} 
\noindent Recall that in Section \ref{ss:SIW} we associated with $\kb$ the  semi-infinite sequences 
$$\kb^{(b)} = (k_1^{(b)},k_2^{(b)},\dots)\qquad (b=1,2,\dots,l),$$
  such that $ k_i^{(b)} = s_b + i - 1 + \lambda_i^{(b)}$ for all $i \in \N.$  The wedge $u_{(\kb)^+_r} = u_{k_1^+}\wedge \cdots \wedge u_{k_r^+}$ may be expressed in terms of these sequences in the following way. For $a \in \{1,\dots,n\},$ $b \in \{1,\dots,l\}$ and $m\in \Z,$ put $u^{(b)}_{a-nm}:= u_{a+n(b-1)-nlm}.$ Then 
\begin{align} 
u_{(\kb)^+_r} = &\textstyle{u^{(1)}_{k_1^{(1)}}\wedge u^{(1)}_{k_2^{(1)}}\wedge\cdots \wedge u^{(1)}_{k_{r_1}^{(1)}}\wedge 
    u^{(2)}_{k_1^{(2)}}\wedge u^{(2)}_{k_2^{(2)}}\wedge\cdots \wedge u^{(2)}_{k_{r_2}^{(2)}}\wedge} \cdots  \label{e:eq3}\\
    & \qquad \qquad \qquad \cdots   \nonumber\\
 &\cdots \wedge u^{(l)}_{k_1^{(l)}}\wedge u^{(l)}_{k_2^{(l)}}\wedge\cdots \wedge u^{(l)}_{k_{r_l}^{(l)}} \nonumber
\end{align}
where for each $b \in \{1,\dots,l\}$ we put $r_b:= \sharp\{ 1\leq i \leq r\:|\: b_i = b\}.$ 

Note that in general the wedge $u_{(\kb)^+_r}$ is not ordered, and using the ordering rules of Proposition \ref{p:RULES} to straighten  $u_{(\kb)^+_r}$ as a linear combination of ordered wedges one typically obtains a linear combination with many terms. The first step towards the proof of the proposition is to show that if the pair $(\lal,\sal)$ we started with is $0$-dominant, then, for any $r \in \N,$ the straightening of $u_{(\kb)^+_r}$ produces only one term, which, up to a power of $q,$ coincides with $u_{(\kb)_r}.$

\begin{lemma} \label{l:STEP}
Let $b_1,b_2 \in \{1,\dots,l\},$ and $a_1,a_2 \in \{1,\dots,n\}$ satisfy the inequalities $b_1 < b_2,$ $a_1 \geq a_2.$ Let $m\in \Z.$ Then, for any $t\in \Z_{\geq 0},$ one has the following relation
\begin{multline*}
u^{(b_1)}_{a_1-nm}\w u^{(b_1)}_{a_1-nm-1}\w u^{(b_1)}_{a_1-nm-2}\w\cdots \w u^{(b_1)}_{a_1-nm-t}\w u^{(b_2)}_{a_2-nm} = \\ 
= q^{\sum_{k=0}^t \delta_{ (a_1 - k) \equiv a_2 \bmod n}} u^{(b_2)}_{a_2-nm}\w u^{(b_1)}_{a_1-nm}\w u^{(b_1)}_{a_1-nm-1}\w u^{(b_1)}_{a_1-nm-2}\w\cdots \w u^{(b_1)}_{a_1-nm-t}.
\end{multline*}
\end{lemma}
\proof This is shown by induction on $t$ using relations (R3) and (R4) of Proposition \ref{p:RULES}, and Lemma \ref{l:ZERO}.

\finpf

\noindent Keeping $(\lal,\sal),$ $\kb,$ $r,$ $(\kb)_r$ and $(\kb)_r^+$ as above, let us define
$$
c_r(\kb):= \sharp\{ 1\leq i < j\leq r \: |\: {b_i > b_j},\; a_i=a_j \},
$$ and $$ 
|\lal,\sal \ket_r^+ : =  u_{(\kb)_r^+}\wedge u_{k_{r+1}}\w u_{k_{r+2}} \w \cdots .
$$
\begin{lemma} \label{l:STEP2}
Suppose $(\lal,\sal)$ is $0$-dominant. Then 
$$ |\lal,\sal\ket_r^+ = q^{c_r(\kb)} |\lal,\sal \ket.$$ 
\end{lemma}
\proof Since $|\lal,\sal \ket = u_{\kb} = u_{(\kb)_r}\wedge u_{k_{r+1}}\w u_{k_{r+2}} \w \cdots ,$ we must prove that 
 $$ u_{(\kb)_r^+} = q^{c_r(\kb)} u_{(\kb)_r}.$$ 
First of all, let us examine what the $0$-dominance of $(\lal,\sal)$ implies for the semi-infinite sequence $\kb.$ For each $b \in \{1,\dots,l\},$ let $p_b$ be the minimal number such that $k_{i}^{(b)} = s_b - i + 1$ for all $i\geq p_b.$ Then $p_b = l(\lambda^{(b)})+1,$ where $l(\lambda)$ denotes the length of a partition $\lambda,$ and we have $k_{p_b}^{(b)} = s_b - l(\lambda^{(b)}).$ On the other hand, $k_1^{(b)} = s_b + \lambda^{(b)}_1.$ Hence, using the assumption that $(\lal,\sal)$ is $0$-dominant (cf. Definition \ref{d:dom}) we find that   for all $b=1,2,\dots,l-1,$ 
$$ \textstyle{k_{p_b}^{(b)} - k_1^{(b+1)} = s_b - s_{b+1} - l(\lambda^{(b)}) - \lambda_1^{(b+1)} \geq s_b-s_{b+1} - |\lal| \geq 0.}$$

\mbox{} The  fact that we have the inequalities $k_{p_b}^{(b)} \geq  k_1^{(b+1)}$ for all $b=1,2,\dots,l-1,$ implies that to straighten $ u_{(\kb)_r^+} $ on the basis of ordered wedges we only need to repeatedly apply Lemma \ref{l:STEP}. The result follows. \finpf    

\begin{example} {\rm 
Let us illustrate the proof of Lemma \ref{l:STEP2} for $|\lal,\sal \ket^+_r$ where $(\lal,\sal)$ and $r$ are the same as in Example \ref{ex:kr}. Note that the pair $(\lal,\sal)$ is $0$-dominant. 

In this case $ u_{(\kb)_r}$ is given by the following expression:
\begin{scriptsize}
\begin{align*}
&u^{(1)}_{2+2\cdot 2}\w u^{(1)}_{1+2\cdot 2}\w u^{(1)}_{1+2}\w u^{(2)}_{1}\w u^{(1)}_{2}\w u^{(1)}_{1}\w u^{(2)}_{1-2\cdot 1}\w u^{(1)}_{2-2\cdot 1}\w u^{(1)}_{1-2\cdot 1}\w u^{(2)}_{2-2\cdot 2}\w u^{(2)}_{1-2\cdot 2}\w u^{(1)}_{2-2\cdot 2}\w u^{(1)}_{1-2\cdot 2}\w  \\ 
&\w u^{(3)}_{1-2\cdot 3}\w u^{(2)}_{2-2\cdot 3}\w u^{(2)}_{1-2\cdot 3}\w u^{(1)}_{2-2\cdot 3}\w u^{(1)}_{1-2\cdot 3}\w u^{(2)}_{2-2\cdot 4}\w u^{(2)}_{1-2\cdot 4}\w u^{(1)}_{2-2\cdot 4}\w u^{(1)}_{1-2\cdot 4}
\w u^{(3)}_{2-2\cdot 5}\w u^{(3)}_{1-2\cdot 5} \w u^{(2)}_{2-2\cdot 5}. 
\end{align*}
\end{scriptsize}
\vspace{-2mm}

\noindent Now let us rearrange (taking care of powers of $q$) the factors in this wedge by repeatedly applying Lemma \ref{l:STEP}. The rearrangement involves the following $7$ steps: 
$$u_{(\kb)_r} =$$ 
\begin{scriptsize}
\begin{align*}
%&u_{(\kb)_r} = \\
%&\mbox{} \\ 
&u^{(1)}_{2+2\cdot 2}\w u^{(1)}_{1+2\cdot 2}\w u^{(1)}_{1+2}\w u^{(2)}_{1}\w u^{(1)}_{2}\w u^{(1)}_{1}\w u^{(2)}_{1-2\cdot 1}\w u^{(1)}_{2-2\cdot 1}\w u^{(1)}_{1-2\cdot 1}\w u^{(2)}_{2-2\cdot 2}\w u^{(2)}_{1-2\cdot 2}\w u^{(1)}_{2-2\cdot 2}\w u^{(1)}_{1-2\cdot 2}\w  \\ 
&\w u^{(3)}_{1-2\cdot 3}\w u^{(2)}_{2-2\cdot 3}\w u^{(2)}_{1-2\cdot 3}\w u^{(1)}_{2-2\cdot 3}\w u^{(1)}_{1-2\cdot 3}\w u^{(2)}_{2-2\cdot 4}\w u^{(2)}_{1-2\cdot 4}\w u^{(1)}_{2-2\cdot 4}\w u^{(1)}_{1-2\cdot 4}
\w \un{u^{(3)}_{2-2\cdot 5}\w u^{(3)}_{1-2\cdot 5} \w u^{(2)}_{2-2\cdot 5}} = \\ 
&\mbox{}  \\ 
&q^{-1} u^{(1)}_{2+2\cdot 2}\w u^{(1)}_{1+2\cdot 2}\w u^{(1)}_{1+2}\w u^{(2)}_{1}\w u^{(1)}_{2}\w u^{(1)}_{1}\w u^{(2)}_{1-2\cdot 1}\w u^{(1)}_{2-2\cdot 1}\w u^{(1)}_{1-2\cdot 1}\w u^{(2)}_{2-2\cdot 2}\w u^{(2)}_{1-2\cdot 2}\w u^{(1)}_{2-2\cdot 2}\w u^{(1)}_{1-2\cdot 2}\w  \\ 
&\w u^{(3)}_{1-2\cdot 3}\w u^{(2)}_{2-2\cdot 3}\w u^{(2)}_{1-2\cdot 3}\w u^{(1)}_{2-2\cdot 3}\w u^{(1)}_{1-2\cdot 3}\w \un{u^{(2)}_{2-2\cdot 4}\w u^{(2)}_{1-2\cdot 4}\w u^{(1)}_{2-2\cdot 4}\w u^{(1)}_{1-2\cdot 4}}
\w u^{(2)}_{2-2\cdot 5} \w u^{(3)}_{2-2\cdot 5}\w u^{(3)}_{1-2\cdot 5}= \\ 
&\mbox{} \\   
&q^{-3} u^{(1)}_{2+2\cdot 2}\w u^{(1)}_{1+2\cdot 2}\w u^{(1)}_{1+2}\w u^{(2)}_{1}\w u^{(1)}_{2}\w u^{(1)}_{1}\w u^{(2)}_{1-2\cdot 1}\w u^{(1)}_{2-2\cdot 1}\w u^{(1)}_{1-2\cdot 1}\w u^{(2)}_{2-2\cdot 2}\w u^{(2)}_{1-2\cdot 2}\w u^{(1)}_{2-2\cdot 2}\w u^{(1)}_{1-2\cdot 2}\w  \\ 
&\w u^{(3)}_{1-2\cdot 3}\w \un{u^{(2)}_{2-2\cdot 3}\w u^{(2)}_{1-2\cdot 3}\w u^{(1)}_{2-2\cdot 3}\w u^{(1)}_{1-2\cdot 3}\w u^{(1)}_{2-2\cdot 4}\w u^{(1)}_{1-2\cdot 4}}\w u^{(2)}_{2-2\cdot 4}\w u^{(2)}_{1-2\cdot 4}
\w u^{(2)}_{2-2\cdot 5} \w u^{(3)}_{2-2\cdot 5}\w u^{(3)}_{1-2\cdot 5} = \\ 
&\mbox{} \\ 
&q^{-7}u^{(1)}_{2+2\cdot 2}\w u^{(1)}_{1+2\cdot 2}\w u^{(1)}_{1+2}\w u^{(2)}_{1}\w u^{(1)}_{2}\w u^{(1)}_{1}\w u^{(2)}_{1-2\cdot 1}\w u^{(1)}_{2-2\cdot 1}\w u^{(1)}_{1-2\cdot 1}\w u^{(2)}_{2-2\cdot 2}\w u^{(2)}_{1-2\cdot 2}\w u^{(1)}_{2-2\cdot 2}\w u^{(1)}_{1-2\cdot 2}\w  \\ 
&\w \un{u^{(3)}_{1-2\cdot 3}\w u^{(1)}_{2-2\cdot 3}\w u^{(1)}_{1-2\cdot 3}\w u^{(1)}_{2-2\cdot 4}\w u^{(1)}_{1-2\cdot 4}\w u^{(2)}_{2-2\cdot 3}\w u^{(2)}_{1-2\cdot 3}\w u^{(2)}_{2-2\cdot 4}\w u^{(2)}_{1-2\cdot 4}
\w u^{(2)}_{2-2\cdot 5}}\w u^{(3)}_{2-2\cdot 5}\w u^{(3)}_{1-2\cdot 5} = \\ 
&\mbox{} \\ 
&q^{-11}u^{(1)}_{2+2\cdot 2}\w u^{(1)}_{1+2\cdot 2}\w u^{(1)}_{1+2}\w u^{(2)}_{1}\w u^{(1)}_{2}\w u^{(1)}_{1}\w u^{(2)}_{1-2\cdot 1}\w u^{(1)}_{2-2\cdot 1}\w u^{(1)}_{1-2\cdot 1}\w \un{u^{(2)}_{2-2\cdot 2}\w u^{(2)}_{1-2\cdot 2}\w u^{(1)}_{2-2\cdot 2}\w u^{(1)}_{1-2\cdot 2}\w } \\ 
&\un{\w u^{(1)}_{2-2\cdot 3}\w u^{(1)}_{1-2\cdot 3} \w u^{(1)}_{2-2\cdot 4}\w u^{(1)}_{1-2\cdot 4}}\w u^{(2)}_{2-2\cdot 3}\w u^{(2)}_{1-2\cdot 3}\w u^{(2)}_{2-2\cdot 4}\w u^{(2)}_{1-2\cdot 4}
\w u^{(2)}_{2-2\cdot 5}\w u^{(3)}_{1-2\cdot 3}\w u^{(3)}_{2-2\cdot 5}\w u^{(3)}_{1-2\cdot 5} =\\ 
&\mbox{} \\  
&q^{-17}u^{(1)}_{2+2\cdot 2}\w u^{(1)}_{1+2\cdot 2}\w u^{(1)}_{1+2}\w u^{(2)}_{1}\w u^{(1)}_{2}\w u^{(1)}_{1}\w \un{u^{(2)}_{1-2\cdot 1}\w u^{(1)}_{2-2\cdot 1}\w u^{(1)}_{1-2\cdot 1}\w u^{(1)}_{2-2\cdot 2}\w u^{(1)}_{1-2\cdot 2}\w u^{(1)}_{2-2\cdot 3}\w u^{(1)}_{1-2\cdot 3}\w } \\ 
&\un{\w u^{(1)}_{2-2\cdot 4}\w u^{(1)}_{1-2\cdot 4}}\w u^{(2)}_{2-2\cdot 2}\w u^{(2)}_{1-2\cdot 2}\w u^{(2)}_{2-2\cdot 3}\w u^{(2)}_{1-2\cdot 3}\w u^{(2)}_{2-2\cdot 4}\w u^{(2)}_{1-2\cdot 4}
\w u^{(2)}_{2-2\cdot 5}\w u^{(3)}_{1-2\cdot 3}\w u^{(3)}_{2-2\cdot 5}\w u^{(3)}_{1-2\cdot 5} =\\ 
&\mbox{} \\ 
&q^{-21}u^{(1)}_{2+2\cdot 2}\w u^{(1)}_{1+2\cdot 2}\w u^{(1)}_{1+2}\w \un{u^{(2)}_{1}\w u^{(1)}_{2}\w u^{(1)}_{1}\w u^{(1)}_{2-2\cdot 1}\w u^{(1)}_{1-2\cdot 1}\w u^{(1)}_{2-2\cdot 2}\w u^{(1)}_{1-2\cdot 2}\w u^{(1)}_{2-2\cdot 3}\w u^{(1)}_{1-2\cdot 3}\w u^{(1)}_{2-2\cdot 4}\w } \\ 
&\un{\w u^{(1)}_{1-2\cdot 4}}\w u^{(2)}_{1-2\cdot 1}\w u^{(2)}_{2-2\cdot 2}\w u^{(2)}_{1-2\cdot 2}\w u^{(2)}_{2-2\cdot 3}\w u^{(2)}_{1-2\cdot 3}\w u^{(2)}_{2-2\cdot 4}\w u^{(2)}_{1-2\cdot 4}
\w u^{(2)}_{2-2\cdot 5}\w u^{(3)}_{1-2\cdot 3}\w u^{(3)}_{2-2\cdot 5}\w u^{(3)}_{1-2\cdot 5} = \\ 
&\mbox{} \\
&q^{-26}u^{(1)}_{2+2\cdot 2}\w u^{(1)}_{1+2\cdot 2}\w u^{(1)}_{1+2}\w { u^{(1)}_{2}\w u^{(1)}_{1}\w u^{(1)}_{2-2\cdot 1}\w u^{(1)}_{1-2\cdot 1}\w u^{(1)}_{2-2\cdot 2}\w u^{(1)}_{1-2\cdot 2}\w u^{(1)}_{2-2\cdot 3}\w u^{(1)}_{1-2\cdot 3}\w u^{(1)}_{2-2\cdot 4}\w u^{(1)}_{1-2\cdot 4}\w } \!\!\!\!\!\!\!\\ 
&\w u^{(2)}_{1}\w u^{(2)}_{1-2\cdot 1}\w u^{(2)}_{2-2\cdot 2}\w u^{(2)}_{1-2\cdot 2}\w u^{(2)}_{2-2\cdot 3}\w u^{(2)}_{1-2\cdot 3}\w u^{(2)}_{2-2\cdot 4}\w u^{(2)}_{1-2\cdot 4}
\w u^{(2)}_{2-2\cdot 5}\w u^{(3)}_{1-2\cdot 3}\w u^{(3)}_{2-2\cdot 5}\w u^{(3)}_{1-2\cdot 5} = %\\
%&\mbox{} \\
%& = q^{-26} u_{(\kb)^+_r}. 
\end{align*}
\end{scriptsize}
$$  = q^{-26} u_{(\kb)^+_r}.$$
\noindent Here at each step we underline the part to which we apply Lemma \ref{l:STEP}. Note that we use this lemma twice at steps $1,2,3,4,5$ and we use it once at steps $6$ and $7.$ \finex 
}
\end{example}

\noindent Let temporarily  $l=1,$ and for partitions $\lambda$ and $\mu$ define the matrix elements $(B_{-m})_{\lambda,\mu}(s)$ by   
$$ B_{-m}|\lambda,s\ket = \sum_{\mu \in \Pi} (B_{-m})_{\lambda,\mu}(s) \, |\mu,s\ket.
$$ 

Now we proceed with the proof of the proposition. Assume that pair $(\lal,\sal)$ is $nm$-dominant for some $m \in \N.$ Then 
$$ B_{-m}|\lal,\sal\ket = B_{-m} u_{\kb} = ( B_{-m} u_{(\kb)_r})\wedge |s-r\ket = q^{-c_r(\kb)} ( B_{-m} u_{(\kb)^+_r})\wedge |s-r\ket ,$$ 
where the second equality is obtained by taking $r$ sufficiently large, and the third equality follows from Lemma \ref{l:STEP2}. Using (\ref{e:eq3}) we have
$$  B_{-m} u_{(\kb)^+_r} = \sum_{b=1}^l \sum_{i=1}^{r_b} u^{(1)}_{(\kb^{(1)})_{r_1}} \w u^{(2)}_{(\kb^{(2)})_{r_2}} \w \cdots \w u^{(b)}_{(\kb^{(b)})_{r_b} - \epsilon_i nm} \w \cdots \w u^{(l)}_{(\kb^{(l)})_{r_l}},
$$ 
where we put $(\kb^{(b)})_{r_b} = (k^{(b)}_1,\dots,k^{(b)}_{r_b}),$ and $\epsilon_i = (0,\dots,0,1,0\dots,0)$ with $1$ on the $i$th position. Let us now straighten on the basis of ordered wedges the expression   
$$ 
{\sum_{i=1}^{r_b} \; u^{(b)}_{(\kb^{(b)})_{r_b} - \epsilon_i nm}.} 
$$ 
It is clear that to do so we only need to use relations (R1) and (R2) of Proposition \ref{p:RULES}. But these two relations are the same as in the case $l=1.$ Hence, assuming (as we may by choosing large enough $r$) that $r_1,r_2,\dots,r_l$ are sufficiently large, we get  
$$ B_{-m} |\lal,\sal\ket = q^{-c_r(\kb)} \sum_{b=1}^l \sum_{\mu \in \Pi} \, (B_{-m})_{\lambda^{(b)},\mu}(s_b) \, | (\lambda^{(1)},\dots, \lambda^{(b-1)}, \mu , \lambda^{(b+1)}, \dots , \lambda^{(l)}), \sal \ket_r^+. 
$$ 
Observe now, that in the above sum we have, for all $b$ and $\mu,$ 
$$ |\mu| = |\lambda^{(b)}| + nm.$$ 
Hence, $nm$-dominance of $(\lal,\sal)$ implies $0$-dominance of each pair 
\begin{equation}
((\lambda^{(1)},\dots, \lambda^{(b-1)}, \mu , \lambda^{(b+1)}, \dots , \lambda^{(l)}), \sal ). \label{e:eq4}
\end{equation} 
It follows that we may apply Lemma \ref{l:STEP2} in order to straighten each wedge under the sum. Doing so, we get 
$$ 
B_{-m} |\lal,\sal\ket =  \sum_{b=1}^l \sum_{\mu \in \Pi} \,q^{c_r(\lb)-c_r(\kb)}\, (B_{-m})_{\lambda^{(b)},\mu}(s_b) \, | (\lambda^{(1)},\dots, \lambda^{(b-1)}, \mu , \lambda^{(b+1)}, \dots , \lambda^{(l)}), \sal \ket, 
$$ 
where $\lb$ is the unique  element of $\PP^{++}(s)$ such that  
$$ 
| (\lambda^{(1)},\dots, \lambda^{(b-1)}, \mu , \lambda^{(b+1)}, \dots , \lambda^{(l)}), \sal \ket = u_{\lb}.
$$ 
Finally, using $0$-dominance of (\ref{e:eq4}) and $0$-dominance of $(\lal,\sal),$ it is not difficult to see that $c_r(\lb) - c_r(\kb) = (b-1)m $ for all large enough   $r.$   Proposition \ref{p:MAIN} (i) follows.

A proof of relation (i') of   Proposition \ref{p:MAIN} is obtained by following the same steps as above, but interchanging everywhere the roles of $n$ and $l,$ and the roles of $p$ and $q.$ Proofs of (ii) and (ii') are similar to proofs of (i) and (i'), and will be omitted. \finpf  

\section{Examples of the matrices $\|\Delta^+_{\lal,\mul}(\sal|q)\|_k.$ }
In this section we display the transition matrices  $\|\Delta^+_{\lal,\mul}(\sal|q)\|_k$ for $n=l=2,$ and $\sal = (0,0),$ up to $k=6.$ These matrices should be read by columns, for example, at $k=3$ we have    
$$ \Gc^+(((3),\vr),\sal) = |((3),\vr),\sal\ket + q |((2),(1)),\sal\ket + q |((1^2),(1)),\sal\ket + q^2 |((1^3),\vr),\sal\ket .$$
Each square  matrix corresponds to a weight subspace of $\FN.$ The large matrices for weights $2\Lambda_0 - 3\alpha_0 - 2\alpha_1$ and $2\Lambda_0 - 3\alpha_0 - 3\alpha_1$ are split into several parts. The rows are labelled by  multipartitions  $\mul$ together with  partitions $\mu$ such that $|\mul,\sal\ket = |\mu,s\ket,$ and are arranged in the reverse lexicographic order with respect to $\mu.$ Finally, we mark by $*$ the rows labelled by those  $\mul$ which belong to $\Pi^l(\sal).$ The corresponding columns, therefore, give expansions of the elements of the lower global crystal basis of the irreducible submodule $\MN.$   

\medskip

%\documentclass[11pt,a4paper]{article}

%\usepackage{amstext,amssymb,amsfonts,amsmath}

%\setcounter{MaxMatrixCols}{20}

%\def\lal{{\boldsymbol {\lambda}}_l}
%\def\vr{\varnothing}
%\def\wt{{\mathrm {wt}}}
%\begin{document}

\noindent $k = 1$
\begin{align*}
&\text{ $\wt = 2\Lambda_0 - \alpha_0$} \\
&\begin{array}{l|l| c c}
\!\!\!\! *\:(3) &  (\vr, (1)) &    1 & \cdot \cr 
(1^3) &  ((1), \vr) &  q & 1    
\end{array}
\end{align*}
\noindent $k=2$
\begin{align*}
&\text{ $\wt = 2\Lambda_0 - \alpha_0-\alpha_1$} \qquad &  &\text{ $\wt = 2\Lambda_0 - 2\alpha_0$} \\
&\begin{array}{l|l| c c c c }
\!\!\!\! *\:(4) &  (\vr, (2)) &  1 & \cdot & \cdot & \cdot \cr 
(3, 1) &  (\vr, (1^2)) & q & 1 & \cdot & \cdot \cr 
(2, 1^2) &  ((2), \vr) & q & \cdot & 1 & \cdot \cr 
(1^4) &  ((1^2), \vr) & {q^2} & q & q & 1 
\end{array}\qquad & & 
\begin{array}{l|l| c } 
\!\!\!\! *\: (3, 2, 1) &  ((1), (1)) & 1  
\end{array}
\end{align*}

\noindent $k = 3$
\begin{align*}
&\text{ $\wt = 2\Lambda_0 - \alpha_0-2\alpha_1$}\qquad & & \text{ $\wt = 2\Lambda_0 - 2\alpha_0-\alpha_1$} \\
&\begin{array}{l|l| c c}
\!\!\!\! *\:(4, 1) &  (\vr, (2, 1)) & 1 & \cdot \cr 
(2, 1^3) &  ((2, 1), \vr) & q & 1    
\end{array}\qquad & & 
\begin{array}{l|l| c c c c c c c c } 
\!\!\!\! *\:(7) &  (\vr, (3)) &   1 & \cdot & \cdot & \cdot & \cdot & \cdot & \cdot & \cdot \cr 
(5, 1^2) &  ((3), \vr) & q & 1 & \cdot & \cdot & \cdot & \cdot & \cdot & \cdot \cr 
\!\!\!\! *\:(4, 2, 1) &  ((1), (2)) & q & \cdot & 1 & \cdot & \cdot & \cdot & \cdot & \cdot \cr 
(3^2 , 1) &  ((2), (1)) & {q^2} & q & q & 1 & \cdot & \cdot & \cdot & \cdot \cr 
(3, 2^2) &  ((1), (1^2)) & q & \cdot & {q^2} & q & 1 & \cdot & \cdot & \cdot \cr 
(3, 2, 1^2) &  ((1^2), (1)) & {q^2} & q & {q^3} & {q^2} & q & 1 & \cdot & \cdot \cr 
(3, 1^4) &  (\vr, (1^3)) & {q^2} & \cdot & \cdot & \cdot & q & \cdot & 1 & \cdot \cr 
(1^7) &  ((1^3), \vr) & {q^3} & {q^2} & \cdot & q & {q^2} & q & q & 1  
\end{array}
\end{align*}
$k = 4$
\begin{align*}
&\text{ $\wt = 2\Lambda_0 - 2\alpha_0 - 2\alpha_1$} \\
&\begin{array}{l | l |   c c c c c c c c c c c c c c c c}
\!\!\!\! *\:(8) &  (\vr, (4))  &    
1 & \cdot  & \cdot  & \cdot  & \cdot  & \cdot  & \cdot  & \cdot  & \cdot  & \cdot  & \cdot  & \cdot  & \cdot  & \cdot  & \cdot  & \cdot  \cr 
\!\!\!\! *\:(7, 1)& (\vr, (3, 1))  &  
q & 1 & \cdot  & \cdot  & \cdot  & \cdot  & \cdot  & \cdot  & \cdot  & \cdot  & \cdot  & \cdot  & \cdot  & \cdot  & \cdot  & \cdot  \cr 
(6, 1^2) &  ((4), \vr)  &  
q & \cdot  & 1 & \cdot  & \cdot  & \cdot  & \cdot  & \cdot  & \cdot  & \cdot  & \cdot  & \cdot  & \cdot  & \cdot  & \cdot  & \cdot  \cr 
(5, 1^3) &  ((3, 1), \vr)  &  
{q^ 2} & q & q & 1 & \cdot  & \cdot  & \cdot  & \cdot  & \cdot  & \cdot  & \cdot  & \cdot  & \cdot  & \cdot  & \cdot  & \cdot  \cr 
(4^2) &  (\vr, (2^2))  &   
\cdot  & q & \cdot  & \cdot  & 1 & \cdot  & \cdot  & \cdot  & \cdot  & \cdot  & \cdot  & \cdot  & \cdot  & \cdot  & \cdot  & \cdot  \cr 
\!\!\!\! *\:(4, 3, 1) &  ((2), (2))  &      
{q^2} & \cdot  & q & \cdot  & \cdot  & 1 & \cdot  & \cdot  & \cdot  & \cdot  & \cdot  & \cdot  & \cdot  & \cdot  & \cdot  & \cdot  \cr 
(4, 2^2) &  ((1), (2, 1))  &    
q & {q^2} & \cdot  & \cdot  & q & q & 1 & \cdot  & \cdot  & \cdot  & \cdot  & \cdot  & \cdot  & \cdot  & \cdot  & \cdot  \cr 
(4, 2, 1^2) &  ((1^2), (2))  &  
{q^2} & \cdot  & q & \cdot  & \cdot  & {q^ 2} & q & 1 & \cdot  & \cdot  & \cdot  & \cdot  & \cdot  & \cdot  & \cdot  & \cdot  \cr 
(4, 1^4) &  (\vr, (2, 1^2))  & 
{q^2} & {q^ 3} & \cdot  & \cdot  & {q^2} & \cdot  & q & \cdot  & 1 & \cdot  & \cdot  & \cdot  & \cdot  & \cdot  & \cdot  & \cdot  \cr 
(3^2, 2) &  ((2), (1^2))  &       
{q^ 2} & \cdot  & q & \cdot  & \cdot  & {q^2} & q & \cdot  & \cdot  & 1 & \cdot  & \cdot  & \cdot  & \cdot  & \cdot  & \cdot  \cr 
(3^2, 1^2) &  ((2, 1), (1))  &     
{q^3} & {q^2} & 2\,{q^2} & q & q & {q^3} & {q^2} & q & \cdot  & q & 1 & \cdot  & \cdot  & \cdot  & \cdot  & \cdot  \cr 
(3, 2^2, 1) &  ((1^2), (1^2))  &     
{q^ 2} & \cdot  & {q^3} & \cdot  & {q^2} & {q^4} & q + {q^3} & {q^2} & \cdot  & {q^2} & q & 1 & \cdot  & \cdot  & \cdot  & \cdot  \cr 
(3, 1^5) &  (\vr, (1^4))  &     
{q^3} & \cdot  & \cdot  & \cdot  & \cdot  & \cdot  & {q^2} & \cdot  & q & \cdot  & \cdot  & q & 1 & \cdot  & \cdot  & \cdot  \cr 
(2^4) &  ((2^2), \vr)  &          
\cdot  & {q^3} & \cdot  & {q^2} & {q^2} & \cdot  & \cdot  & \cdot  & \cdot  & \cdot  & q & \cdot  & \cdot  & 1 & \cdot  & \cdot  \cr 
(2, 1^6) &  ((2, 1^2), \vr)  &     
{q^3} & {q^4} & {q^2} & {q^3} & {q^3} & \cdot  & {q^2} & q & q & q & {q^2} & \cdot  & \cdot  & q & 1 & \cdot  \cr 
(1^8) &  ((1^4), \vr)  &         
{q^4} & \cdot  & {q^3} & \cdot  & \cdot  & \cdot  & {q^3} & {q^2} & {q^2} & {q^2} & q & {q^2} & q & \cdot  & q & 1  
\end{array}
\end{align*}

\begin{align*}
&\text{ $\wt = 2\Lambda_0 - 3\alpha_0 - \alpha_1$}\\
&\begin{array}{l|l| c c c c  }
\!\!\!\! *\:(7, 2, 1) &  ((1), (3)) &  1 & \cdot &  \cdot & \cdot \cr 
(5, 4, 1) &  ((3), (1)) &  q & 1 &  \cdot & \cdot \cr 
(3, 2^3, 1) &  ((1), (1^3))  & q & \cdot &  1 & \cdot \cr 
(3, 2, 1^5) &  ((1^3), (1)) &  {q^2} & q &  q & 1  
\end{array}
\end{align*}

\mbox{} \\
\vspace{10mm}
\mbox{}

\noindent $k=5$
\begin{align*}
&\text{ $\wt = 2\Lambda_0 - 2\alpha_0 - 3\alpha_1$} \\
&\begin{array}{l | l | c c c c c c c c }
\!\!\!\! *\:(8, 1) & (\vr, (4, 1)) &  1 & \cdot & \cdot & \cdot & \cdot & \cdot & \cdot & \cdot \cr 
(6, 1^3) &  ((4, 1), \vr) & q & 1 & \cdot & \cdot & \cdot & \cdot & \cdot & \cdot \cr 
\!\!\!\! *\:(4, 3, 2) &  ((2), (2, 1))& q & \cdot & 1 & \cdot & \cdot & \cdot & \cdot & \cdot \cr 
(4, 3, 1^2) & ((2, 1), (2)) & {q^2} & q & q & 1 & \cdot & \cdot & \cdot & \cdot \cr 
(4, 2^2, 1) &  ((1^2), (2, 1)) & q & \cdot & {q^2} & q & 1 & \cdot & \cdot & \cdot \cr 
(4, 1^5) & (\vr, (2, 1^3))  & {q^2} & \cdot & \cdot & \cdot & q & 1 & \cdot & \cdot \cr 
(3^2, 2, 1) & ((2, 1), (1^2)) &  {q^2} & q & {q^3} & {q^2} & q & \cdot & 1 & \cdot \cr 
(2, 1^7) &  ((2, 1^3), \vr) & {q^3} & {q^2} & \cdot & q & {q^2} & q & q & 1  
\end{array}
\end{align*}
%\end{document}

%\hoffset=-1.in
\begin{align*}
&\!\!\!\!\!\!\!\!\text{ $\wt = 2\Lambda_0 - 3\alpha_0 - 2\alpha_1$ $(\dim =28)$}\\
&\!\!\!\!\!\!\!\!\begin{array}{l|l| c c c c c c c c c c c c c c}
\!\!\!\! *\: ( 11)  & ( \vr ,( 5) )  & 1 & \cdot & \cdot & \cdot & \cdot & \cdot & \cdot & \cdot & \cdot & \cdot & \cdot & \cdot & \cdot & \cdot \cr 
( 9,1^2)  & ( ( 5) , \vr )  & q & 1 & \cdot & \cdot & \cdot & \cdot & \cdot & \cdot & \cdot & \cdot & \cdot & \cdot & \cdot & \cdot \cr 
\!\!\!\! *\:( 8,   2,1)  & ( ( 1) , ( 4) )  & q & \cdot & 1 & \cdot & \cdot & \cdot & \cdot & \cdot & \cdot & \cdot & \cdot & \cdot & \cdot & \cdot \cr 
\!\!\!\! *\:( 7,4)  & ( \vr , ( 3,2) )  & \cdot & \cdot & \cdot & 1 & \cdot & \cdot & \cdot & \cdot & \cdot & \cdot & \cdot & \cdot & \cdot & \cdot \cr 
\!\!\!\! *\:( 7,3,1)  & ( ( 2) ,( 3) )  & {q^2} & q & q & \cdot & 1 & \cdot & \cdot & \cdot & \cdot & \cdot & \cdot & \cdot & \cdot & \cdot \cr 
( 7,2^2)  & ( ( 1) ,( 3,1) )  & q & \cdot & {q^2} & q & q & 1 & \cdot & \cdot & \cdot & \cdot & \cdot & \cdot & \cdot & \cdot \cr 
( 7,2,1^2)  & (( 1^2) ,( 3) )  & {q^2} & q & {q^3} & \cdot & {q^2} & q & 1 & \cdot & \cdot & \cdot & \cdot & \cdot & \cdot & \cdot \cr 
( 7,1^4)  & ( \vr ,( 3,1^2) )  & {q^2} & \cdot & \cdot & {q^2} & \cdot & q & \cdot & 1 & \cdot & \cdot & \cdot & \cdot & \cdot & \cdot \cr 
( 6,4,1)  & ( ( 4) , ( 1) )  & {q^ 2} & q & q & \cdot & \cdot & \cdot & \cdot & \cdot & 1 & \cdot & \cdot & \cdot & \cdot & \cdot \cr 
( 5,5,1)  & ( ( 3) ,( 2) )  & {q^3} & {q^2} & {q^2} & \cdot & q & \cdot & \cdot & \cdot & q & 1 & \cdot & \cdot & \cdot & \cdot \cr 
( 5,4,2)  & ( ( 3) ,( 1^2) )  & {q^2} & q & {q^3} & \cdot & {q^2} & q & \cdot & \cdot & {q^2} & q & 1 & \cdot & \cdot & \cdot \cr 
( 5,4,1^2)  & ( ( 3,1) ,( 1) )  & {q^3} & 2\,{q^2} & {q^4} & q & {q^3} & {q^2} & q & \cdot & {q^3} & {q^2} & q & 1 & \cdot & \cdot \cr ( 5,2^3)  & ( ( 3,2) ,
   \vr )  & \cdot & \cdot & \cdot & {q^
    2} & \cdot & \cdot & \cdot & \cdot & \cdot & \cdot & \cdot & q & 1 & \cdot \cr ( 5,1^6)  & ( 
   ( 3,1^2) ,\vr )  & {q^3} & {q^2} & \cdot & {q^3} & \cdot & {q^
    2} & q & q & \cdot & \cdot & q & {q^2} & q & 1 \cr ( 4^2,3)  & ( ( 1) ,
   ( 2^2) )  & \cdot & \cdot & \cdot & {q^2} & {q^
    2} & q & \cdot & \cdot & \cdot & q & \cdot & \cdot & \cdot & \cdot \cr ( 4^2,1^3)  & ( \vr ,
   ( 2^2,1) )  & \cdot & \cdot & \cdot & {q^3} & \cdot & {q^
    2} & \cdot & q & \cdot & \cdot & \cdot & \cdot & \cdot & \cdot \cr ( 4,2^3,1)  & ( ( 1) ,
   ( 2,1^2) )  & {q^2} & \cdot & q & {q^4} & {q^2} & q + {q^3} & \cdot & {q^
    2} & \cdot & q & \cdot & \cdot & \cdot & \cdot \cr ( 4,2,1^5)  & ( ( 1^3) ,
   ( 2) )  & {q^3} & {q^2} & {q^2} & \cdot & {q^3} & {q^2} & q & \cdot & q & {q^
    2} & \cdot & \cdot & \cdot & \cdot \cr ( 3^3,2)  & ( ( 2^2) ,( 1) )  & \cdot & {q^
    2} & \cdot & {q^3} & {q^3} & {q^2} & q & \cdot & \cdot & {q^2} & q & {q^
    2} & q & \cdot \cr ( 3^2,2^2,1)  & ( ( 2) ,( 1^3) )  & {q^3} & {q^
    2} & {q^2} & \cdot & {q^3} & {q^2} & \cdot & \cdot & \cdot & {q^
    2} & q & \cdot & \cdot & \cdot \cr ( 3^2,1^5)  & ( ( 2,1^2) ,
   ( 1) )  & {q^4} & 2\,{q^3} & {q^3} & {q^4} & {q^4} & 2\,{q^3} & 2\,
   {q^2} & {q^2} & {q^2} & q + {q^3} & 2\,{q^2} & {q^3} & {q^2} & q \cr ( 3,
   2^4)  & ( ( 1^2) ,( 1^3) )  & {q^2} & {q^3} & {q^
    3} & \cdot & {q^4} & q + {q^3} & {q^2} & \cdot & \cdot & {q^3} & {q^
    2} & \cdot & \cdot & \cdot \cr ( 3,2^3,1^2)  & ( ( 1) ,( 1^4) )  & {q^
    3} & \cdot & {q^4} & \cdot & \cdot & {q^2} & \cdot & q & \cdot & \cdot & \cdot & \cdot & \cdot & \cdot \cr ( 3,2^2,1^4)  & ( ( 1^3) ,( 1^2) )  & {q^3} & {q^4} & {q^
    4} & \cdot & {q^5} & {q^2} + {q^4} & {q^3} & \cdot & {q^3} & {q^2} + {q^4} & {q^
    3} & \cdot & \cdot & \cdot \cr ( 3,2,1^6)  & ( ( 1^4) ,
   ( 1) )  & {q^4} & {q^3} & {q^5} & \cdot & \cdot & {q^3} & {q^2} & {q^2} & {q^
    4} & {q^3} & {q^2} & \cdot & \cdot & q \cr ( 3,1^8)  & ( \vr ,
   ( 1^5) )  & {q^4} & \cdot & \cdot & \cdot & \cdot & {q^3} & \cdot & {q^
    2} & \cdot & \cdot & \cdot & \cdot & \cdot & \cdot \cr ( 2^4,1^3)  & ( ( 2^2,1) ,
   \vr )  & \cdot & \cdot & \cdot & {q^5} & \cdot & {q^4} & {q^3} & {q^3} & \cdot & {q^
    2} & {q^3} & {q^4} & q + {q^3} & {q^2} \cr ( 1^{11})  & ( ( 1^5) ,\vr )  & {q^5} & {q^4} & \cdot & \cdot & \cdot & {q^
    4} & {q^3} & {q^3} & \cdot & \cdot & {q^3} & \cdot & \cdot & {q^2}  
\end{array}
\end{align*}
%\noindent $\bullet$ ${\mathrm {wt}} = 2\Lambda_0 - 3\alpha_0 - 2\alpha_1$ (continued).
\begin{align*}
&\text{ $\wt = 2\Lambda_0 - 3\alpha_0 - 2\alpha_1$ (continued)} \\
&\begin{array}{l|l| c c c c c c c c c c c c c c}
( 4^2,3)  & ( ( 1) ,
   ( 2^2) )  & 1 & \cdot & \cdot & \cdot & \cdot & \cdot & \cdot & \cdot & \cdot & \cdot & \cdot & \cdot & \cdot & \cdot 
\cr ( 4^2,1^3)  & ( \vr ,
   ( 2^2,1) )  & q & 1 & \cdot & \cdot & \cdot & \cdot & \cdot & \cdot & \cdot & \cdot & \cdot & \cdot & \cdot & \cdot 
\cr ( 4,2^3,1)  & ( ( 1) ,( 2,1^2) )  & {q^
    2} & q & 1 & \cdot & \cdot & \cdot & \cdot & \cdot & \cdot & \cdot & \cdot & \cdot & \cdot & \cdot \cr ( 4,2,1^5)  & ( ( 1^3) ,
   ( 2) )  & \cdot & \cdot & q & 1 & \cdot & \cdot & \cdot & \cdot & \cdot & \cdot & \cdot & \cdot & \cdot & \cdot 
\cr ( 3^3,2)  & ( ( 2^2) ,
   ( 1) )  & q & \cdot & \cdot & \cdot & 1 & \cdot & \cdot & \cdot & \cdot & \cdot & \cdot & \cdot & \cdot & \cdot 
\cr ( 3^2,2^2,1)  & ( ( 2) ,
   ( 1^3) )  & \cdot & \cdot & q & \cdot & \cdot & 1 & \cdot & \cdot & \cdot & \cdot & \cdot & \cdot & \cdot & \cdot 
\cr ( 3^2,1^5)  & ( ( 2,1^2) ,( 1) )  & {q^2} & q & {q^
    2} & q & q & q & 1 & \cdot & \cdot & \cdot & \cdot & \cdot & \cdot & \cdot \cr ( 3,2^4)  & ( 
   ( 1^2) ,( 1^3) )  & {q^2} & q & {q^
    2} & \cdot & q & q & \cdot & 1 & \cdot & \cdot & \cdot & \cdot & \cdot & \cdot \cr ( 3,2^3,1^2)  & ( ( 1) ,( 1^4) )  & \cdot & {q^2} & {q^3} & \cdot & \cdot & {q^
    2} & \cdot & q & 1 & \cdot & \cdot & \cdot & \cdot & \cdot \cr ( 3,2^2,1^4)  & ( 
   ( 1^3) ,( 1^2) )  & {q^3} & {q^2} & {q^3} & {q^2} & {q^2} & {q^
    2} & q & q & \cdot & 1 & \cdot & \cdot & \cdot & \cdot \cr ( 3,2,1^6)  & ( 
   ( 1^4) ,( 1) )  & \cdot & {q^3} & {q^4} & {q^3} & q & {q^3} & {q^
    2} & {q^2} & q & q & 1 & \cdot & \cdot & \cdot \cr ( 3,1^8)  & ( 
   \vr ,( 1^5) )  & \cdot & q & \cdot & \cdot & \cdot & \cdot & \cdot & {q^
    2} & q & q & \cdot & 1 & \cdot & \cdot \cr ( 2^4,1^3)  & ( ( 2^2,1) ,
   \vr )  & {q^3} & {q^2} & \cdot & \cdot & {q^
    2} & \cdot & q & \cdot & \cdot & \cdot & \cdot & \cdot & 1 & \cdot \cr ( 1^{11})  & ( ( 1^5) ,\vr )  & \cdot & {q^2} & \cdot & \cdot & {q^2} & {q^
    2} & q & {q^3} & {q^2} & {q^2} & q & q & \cdot & 1  
\end{array}
\end{align*}
%\newpage
%\input{G12}

%\hoffset=-.5in

\bigskip

\noindent $k=6$
%\noindent$\bullet$ ${\mathrm {wt}}=2\Lambda_0 - 2\alpha_0-4\alpha_1.$ $\dim = 1.$ 
\begin{align*}
&\text{ $\wt=2\Lambda_0 - 2\alpha_0-4\alpha_1$}\\ 
&\begin{array}{l|l| c}
\!\!\!\! *\:(4,3,2,1) & ((2,1),(2,1)) & 1 
\end{array}
\end{align*}
%\noindent$\bullet$ ${\mathrm {wt}}=2\Lambda_0 - 4\alpha_0-2\alpha_1.$ $\dim = 16.$ 
%\input{G14}
\begin{align*}
&\text{ $\wt=2\Lambda_0 - 4\alpha_0-2\alpha_1$ }\\ 
&\begin{array}{l|l| c c c c c c c c c c c c c c c c} 
\!\!\!\! *\:( 11,2,1)  & ( ( 1) ,
   ( 5) )  & 1 & \cdot & \cdot & \cdot & \cdot & \cdot & \cdot & \cdot & \cdot & \cdot & \cdot & \cdot & \cdot & \cdot & \cdot & 
\cdot \cr ( 9,4,1)  & ( ( 5) ,
   ( 1) )  & q & 1 & \cdot & \cdot & \cdot & \cdot & \cdot & \cdot & \cdot & \cdot & \cdot & \cdot & \cdot & \cdot & \cdot & 
\cdot \cr \!\!\!\! *\:( 7,6,1)  & ( ( 3) ,( 3) )  & {q^
    2} & q & 1 & \cdot & \cdot & \cdot & \cdot & \cdot & \cdot & \cdot & \cdot & \cdot & \cdot & \cdot & \cdot & \cdot \cr \!\!\!\! *\:( 7,4,3)  & ( ( 1) ,
   ( 3,2) )  & \cdot & \cdot & q & 1 & \cdot & \cdot & \cdot & \cdot & \cdot & \cdot & \cdot & \cdot & \cdot & \cdot & \cdot 
& \cdot \cr ( 7,4,1^3)  & ( \vr ,
   ( 3,2,1) )  & \cdot & \cdot & \cdot & q & 1 & \cdot & \cdot & \cdot & \cdot & \cdot & \cdot & \cdot & \cdot & \cdot & 
\cdot & \cdot \cr ( 7,2^3,1)  & ( ( 1) ,( 3,1^2) )  & q & \cdot & q & {q^
    2} & q & 1 & \cdot & \cdot & \cdot & \cdot & \cdot & \cdot & \cdot & \cdot & \cdot & \cdot \cr ( 7,2,1^5)  & ( ( 1^3) ,( 3) )  & {q^2} & q & {q^
    2} & \cdot & \cdot & q & 1 & \cdot & \cdot & \cdot & \cdot & \cdot & \cdot & \cdot & \cdot & \cdot \cr ( 5,4,3,2)  & ( ( 3,2) ,( 1) )  & \cdot & q & {q^
    2} & q & \cdot & \cdot & \cdot & 1 & \cdot & \cdot & \cdot & \cdot & \cdot & \cdot & \cdot & \cdot \cr ( 5,4,2^2,1)  & ( ( 3) ,( 1^3) )  & {q^2} & q & {q^
    2} & \cdot & \cdot & q & \cdot & \cdot & 1 & \cdot & \cdot & \cdot & \cdot & \cdot & \cdot & \cdot \cr ( 5,4,1^5)  & ( ( 3,1^2) ,( 1) )  & {q^3} & 2\,{q^2} & {q^3} & {q^
    2} & q & {q^2} & q & q & q & 1 & \cdot & \cdot & \cdot & \cdot & \cdot & \cdot \cr ( 5,2^3,1^3)  & ( ( 3,2,1) ,\vr )  & \cdot & \cdot & \cdot & {q^3} & {q^2} & \cdot & \cdot & {q^
    2} & \cdot & q & 1 & \cdot & \cdot & \cdot & \cdot & \cdot \cr ( 4^2,3,2,1)  & ( ( 1) ,
   ( 2^2,1) )  & \cdot & \cdot & {q^2} & {q^3} & {q^
    2} & q & \cdot & \cdot & \cdot & \cdot & \cdot & 1 & \cdot & \cdot & \cdot & \cdot \cr ( 3^3,2,1^3)  & ( ( 2^2,1) ,( 1) )  & \cdot & {q^2} & {q^3} & {q^4} & {q^
    3} & {q^2} & q & {q^3} & q & {q^2} & q & q & 1 & \cdot & \cdot & \cdot \cr ( 3,2^5,1)  & ( ( 1^3) ,( 1^3) )  & {q^2} & {q^3} & {q^
    4} & \cdot & \cdot & q + {q^3} & {q^2} & \cdot & {q^2} & \cdot & \cdot & {q^
    2} & q & 1 & \cdot & \cdot \cr ( 3,2^3,1^5)  & ( ( 1) ,
   ( 1^5) )  & {q^3} & \cdot & \cdot & \cdot & q & {q^
    2} & \cdot & \cdot & \cdot & \cdot & \cdot & \cdot & \cdot & q & 1 & \cdot \cr ( 3,2,1^9)  & ( ( 1^5) ,( 1) )  & {q^4} & {q^3} & \cdot & \cdot & {q^
    2} & {q^3} & {q^2} & \cdot & {q^2} & q & \cdot & \cdot & q & {q^2} & q & 1  
\end{array}
\end{align*}

%\noindent$\bullet$ ${\mathrm {wt}}=2\Lambda_0 - 3\alpha_0-3\alpha_1.$ $\dim = 48.$ 

\newpage \thispagestyle{empty}
\begin{align*}
&\!\!\!\!\!\!\!\!\text{ $\wt =2\Lambda_0 - 3\alpha_0-3\alpha_1$ $(\dim = 48)$} \\
&\!\!\!\!\!\!\!\!\begin{array}{l|l| c c c c c c c c c c c c}
   \!\!\!\! *\:( 12)  & ( \vr ,
   ( 6) )  & 1 & \cdot & \cdot & \cdot & \cdot & \cdot & \cdot & \cdot & \cdot & \cdot & \cdot & \cdot \cr \!\!\!\! *\:( 11,
   1)  & ( \vr ,
   ( 5,1) )  & q & 1 & \cdot & \cdot & \cdot & \cdot & \cdot & \cdot & \cdot & \cdot & \cdot & \cdot \cr ( 10,1^2)  & ( ( 6) ,
   \vr )  & q & \cdot & 1 & \cdot & \cdot & \cdot & \cdot & \cdot & \cdot & \cdot & \cdot & \cdot \cr ( 9,1^3)  & ( ( 5,1) ,\vr )  & {q^
    2} & q & q & 1 & \cdot & \cdot & \cdot & \cdot & \cdot & \cdot & \cdot & \cdot \cr \!\!\!\! *\:( 8,4)  & ( \vr ,
   ( 4,2) )  & \cdot & q & \cdot & \cdot & 1 & \cdot & \cdot & \cdot & \cdot & \cdot & \cdot & \cdot \cr \!\!\!\! *\:( 8,3,
   1)  & ( ( 2) ,( 4) )  & {q^
    2} & \cdot & q & \cdot & \cdot & 1 & \cdot & \cdot & \cdot & \cdot & \cdot & \cdot \cr ( 8,2^2)  & ( 
   ( 1) ,( 4,1) )  & q & {q^
    2} & \cdot & \cdot & q & q & 1 & \cdot & \cdot & \cdot & \cdot & \cdot \cr ( 8,2,1^2)  & ( 
   ( 1^2) ,( 4) )  & {q^2} & \cdot & q & \cdot & \cdot & {q^
    2} & q & 1 & \cdot & \cdot & \cdot & \cdot \cr ( 8,1^4)  & ( \vr ,
   ( 4,1^2) )  & {q^2} & {q^3} & \cdot & \cdot & {q^
    2} & \cdot & q & \cdot & 1 & \cdot & \cdot & \cdot \cr ( 7,5)  & ( \vr ,
   ( 3^2) )  & \cdot & \cdot & \cdot & \cdot & q & \cdot & \cdot & \cdot & \cdot & 1 & \cdot & \cdot \cr \!\!\!\! *\:( 7,3,
   2)  & ( ( 2) ,( 3,1) )  & {q^2} & q & q & \cdot & {q^2} & {q^
    2} & q & \cdot & \cdot & q & 1 & \cdot \cr ( 7,3,1^2)  & ( ( 2,1) ,
   ( 3) )  & {q^3} & {q^2} & 2\,{q^2} & q & \cdot & {q^3} & {q^
    2} & q & \cdot & \cdot & q & 1 \cr ( 7,2^2,1)  & ( ( 1^2) ,
   ( 3,1) )  & {q^2} & q & {q^3} & \cdot & {q^2} & {q^4} & q + {q^3} & {q^
    2} & \cdot & q & {q^2} & q \cr ( 7,1^5)  & ( \vr ,
   ( 3,1^3) )  & {q^3} & {q^2} & \cdot & \cdot & {q^3} & \cdot & {q^
    2} & \cdot & q & {q^2} & \cdot & \cdot \cr ( 6,5,1)  & ( ( 4) ,( 2) )  & {q^
    3} & \cdot & {q^2} & \cdot & \cdot & q & \cdot & \cdot & \cdot & \cdot & \cdot & \cdot \cr ( 6,4,2)  & ( 
   ( 4) ,( 1^2) )  & {q^2} & \cdot & q & \cdot & \cdot & {q^
    2} & q & \cdot & \cdot & \cdot & \cdot & \cdot \cr ( 6,4,1^2)  & ( ( 4,1) ,
   ( 1) )  & {q^3} & {q^2} & 2\,{q^2} & q & q & {q^3} & {q^
    2} & q & \cdot & \cdot & \cdot & \cdot \cr ( 6,2^3)  & ( ( 4,2) ,
   \vr )  & \cdot & {q^3} & \cdot & {q^2} & {q^
    2} & \cdot & \cdot & \cdot & \cdot & \cdot & \cdot & \cdot \cr ( 6,1^6)  & ( ( 4,1^2) ,
   \vr )  & {q^3} & {q^4} & {q^2} & {q^3} & {q^3} & \cdot & {q^
    2} & q & q & \cdot & \cdot & \cdot \cr ( 5^2,2)  & ( ( 3) ,( 2,1) )  & {q^
    3} & {q^2} & {q^2} & \cdot & \cdot & {q^3} & {q^2} & \cdot & \cdot & \cdot & q & \cdot \cr ( 5^2,
   1^2)  & ( ( 3,1) ,( 2) )  & {q^4} & {q^3} & 2\,{q^3} & {q^2} & {q^
    2} & {q^4} & {q^3} & {q^2} & \cdot & q & {q^2} & q \cr ( 5,4,2,1)  & ( 
   ( 3,1) ,( 1^2) )  & {q^3} & {q^2} & {q^2} + {q^4} & q & {q^3} & {q^
    5} & {q^2} + {q^4} & {q^3} & \cdot & {q^2} & {q^3} & {q^2} \cr ( 5,3,2^2)  & ( ( 3^2) ,\vr )  & \cdot & \cdot & \cdot & \cdot & {q^
    3} & \cdot & \cdot & \cdot & \cdot & {q^2} & \cdot & \cdot \cr ( 5,1^7)  & ( 
   ( 3,1^3) ,\vr )  & {q^4} & {q^3} & {q^3} & {q^2} & {q^4} & \cdot & {q^
    3} & {q^2} & {q^2} & {q^3} & \cdot & q \cr ( 4^3)  & ( ( 2) ,
   ( 2^2) )  & \cdot & {q^3} & \cdot & \cdot & {q^2} & \cdot & q & \cdot & \cdot & q & {q^
    2} & \cdot \cr ( 4^2,3,1)  & ( ( 1^2) ,( 2^2) )  & \cdot & {q^
    2} & \cdot & \cdot & {q^3} & \cdot & {q^2} & \cdot & \cdot & 2\,{q^2} & {q^3} & {q^
    2} \cr ( 4^2,2,1^2)  & ( \vr ,( 2^3) )  & \cdot & \cdot & \cdot & \cdot & {q^
    3} & \cdot & {q^2} & \cdot & q & {q^2} & \cdot & \cdot \cr ( 4^2,1^4)  & ( \vr ,
   ( 2^2,1^2) )  & \cdot & {q^3} & \cdot & \cdot & {q^4} & \cdot & {q^3} & \cdot & {q^
    2} & {q^3} & \cdot & \cdot \cr ( 4,3^2,2)  & ( ( 2^2) ,( 2) )  & \cdot & {q^
    4} & {q^2} & {q^3} & {q^3} & \cdot & {q^2} & q & \cdot & {q^2} & {q^3} & {q^
    2} \cr ( 4,3,2^2,1)  & ( ( 2) ,( 2,1^2) )  & {q^3} & {q^4} & {q^
    2} & \cdot & {q^3} & q & 2\,{q^2} & \cdot & q & {q^2} & {q^3} & \cdot \cr ( 4,3,1^5)  & ( ( 2,1^2) ,( 2) )  & {q^4} & {q^5} & 2\,{q^3} & {q^
    4} & {q^4} & {q^2} & 2\,{q^3} & 2\,{q^2} & {q^2} & {q^3} & {q^4} & {q^
    3} \cr ( 4,2^4)  & ( ( 1^2) ,( 2,1^2) )  & {q^2} & {q^
    3} & {q^3} & \cdot & {q^4} & {q^2} & q + 2\,{q^3} & {q^2} & {q^2} & 2\,
   {q^3} & {q^4} & {q^3} \cr ( 4,2^3,1^2)  & ( ( 1) ,
   ( 2,1^3) )  & {q^3} & {q^4} & \cdot & \cdot & {q^5} & {q^3} & {q^2} + 
   {q^4} & \cdot & q + {q^3} & {q^4} & \cdot & \cdot \cr ( 4,2^2,1^4)  & ( 
   ( 1^3) ,( 2,1) )  & {q^3} & {q^4} & {q^4} & \cdot & \cdot & {q^3} & {q^2} + 
   {q^4} & {q^3} & \cdot & {q^4} & {q^5} & {q^4} \cr ( 4,2,1^6)  & ( 
   ( 1^4) ,( 2) )  & {q^4} & \cdot & {q^3} & \cdot & \cdot & {q^4} & {q^
    3} & {q^2} & {q^2} & \cdot & \cdot & \cdot \cr ( 4,1^8)  & ( \vr ,
   ( 2,1^4) )  & {q^4} & {q^5} & \cdot & \cdot & \cdot & \cdot & {q^3} & \cdot & {q^
    2} & \cdot & \cdot & \cdot \cr ( 3^4)  & ( ( 2^2) ,( 1^2) )  & \cdot & {q^
    3} & {q^3} & {q^2} & {q^4} & \cdot & {q^3} & {q^2} & \cdot & 2\,{q^3} & {q^
    4} & {q^3} \cr ( 3^2,2^3)  & ( ( 2,1) ,( 1^3) )  & {q^
    3} & {q^4} & {q^2} + {q^4} & {q^3} & \cdot & {q^3} & {q^2} + {q^4} & {q^
    3} & \cdot & {q^4} & {q^5} & {q^4} \cr ( 3^2,2^2,1^2)  & ( ( 2) ,
   ( 1^4) )  & {q^4} & \cdot & {q^3} & \cdot & \cdot & {q^4} & {q^3} & \cdot & {q^
    2} & \cdot & \cdot & \cdot \cr ( 3^2,2,1^4)  & ( ( 2,1^2) ,
   ( 1^2) )  & {q^4} & {q^5} & {q^3} + {q^5} & {q^4} & {q^4} & {q^
    4} & 2\,{q^3} + {q^5} & {q^2} + {q^4} & {q^2} & {q^3} + {q^5} & {q^
    6} & {q^5} \cr ( 3^2,1^6)  & ( ( 2,1^3) ,( 1) )  & {q^
    5} & {q^4} & 2\,{q^4} & {q^3} & {q^5} & {q^5} & 2\,{q^4} & 2\,
   {q^3} & 2\,{q^3} & {q^4} & \cdot & {q^2} \cr ( 3,2^4,1)  & ( ( 1^2) ,
   ( 1^4) )  & {q^3} & \cdot & {q^4} & \cdot & \cdot & {q^5} & {q^2} + 
   {q^4} & {q^3} & {q^3} & \cdot & \cdot & \cdot \cr ( 3,2^2,1^5)  & ( 
   ( 1^4) ,( 1^2) )  & {q^4} & \cdot & {q^5} & \cdot & \cdot & {q^6} & {q^3} + 
   {q^5} & {q^4} & {q^4} & \cdot & \cdot & \cdot \cr ( 3,1^9)  & ( 
   \vr ,( 1^6) )  & {q^5} & \cdot & \cdot & \cdot & \cdot & \cdot & {q^4} & \cdot & {q^
    3} & \cdot & \cdot & \cdot \cr ( 2^5,1^2)  & ( ( 2^3) ,
   \vr )  & \cdot & \cdot & \cdot & \cdot & {q^5} & \cdot & {q^4} & {q^3} & {q^3} & {q^
    4} & \cdot & \cdot \cr ( 2^4,1^4)  & ( ( 2^2,1^2) ,
   \vr )  & \cdot & {q^5} & \cdot & {q^4} & {q^6} & \cdot & {q^5} & {q^4} & {q^
    4} & {q^5} & \cdot & {q^3} \cr ( 2,1^{10})  & ( 
   ( 2,1^4) ,\vr )  & {q^5} & {q^6} & {q^4} & {q^5} & \cdot & \cdot & {q^
    4} & {q^3} & {q^3} & \cdot & \cdot & {q^4} \cr ( 1^{12})  & ( ( 1^6) ,\vr )  & {q^6} & \cdot & {q^
    5} & \cdot & \cdot & \cdot & {q^5} & {q^4} & {q^4} & \cdot & \cdot & \cdot  
\end{array}
\end{align*}
%\noindent$\bullet$ ${\mathrm {wt}}=2\Lambda_0 - 3\alpha_0-3\alpha_1$ (continued). 
\begin{align*}
&\!\!\!\!\!\!\!\!\text{ $\wt =2\Lambda_0 - 3\alpha_0-3\alpha_1$ (continued)} \\
&\!\!\!\!\!\!\!\!\begin{array}{l|l| c c c c c c c c c c c c}
 ( 7,2^2,1)  & ( ( 1^2) ,
   ( 3,1) )  & 1 & \cdot & \cdot & \cdot & \cdot & \cdot & \cdot & \cdot & \cdot & \cdot & \cdot & \cdot \cr ( 7,1^5)  & ( \vr ,
   ( 3,1^3) )  & q & 1 & \cdot & \cdot & \cdot & \cdot & \cdot & \cdot & \cdot & \cdot & \cdot & \cdot \cr ( 6, 5,1)  & ( ( 4) ,
   ( 2) )  & \cdot & \cdot & 1 & \cdot & \cdot & \cdot & \cdot & \cdot & \cdot & \cdot & \cdot & \cdot \cr ( 6,4,2)  & ( ( 4) ,
   ( 1^2) )  & \cdot & \cdot & q & 1 & \cdot & \cdot & \cdot & \cdot & \cdot & \cdot & \cdot & \cdot \cr ( 6,4,1^2)  & ( ( 4,1) ,( 1) )  & \cdot & \cdot & {q^
    2} & q & 1 & \cdot & \cdot & \cdot & \cdot & \cdot & \cdot & \cdot \cr ( 6,2^3)  & ( ( 4,2) ,
   \vr )  & \cdot & \cdot & \cdot & \cdot & q & 1 & \cdot & \cdot & \cdot & \cdot & \cdot & \cdot \cr ( 6,1^6)  & ( ( 4,1^2) ,\vr )  & \cdot & \cdot & \cdot & q & {q^
    2} & q & 1 & \cdot & \cdot & \cdot & \cdot & \cdot \cr ( 5^2,2)  & ( ( 3) ,
   ( 2,1) )  & \cdot & \cdot & {q^2} & q & \cdot & \cdot & \cdot & 1 & \cdot & \cdot & \cdot & \cdot \cr ( 5^2,1^2)  & ( ( 3,1) ,( 2) )  & \cdot & \cdot & {q^3} & {q^
    2} & q & \cdot & \cdot & q & 1 & \cdot & \cdot & \cdot \cr ( 5,4,2,1)  & ( ( 3,1) ,
   ( 1^2) )  & q & \cdot & {q^4} & q + {q^3} & {q^2} & \cdot & \cdot & {q^
    2} & q & 1 & \cdot & \cdot \cr ( 5,3,2^2)  & ( ( 3^2) ,
   \vr )  & \cdot & \cdot & \cdot & \cdot & {q^2} & q & \cdot & \cdot & q & \cdot & 1 & \cdot \cr ( 5,1^7)  & ( ( 3,1^3) ,\vr )  & {q^2} & q & \cdot & {q^2} & {q^
    3} & {q^2} & q & \cdot & {q^2} & q & q & 1 \cr ( 4^3)  & ( ( 2) ,
   ( 2^2) )  & \cdot & \cdot & \cdot & \cdot & \cdot & \cdot & \cdot & q & \cdot & \cdot & \cdot & \cdot \cr ( 4^2,3,1)  & ( ( 1^2) ,( 2^2) )  & q & \cdot & \cdot & \cdot & \cdot & \cdot & \cdot & {q^
    2} & q & \cdot & \cdot & \cdot \cr ( 4^2,2,1^2)  & ( \vr ,
   ( 2^3) )  & \cdot & \cdot & \cdot & \cdot & \cdot & \cdot & \cdot & \cdot & \cdot & \cdot & \cdot & \cdot \cr ( 4^2,1^4)  & ( \vr ,( 2^2,1^2) )  & {q^
    2} & q & \cdot & \cdot & \cdot & \cdot & \cdot & \cdot & \cdot & \cdot & \cdot & \cdot \cr ( 4,3^2,2)  & ( 
   ( 2^2) ,( 2) )  & \cdot & \cdot & \cdot & q & {q^2} & q & \cdot & {q^
    2} & q & \cdot & \cdot & \cdot \cr ( 4,3,2^2,1)  & ( ( 2) ,
   ( 2,1^2) )  & \cdot & \cdot & \cdot & q & \cdot & \cdot & \cdot & {q^
    2} & \cdot & \cdot & \cdot & \cdot \cr ( 4,3,1^5)  & ( ( 2,1^2) ,
   ( 2) )  & \cdot & \cdot & q & 2\,{q^2} & {q^3} & {q^2} & q & {q^3} & {q^
    2} & \cdot & \cdot & \cdot \cr ( 4,2^4)  & ( ( 1^2) ,( 2,1^2) )  & {q^
    2} & \cdot & \cdot & {q^2} & \cdot & \cdot & \cdot & {q^3} & {q^2} & \cdot & \cdot & \cdot \cr ( 4,2^3,
   1^2)  & ( ( 1) ,( 2,1^3) )  & {q^3} & {q^
    2} & \cdot & \cdot & \cdot & \cdot & \cdot & \cdot & \cdot & \cdot & \cdot & \cdot \cr ( 4,2^2,1^4)  & ( 
   ( 1^3) ,( 2,1) )  & {q^3} & \cdot & {q^2} & {q^3} & \cdot & \cdot & \cdot & {q^
    4} & {q^3} & \cdot & \cdot & \cdot \cr ( 4,2,1^6)  & ( ( 1^4) ,
   ( 2) )  & \cdot & \cdot & {q^3} & {q^
    2} & \cdot & \cdot & q & \cdot & \cdot & \cdot & \cdot & \cdot \cr ( 4,1^8)  & ( 
   \vr ,( 2,1^4) )  & {q^4} & {q^
    3} & \cdot & \cdot & \cdot & \cdot & \cdot & \cdot & \cdot & \cdot & \cdot & \cdot \cr ( 3^4)  & ( 
   ( 2^2) ,( 1^2) )  & {q^2} & \cdot & \cdot & {q^2} & {q^3} & {q^2} & \cdot & {q^
    3} & 2\,{q^2} & q & q & \cdot \cr ( 3^2,2^3)  & ( ( 2,1) ,
   ( 1^3) )  & {q^3} & \cdot & \cdot & q + {q^3} & \cdot & \cdot & \cdot & {q^4} & {q^
    3} & {q^2} & \cdot & \cdot \cr ( 3^2,2^2,1^2)  & ( ( 2) ,
   ( 1^4) )  & \cdot & \cdot & \cdot & {q^
    2} & \cdot & \cdot & \cdot & \cdot & \cdot & \cdot & \cdot & \cdot \cr ( 3^2,2,1^4)  & ( 
   ( 2,1^2) ,( 1^2) )  & {q^4} & \cdot & {q^3} & 2\,{q^2} + {q^4} & {q^
    3} & {q^2} & q & q + {q^5} & {q^2} + {q^4} & {q^3} & q & \cdot \cr ( 3^2,1^6)  & ( ( 2,1^3) ,( 1) )  & {q^3} & {q^2} & {q^4} & 3\,
   {q^3} & {q^4} & {q^3} & 2\,{q^2} & {q^2} & q + {q^3} & {q^2} & {q^
    2} & q \cr ( 3,2^4,1)  & ( ( 1^2) ,
   ( 1^4) )  & \cdot & \cdot & \cdot & {q^
    3} & \cdot & \cdot & \cdot & \cdot & \cdot & \cdot & \cdot & \cdot \cr ( 3,2^2,1^5)  & ( 
   ( 1^4) ,( 1^2) )  & \cdot & \cdot & {q^5} & 2\,{q^4} & \cdot & \cdot & {q^
    3} & {q^3} & {q^2} & \cdot & \cdot & \cdot \cr ( 3,1^9)  & ( \vr ,
   ( 1^6) )  & \cdot & \cdot & \cdot & \cdot & \cdot & \cdot & \cdot & \cdot & \cdot & \cdot & \cdot & \cdot \
\cr ( 2^5,1^2)  & ( ( 2^3) ,\vr )  & \cdot & \cdot & \cdot & {q^3} & {q^
    4} & q + {q^3} & {q^2} & {q^2} & {q^3} & \cdot & {q^2} & \cdot \cr ( 2^4,1^4)  & ( ( 2^2,1^2) ,\vr )  & {q^4} & {q^3} & \cdot & {q^4} & {q^
    5} & {q^2} + {q^4} & {q^3} & {q^3} & {q^2} + {q^4} & {q^3} & q + 
   {q^3} & {q^2} \cr ( 2,1^{10})  & ( ( 2,1^4) ,
   \vr )  & {q^5} & {q^4} & \cdot & {q^3} & \cdot & \cdot & {q^2} & \cdot & {q^3} & {q^
    4} & {q^2} & {q^3} \cr ( 1^{12})  & ( ( 1^6) ,\vr )  & \cdot & \cdot & \cdot & {q^4} & \cdot & \cdot & {q^
    3} & \cdot & \cdot & \cdot & \cdot & \cdot  
\end{array}
\end{align*}
%\noindent$\bullet$ ${\mathrm {wt}}=2\Lambda_0 - 3\alpha_0-3\alpha_1$ (continued).  
\begin{align*}&\text{ $\wt =2\Lambda_0 - 3\alpha_0-3\alpha_1$ (continued)} \\
&\begin{array}{l|l| c c c c c c c c c c c c}
( 4^3)  & ( ( 2) ,
   ( 2^2) )  & 1 & \cdot & \cdot & \cdot & \cdot & \cdot & \cdot & \cdot & \cdot & \cdot & \cdot & \cdot \cr ( 4^2,3,1)  & ( ( 1^2) ,
   ( 2^2) )  & q & 1 & \cdot & \cdot & \cdot & \cdot & \cdot & \cdot & \cdot & \cdot & \cdot & \cdot \cr ( 4^2,2,1^2)  & ( \vr ,
   ( 2^3) )  & q & \cdot & 1 & \cdot & \cdot & \cdot & \cdot & \cdot & \cdot & \cdot & \cdot & \cdot \cr ( 4^2,1^4)  & ( \vr ,( 2^2,1^2) )  & {q^
    2} & q & q & 1 & \cdot & \cdot & \cdot & \cdot & \cdot & \cdot & \cdot & \cdot \cr ( 4,3^2,2)  & ( 
   ( 2^2) ,( 2) )  & q & \cdot & \cdot & \cdot & 1 & \cdot & \cdot & \cdot & \cdot & \cdot & \cdot & \cdot \
\cr ( 4,3,2^2,1)  & ( ( 2) ,
   ( 2,1^2) )  & q & \cdot & \cdot & \cdot & \cdot & 1 & \cdot & \cdot & \cdot & \cdot & \cdot & \cdot \cr ( 4,3,1^5)  & ( ( 2,1^2) ,( 2) )  & {q^
    2} & \cdot & \cdot & \cdot & q & q & 1 & \cdot & \cdot & \cdot & \cdot & \cdot \cr ( 4,2^4)  & ( 
   ( 1^2) ,( 2,1^2) )  & 2\,
   {q^2} & q & q & \cdot & q & q & \cdot & 1 & \cdot & \cdot & \cdot & \cdot \cr ( 4,2^3,1^2)  & ( ( 1) ,( 2,1^3) )  & {q^3} & {q^2} & {q^2} & q & \cdot & {q^
    2} & \cdot & q & 1 & \cdot & \cdot & \cdot \cr ( 4,2^2,1^4)  & ( ( 1^3) ,
   ( 2,1) )  & {q^3} & {q^2} & \cdot & \cdot & {q^2} & {q^
    2} & q & q & \cdot & 1 & \cdot & \cdot \cr ( 4,2,1^6)  & ( ( 1^4) ,
   ( 2) )  & \cdot & \cdot & \cdot & \cdot & q & {q^3} & {q^2} & {q^
    2} & q & q & 1 & \cdot \cr ( 4,1^8)  & ( \vr ,
   ( 2,1^4) )  & \cdot & {q^3} & q & {q^2} & \cdot & \cdot & \cdot & {q^
    2} & q & q & \cdot & 1 \cr ( 3^4)  & ( ( 2^2) ,( 1^2) )  & {q^
    2} & q & \cdot & \cdot & q & \cdot & \cdot & \cdot & \cdot & \cdot & \cdot & \cdot \cr ( 3^2,2^3)  & ( 
   ( 2,1) ,( 1^3) )  & {q^3} & {q^2} & \cdot & \cdot & {q^2} & {q^
    2} & \cdot & q & \cdot & \cdot & \cdot & \cdot \cr ( 3^2,2^2,1^2)  & ( ( 2) ,
   ( 1^4) )  & \cdot & \cdot & \cdot & \cdot & \cdot & {q^3} & \cdot & {q^
    2} & q & \cdot & \cdot & \cdot \cr ( 3^2,2,1^4)  & ( ( 2,1^2) ,
   ( 1^2) )  & {q^2} + {q^4} & {q^3} & q & \cdot & q + {q^3} & {q^3} & {q^
    2} & {q^2} & \cdot & q & \cdot & \cdot \cr ( 3^2,1^6)  & ( ( 2,1^3) ,
   ( 1) )  & {q^3} & {q^2} & {q^2} & q & 2\,{q^2} & {q^4} & {q^3} & {q^
    3} & {q^2} & {q^2} & q & \cdot \cr ( 3,2^4,1)  & ( ( 1^2) ,
   ( 1^4) )  & {q^3} & {q^2} & {q^2} & q & {q^2} & {q^4} & \cdot & q + 
   {q^3} & {q^2} & \cdot & \cdot & \cdot \cr ( 3,2^2,1^5)  & ( ( 1^4) ,
   ( 1^2) )  & {q^4} & {q^3} & {q^3} & {q^2} & 2\,{q^3} & {q^5} & {q^
    4} & {q^2} + {q^4} & {q^3} & q + {q^3} & {q^2} & \cdot \cr ( 3,1^9)  & ( \vr ,( 1^6) )  & \cdot & \cdot & {q^
    2} & q & \cdot & \cdot & \cdot & {q^3} & {q^2} & {q^2} & \cdot & q \cr ( 2^5,1^2)  & ( ( 2^3) ,\vr )  & {q^3} & \cdot & {q^2} & \cdot & {q^
    2} & \cdot & \cdot & \cdot & \cdot & \cdot & \cdot & \cdot \cr ( 2^4,1^4)  & ( 
   ( 2^2,1^2) ,\vr )  & {q^4} & {q^3} & {q^3} & {q^2} & {q^
    3} & \cdot & \cdot & \cdot & \cdot & \cdot & \cdot & \cdot \cr ( 2,1^{10})  & ( 
   ( 2,1^4) ,\vr )  & \cdot & {q^4} & {q^2} & {q^3} & {q^
    2} & \cdot & \cdot & {q^3} & {q^2} & {q^2} & q & q \cr ( 1^{12})  & ( ( 1^6) ,\vr )  & \cdot & \cdot & {q^3} & {q^2} & {q^
    3} & \cdot & \cdot & {q^4} & {q^3} & {q^3} & {q^2} & {q^2}  
\end{array}
\end{align*}
%\noindent$\bullet$ ${\mathrm {wt}}=2\Lambda_0 - 3\alpha_0-3\alpha_1$ (continued). 
\begin{align*}
&\text{ $\wt =2\Lambda_0 - 3\alpha_0-3\alpha_1$ (continued)} \\
&\begin{array}{l|l| c c c c c c c c c c c c}
( 3^4)  & ( ( 2^2) ,
   ( 1^2) )  & 1 & \cdot & \cdot & \cdot & \cdot & \cdot & \cdot & \cdot & \cdot & \cdot & \cdot & \cdot \cr ( 3^2,2^3)  & ( ( 2,1) ,
   ( 1^3) )  & q & 1 & \cdot & \cdot & \cdot & \cdot & \cdot & \cdot & \cdot & \cdot & \cdot & \cdot \cr ( 3^2,2^2,1^2)  & ( ( 2) ,
   ( 1^4) )  & \cdot & q & 1 & \cdot & \cdot & \cdot & \cdot & \cdot & \cdot & \cdot & \cdot & \cdot \cr ( 3^2,2,1^4)  & ( ( 2,1^2) ,( 1^2) )  & {q^
    2} & q & \cdot & 1 & \cdot & \cdot & \cdot & \cdot & \cdot & \cdot & \cdot & \cdot \cr ( 3^2,1^6)  & ( ( 2,1^3) ,( 1) )  & q & {q^
    2} & q & q & 1 & \cdot & \cdot & \cdot & \cdot & \cdot & \cdot & \cdot \cr ( 3,2^4,1)  & ( 
   ( 1^2) ,( 1^4) )  & q & {q^
    2} & q & \cdot & \cdot & 1 & \cdot & \cdot & \cdot & \cdot & \cdot & \cdot \cr ( 3,2^2,1^5)  & ( ( 1^4) ,( 1^2) )  & {q^2} & {q^3} & {q^2} & {q^
    2} & q & q & 1 & \cdot & \cdot & \cdot & \cdot & \cdot \cr ( 3,1^9)  & ( 
   \vr ,( 1^6) )  & \cdot & \cdot & q & \cdot & \cdot & {q^
    2} & q & 1 & \cdot & \cdot & \cdot & \cdot \cr ( 2^5,1^2)  & ( ( 2^3) ,
   \vr )  & \cdot & \cdot & \cdot & q & \cdot & \cdot & \cdot & \cdot & 1 & \cdot & \cdot & \cdot \cr ( 2^4,1^4)  & ( ( 2^2,1^2) ,\vr )  & {q^2} & \cdot & \cdot & {q^
    2} & q & \cdot & \cdot & \cdot & q & 1 & \cdot & \cdot \cr ( 2,1^{10})  & ( 
   ( 2,1^4) ,\vr )  & {q^3} & {q^2} & q & q & {q^
    2} & \cdot & \cdot & \cdot & \cdot & q & 1 & \cdot \cr ( 1^{12})  & ( 
   ( 1^6) ,\vr )  & {q^2} & {q^3} & 2\,{q^2} & {q^2} & q & {q^
    3} & {q^2} & q & \cdot & \cdot & q & 1  
\end{array}
\end{align*}

%\end{document}

\bigskip 

\noindent {\large {\bf Acknowledgments}}  I would like to thank S. Ariki, T. Baker, M. Kashiwara, B. Leclerc, T. Miwa and J.-Y. Thibon for illuminating discussions. I am especially indebted to B. Leclerc for explanations concerning the papers \cite{Fo}, \cite{LT1}, \cite{LT2}. The influence of \cite{LT1} and \cite{LT2} on the present article will be obvious.

%------------------------------------------------------------------------------
% References
%------------------------------------------------------------------------------
\newcommand{\BOOK}[6]{\bibitem[{#6}]{#1}{\sc #2}, {\it #3} (#4) #5.}
\newcommand{\JPAPER}[8]{\bibitem[{#8}]{#1}{\sc #2}, {\it #3},
{#4} {\bf #5} (#6) #7.}
\newcommand{\PRE}[6]{\bibitem[{#6}]{#1}{\sc #2}, {\it #3},
{#4}  (#5).}
\newcommand{\JPAPERS}[9]{\bibitem[{\bf #9}]{#1}{\sc #2}, `#3', {\it #4} #5 #6,
#7 #8.}

\end{document}